\documentclass[10pt,a4paper,reqno]{amsart}
\usepackage{acronym}
\usepackage{amsfonts}
\usepackage{amsmath}
\usepackage{amssymb}
\usepackage{amsthm}
\usepackage{bm}
\usepackage{braket}
\usepackage{cite}
\usepackage{color}
\usepackage{geometry}
\usepackage{graphicx}
\usepackage{hyperref}
\usepackage{mathrsfs}
\usepackage{mathtools}
\usepackage{setspace}
\usepackage{epstopdf}

\theoremstyle{definition}

\theoremstyle{definition}

\usepackage{booktabs}
\usepackage{threeparttable}
\usepackage{diagbox}
\usepackage{multirow}
\usepackage{float}
\usepackage{array}

\usepackage{stmaryrd}
\usepackage{subcaption}
\usepackage{tikz}

\usepackage{epsfig}
\usepackage{setspace}
\usepackage{epstopdf}
\usepackage{booktabs}
\usepackage{threeparttable}
\usepackage{diagbox}

\newgeometry{top=2cm,bottom=2cm,outer=1.5cm,inner=1.5cm}
\newcommand{\bse}{\begin{subequations}}
\newcommand{\ese}{\end{subequations}}

\numberwithin{equation}{section}
\DeclareSymbolFont{largesymbols}{OMX}{yhex}{m}{n}
\DeclareMathAccent{\Widehat}{\mathord}{largesymbols}{"62}

\title[Gradient-enhanced physics-informed neural networks based on transfer learning for inverse problems of the variable coefficient differential equations]{Gradient-enhanced physics-informed neural networks based on transfer learning for inverse problems of the variable coefficient differential equations}
\author{Shuning Lin}
\address[SL]{School of Mathematical Sciences, Shanghai Key Laboratory of Pure Mathematics and Mathematical Practice, and Shanghai Key Laboratory of Trustworthy Computing \\
East China Normal University \\ Shanghai 200241 \\ China}
\author{Yong Chen$^*$}
\address[YC]{School of Mathematical Sciences, Shanghai Key Laboratory of Pure Mathematics and Mathematical Practice, and Shanghai Key Laboratory of Trustworthy Computing \\
East China Normal University \\ Shanghai 200241 \\ China}
\address[YC]{College of Mathematics and Systems Science \\ Shandong University of Science and Technology \\ Qingdao 266590 \\ China}
\email{ychen@sei.ecnu.edu.cn}

\begin{document}

\begin{abstract}

We propose gradient-enhanced PINNs based on transfer learning (TL-gPINNs) for inverse problems of the function coefficient discovery in order to overcome deficiency of the discrete characterization of the PDE loss in neural networks and improve accuracy of function feature description, which offers a new angle of view for gPINNs. The TL-gPINN algorithm is applied to infer the unknown variable coefficients of various forms (the polynomial, trigonometric function, hyperbolic function and fractional polynomial) and multiple variable coefficients simultaneously with abundant soliton solutions for the well-known variable coefficient nonlinear Schr\"{o}odinger equation. Compared with the PINN and gPINN, TL-gPINN yields considerable improvement in accuracy. Moreover, our method leverages the advantage of the transfer learning technique, which can help to mitigate the problem of inefficiency caused by extra loss terms of the gradient. Numerical results fully demonstrate the effectiveness of the TL-gPINN method in  significant accuracy enhancement, and it also outperforms gPINN in efficiency even when the training data was corrupted with different levels of noise or hyper-parameters of neural networks are arbitrarily changed.

\noindent{Keywords: TL-gPINN; transfer learning; variable coefficients; inverse problem.}

\end{abstract}
\maketitle

\section{Introduction}

With the vigorous development of nonlinear science, nonlinear models have been applied in more and more fields \cite{Peregrine1983,Gustafsson2015,Mio1976,Dalfovo1999,Hasegawa2013}, such as optical fiber communication, fluid mechanics, biophysics and information science. Among them, nonlinear evolution equations, especially the variable coefficient ones,  are an important class of nonlinear models and have attracted widespread attention \cite{TianGao2005,TianGao2006,TianShan2005} since the model with variable-coefficients is often preferable and suitable in describing real phenomena in many physical and engineering situations. For example, variable coefficient nonlinear Schr\"{o}odinger models plays an important role in the study of optical fiber system or the Rossby waves \cite{Rossby}. Variable coefficient higher-order nonlinear Schr\"{o}dinger and the variable coefficient Hirota equation can be used to describe the femtosecond pulse propagation \cite{vcHNLS} and the certain ultrashort optical pulses propagating in a nonlinear inhomogeneous fiber \cite{vcHirota1985}, respectively. Besides, many classical methods in the field of integrable systems have been widely used to derive exact solutions of variable coefficient equations, e.g., auto-B\"{a}cklund transformation \cite{TianGaoZhu2007,FanvcKdV}, the Riemann–Hilbert method \cite{ZhouvcHirota}, the Hirota bilinear method \cite{vcSG2014,vcHirota2017,DCQ2022}, the Darboux transformation \cite{DCQ2006,vcHirota2021,MvcNLS2022}, etc.

As early as the 1990s, the idea of solving partial differential equations (PDEs) by using the technique of neural networks was put forward \cite{Dissanayake1994}. However, limited by the level of science and technology at that time, it failed to get further development. With the explosive growth of computing resources, there has been renewed interest in researches of the numerical methods based on neural networks in recent years. This idea was revived by Raissi, Perdikaris and Karniadakis \cite{PINN2019} in 2019, and the physics-informed neural network (PINN) method was proposed to solve forward and inverse problems involving nonlinear partial differential equations. Based on the universal approximation theorems of neural networks \cite{Hornik1989}, it can accurately approximate functions with extraordinarily less data by embodying underlying physical constraints into neural networks. Due to its high accuracy and efficiency, the PINN method opens up a new approach for numerically solving nonlinear PDEs and immediately sets off a new research upsurge. On this foundation, variants and extensions targeted at different application scenarios also subsequently emerged in multitude, like fPINN \cite{fPINN2019} for solving fractional PDEs, NN-arbitrary polynomial chaos (NN-aPC) for solving stochastic problems \cite{stochastic2019}, XPINN \cite{XPINNs} and FBPINN \cite{FBPINN} for multiscale problems, B-PINN \cite{B-PINN} for forward and inverse PDE problems with noisy data, hp-VPINN \cite{hp-VPINNs} for rough solutions/input data such as singularities, steep solution, and sharp changes, etc. In addition, there have been many attempts to improve accuracy of the PINN method, such as locally adaptive activation functions with slope recovery \cite{Jagtap2020}, residual-based adaptive sampling \cite{DeepXDE2021,RAR2023}, gradient-enhanced PINN (gPINN)\cite{gPINN}, PINN with multi-scale Fourier features \cite{Fourierfeatures} and so on. Overall, the framework of PINNs, a model integrating data and mathematical physics seamlessly, is groundbreaking and has had a significant impact on the field of scientific computing and beyond.

Integrable deep learning, a concept first brought forward by Chen, deals with the combination of deep neural networks with integrable systems. In 2020, Li and Chen \cite{LJ1,LJ2} pioneered the use of the PINNs method in the field of integrable systems. Later, the dynamic behavior of rogue wave solution for the nonlinear Schr\"{o}dinger equation \cite{rouge} was reproduced for the first time by PINN. Abundant localized wave solutions are also simulated including the rogue periodic wave solution for the Chen–Lee–Liu equation \cite{periodicrouge}, vector localized waves for Manakov system \cite{PJCchaos}, data-driven soliton solutions for the nonlocal Hirota equation \cite{PWQnonlocal} and so on. Then the framework of the PINN method is extended to solve the ($N+1$)-dimensional initial-boundary value problem with $2N + 1$ hyperplane boundaries and is applied to reproduce dynamic behaviors (e.g., breathers, lump and resonance rogue) of high-dimensional integrable systems \cite{MZW2022}. Since integrable systems possess outstanding properties such as abundant symmetry, infinite conservation laws, the Lax pair and transformations, we devote to further improving the neural network method with the advantages of integrable systems. In 2022, we proposed new training algorithms based on the theory of integrable systems, namely a two-stage PINN method based on conserved quantities \cite{LSNJCP}. The novelty of this study lies in that constraints of conserved quantities, one of the most important features of integrable systems, have been successfully incorporated into neural networks to remarkably improve prediction accuracy and enhance the ability of generalization compared to the original PINN method. An implementation method of unsupervised learning——the PINN method based on Miura transformations \cite{LSNPD},  which is a significant transformation of integrable systems, was also put forward to solve nonlinear PDEs. With the aid of this method, we can simply exploit the initial–boundary data of a solution for a certain nonlinear equation to obtain the data-driven solution for another evolution equation. It was applied to discover a new type of numerical solution, i.e., kink-bell type solution of the defocusing mKdV equation, by fully leveraging the many-to-one relationship between solutions before and after Miura transformations. Collectively, how to devise significant integrable-deep learning algorithms and utilize the PINN method to pertinently solve problems arising in the field of integrable systems that cannot be solved by classical methods is our target to aim at.

Despite some progress, solving inverse problems by traditional numerical methods still requires complex mathematics deduction and extensive calculations. Comparatively, deep learning algorithm has great advantages in solving inverse problems of partial differential equations. Inverse problem refers to the task that, given some information of the solution, it is expected to deduce the unknown quantity (which may be the unknown coefficient or unknown term) in PDE and the solution itself. Traditional numerical methods have many limitations in solving inverse problems, especially in dealing with noisy data, complex regions, and high-dimensional problems. However, the PINN method is mesh-free and has been proved to be robust to noise in many cases. Even for high-dimensional tasks, it also shows the outstanding performance in both accuracy and computing efficiency. What's more, physics-informed machine learning has better interpretability and can achieve satisfactory accuracy and better convergence with a small amount of data by embedding physics into machine learning.

Most past researches on inverse problems were concentrated on the parameter discovery (such as the inverse problem of constant coefficient equation) rather than the function discovery, so this paper mainly focused on studying the inverse problem of variable coefficient PDEs, an important class of equations in integrable systems, by using deep learning algorithm. The classical methods for studying integrable systems, e.g., the aforementioned Hirota bilinear method, the Darboux transformation and Riemann-Hilbert method, can only be utilized to derive the solution of the variable coefficient equation, while the PINN method can obtain not only the solution itself but also the corresponding unknown variable coefficients. Meanwhile, previous studies on variable coefficient equations using PINN method are also relatively few \cite{zhjvcHirota,MZWvcPINN}, especially on the improvement of the algorithm for accuracy enhancement.

To fill in the gap, the gradient-enhanced PINN method is considered here. The original gPINN stems from that PINNs only enforce the PDE residual $f$ to be $0$ and one can further utilize the property that when a function equals to $0$, its derivatives of all orders are also $0$, while our motivation and conception of this method are entirely different. We aim to overcome deficiency of the discrete characterization of the PDE loss in neural networks, and to improve accuracy of function feature description, which offers a new angle of view for gPINNs. Specifically, the PDE constraint is characterized by some discrete points in PINNs. Since it is impossible to consider all collocation points, the term of PDE loss can only ensure that values of the variable coefficient are close to the true ones at these selected points, while the accuracy at points outside the given points lacks adequate attention. Moreover, it is biased to characterize a function (i.e., the unknown variable coefficient here) solely by the values at discrete points. Thus, loss terms of the gradients are introduced to enhance the accuracy of the identified variable coefficient from the perspective of gradients. Due to the lack of partial derivative values of variable coefficient at configuration points, a straightforward and simple way is to enforce the partial derivatives of it to satisfy the corresponding equation. However, we are inevitably faced with the challenge of slow computation speed caused by extra constraints of the gradient. One viable path towards accelerating the convergence of training could come by adopting the technique of transfer learning. Therefore, the gradient-enhanced PINN method based on transfer learning (TL-gPINNs) is brought forward in this thesis.

This paper is organized as follows. In Sec. \ref{Methodology}, we propose gradient-enhanced PINNs based on transfer learning for inverse problems of the variable coefficient equations after a brief review of the PINN and gradient-enhanced PINN (gPINN) methods. Then in Sec. \ref{example}, the TL-gPINN method is applied to identify the unknown variable coefficients of various forms (e.g., the linear, quadratic, sine, hyperbolic tangent and fractional forms) together with the soliton solution for the variable coefficient nonlinear Schr\"{o}odinger equation. We also systematically compare the performance of PINNs, gPINNs and TL-gPINNs, and numerical results demonstrate the ability of TL-gPINNs in significant accuracy enhancement. Further error analyses including the robustness analysis and parametric sensitivity analysis of TL-gPINNs are conducted in Sec. \ref{analysis}, where the heat map serves as an effective visualization tool. Finally, the conclusion and expectation are given in the last section.

\section{Methodology}\label{Methodology}
\subsection{Introduction of PINNs and gradient-enhanced PINNs}
\quad

The first part gives a brief overview of physics-informed neural networks (PINNs), an effective tool in solving forward and inverse problems of partial differential equations.

Let's consider the general form of a ($N+1$)-dimensional partial differential equation with parameters $\boldsymbol{\lambda}$
\begin{align}
f\left(\mathbf{x},t ; \frac{\partial u}{\partial x_1}, \ldots, \frac{\partial u}{\partial x_N}, \frac{\partial u}{\partial t} ; \frac{\partial^2 u}{\partial x_1^2}, \ldots, \frac{\partial^2 u}{\partial x_1 \partial x_N}, \frac{\partial^2 u}{\partial x_1 \partial t} ; \ldots ; \boldsymbol{\lambda}\right)=0, \quad \mathbf{x}=\left(x_1, \cdots, x_N\right) \in \Omega, \quad t \in [t_0, t_1],
\end{align}
where $u(\mathbf{x},t)$ is the solution and $\Omega$ is a subset of $\mathbb{R}^N$. 

To solve the above PDE with the first kind of boundary condition (Dirichlet boundary condition)
\begin{equation}
\begin{split}
\begin{cases}
u(\mathbf{x},t_0)=u_0(\mathbf{x}), \quad \forall \mathbf{x} \in \Omega\\
u(\mathbf{x},t)=\mathcal{B}(\mathbf{x},t), \quad \forall \mathbf{x} \in \partial \Omega, t \in [t_0, t_1],
\end{cases}
\end{split}
\end{equation}
we construct a neural network of depth $L$ consisting of one input layer, $L-1$ hidden layers and one output layer. Suppose that the $l$th ($l=0,1,\cdots,L$) layer has $N_l$ neurons, and then the connection between layers can be achieved by the following affine transformation $\mathcal{A}$ and activation function $\sigma(\cdot)$
\begin{align}
\mathbf{x}^l=\sigma(\mathcal{A}_l(\mathbf{x}^{l-1}))=\sigma(\mathbf{w}^{l} \mathbf{x}^{l-1}+\mathbf{b}^{l}),	
\end{align}
where $\mathbf{w}^{l}\in \mathbb{R}^{N_{l} \times N_{l-1}}$ and $\mathbf{b}^{l}\in \mathbb{R}^{N_{l}}$ denote the weight matrix and bias vector separately. Especially, the input is $\mathbf{x}^0=\left(x_1, \cdots, x_N,t\right)$ and output $\mathbf{o}(\mathbf{x}^0,\boldsymbol{\Theta})$ is given by
\begin{align}
\mathbf{o}(\mathbf{x}^0,\boldsymbol{\Theta})=(\mathcal{A}_L \circ \sigma \circ \mathcal{A}_{L-1} \circ \cdots \circ \sigma \circ \mathcal{A}_1)(\mathbf{x}^0),
\end{align}
which is used to approximate the solution $u(\mathbf{x},t)$, and $\boldsymbol{\Theta}=\left\{\mathbf{w}^{l}, \mathbf{b}^{l}\right\}_{l=1}^{L}$ represents the trainable parameters of PINN. With the initial-boundary dataset $\{\mathbf{x}^i_u,t^i_u,u^i\}^{N_u}_{i=1}$ and the set of collocation points of $f(\mathbf{x},t)$, denoted by $\{\mathbf{x}_{f}^i,t_{f}^i\}^{N_{f}}_{i=1}$, the loss function can be defined to measure the difference between the predicted values and the true values of each iteration
\begin{equation}
MSE_{forward}=MSE_u+MSE_f,
\end{equation}
where
\begin{equation}
MSE_u=\frac{1}{N_u}\sum^{N_u}_{i=1}|\Widehat{u}(\mathbf{x}_u^i,t_u^i)-u^i|^2,
\end{equation}
\begin{equation}
MSE_{f}=\frac{1}{N_f}\sum^{N_f}_{i=1}|f(\mathbf{x}_{f}^i,t_{f}^i)|^2.
\end{equation}

With regard to the inverse problem, namely, the situation that the parameters $\boldsymbol{\lambda}$ are unknown, some extra measurements $\{\mathbf{x}^i_{in},t^i_{in},u^i_{in}\}^{N_{u_{in}}}_{i=1}$ of the internal area should be obtained and utilized to define a new loss function to learn the unknown parameters $\boldsymbol{\lambda}$
\begin{equation}
MSE_{inverse}=MSE_u+MSE_f+MSE_{u_{in}},
\end{equation}
where
\begin{equation}
MSE_{u_{in}}=\frac{1}{N_{u_{in}}}\sum^{N_{u_{in}}}_{i=1}|\Widehat{u}(\mathbf{x}_{in}^i,t_{in}^i)-u^i_{in}|^2.
\end{equation}

Later, a deep learning method, gradient-enhanced physics-informed neural networks (gPINNs) \cite{gPINN} was proposed for improving the accuracy and training efficiency of PINNs by leveraging gradient information of the PDE residual and embedding the gradient into the loss function. The basic idea of gPINNs is that it enforces the derivatives of the PDE residual $f$ to be zero since $f(\mathbf{x},t)$ is zero for any $\mathbf{x}$ and $t$, i.e.,
\begin{equation}
\nabla f(\mathbf{x})=\left(\frac{\partial f}{\partial x_1}, \frac{\partial f}{\partial x_2}, \cdots, \frac{\partial f}{\partial x_N},\frac{\partial f}{\partial t}\right)=\mathbf{0}, \quad \mathbf{x} \in \Omega, \quad t \in [t_0, t_1].
\end{equation}
Then, based on the set of residual points $\{\mathbf{x}_{g}^i,t_{g}^i\}^{N_{g}}_{i=1}$ for the derivatives, the loss functions of the forward and inverse problems are separately defined as
\begin{align}
&MSE_{forward}^g=MSE_u+MSE_f+MSE_g,\\
&MSE_{inverse}^g=MSE_u+MSE_f+MSE_{u_{in}}+MSE_g,
\end{align}
where 
\begin{align}
MSE_g=\frac{1}{N_g} \left( \sum^{N}_{j=1} \sum^{N_g}_{i=1}|\frac{\partial f}{\partial x_j}(\mathbf{x}_{g}^i,t_{g}^i)|^2 + \sum^{N_g}_{i=1}|\frac{\partial f}{\partial t}(\mathbf{x}_{g}^i,t_{g}^i)|^2 \right).	
\end{align}
The set of residual points $\{\mathbf{x}_{g}^i,t_{g}^i\}^{N_{g}}_{i=1}$ for the derivatives can be different from the set of collocation points $\{\mathbf{x}_{f}^i,t_{f}^i\}^{N_{f}}_{i=1}$ of $f(\mathbf{x},t)$, but we usually choose the same set for convenience.

\subsection{Gradient-enhanced PINNs based on transfer learning for data-driven variable coefficients}
\quad

For the inverse PDE problems with the aid of PINNs and its variants, the existing researches are mainly focused on the parameter discovery rather than the function discovery. Here, we propose the gradient-enhanced PINNs based on transfer learning (TL-gPINNs) to infer variable coefficients.

$\bullet$ \textbf{Motivation}

For the inverse problem of identifying the variable coefficients, we aim to improve the neural network method to enhance prediction precision.

Consider the ($N+1$)-dimensional PDE with variable coefficients $\boldsymbol{\Lambda}(\mathbf{x},t)$
\begin{align}
f\left(\mathbf{x},t ; \frac{\partial u}{\partial x_1}, \ldots, \frac{\partial u}{\partial x_N}, \frac{\partial u}{\partial t} ; \frac{\partial^2 u}{\partial x_1^2}, \ldots, \frac{\partial^2 u}{\partial x_1 \partial x_N}, \frac{\partial^2 u}{\partial x_1 \partial t} ; \ldots ; \boldsymbol{\Lambda}\right)=0, \quad \mathbf{x}=\left(x_1, \cdots, x_N\right) \in \Omega, \quad t \in [t_0, t_1].
\end{align}
Obviously, neural networks are constrained by the PDE loss $MSE_{f}=\frac{1}{N_f}\sum^{N_f}_{i=1}|f(\mathbf{x}_{f}^i,t_{f}^i)|^2$ to satisfy the equation above. In other words, PINNs enforce the PDE residual $f$ to be $0$ to make the data-driven variable coefficients $\boldsymbol{\Lambda}(\mathbf{x},t)$ close to the exact ones.

Since this constraint is characterized by discrete points $\{\mathbf{x}_{f}^i,t_{f}^i\}^{N_{f}}_{i=1}$, it can only ensure that values of the obtained $\boldsymbol{\Lambda}(\mathbf{x},t)$ are close to the true ones at these selected points. However, even if the predicted values of the variable coefficients $\boldsymbol{\Lambda}(\mathbf{x},t)$ are equal to the true values at these discrete points, the values of the identified $\boldsymbol{\Lambda}(\mathbf{x},t)$ outside the given point set $\{\mathbf{x}_{f}^i,t_{f}^i\}^{N_{f}}_{i=1}$ have little direct effect on $MSE_{f}$ and thus the data-driven variable coefficients may be a considerable departure from the exact ones.

Consequently, it is biased and inaccurate to characterize a function (i.e., the unknown variable coefficient here) solely by the values at discrete points. It is necessary to introduce much more rigorous constraints considering the partial derivative values of variable coefficients from the perspective of enhancing gradients. Due to the lack of partial derivative values of variable coefficient at configuration points, a way with similar effect is to enforce the partial derivatives of it to satisfy the corresponding equation. Then residuals of equations satisfied by not only variable coefficients but also the partial derivatives of them should be taken into account in the design of loss functions. Specifically, neural networks enforce the residual of equations satisfied by both variable coefficients and the partial derivatives to be $0$ to require the values of them to approximate the true ones at the discrete points. It happens to coincide with gradient-enhanced PINNs, the idea of which stems from that derivatives of the zero-valued function, (i.e., the PDE residual $f$) should also be 0.

Note that the equations satisfied by the partial derivatives of variable coefficients ($\frac{\partial \boldsymbol{\Lambda}}{\partial t}, \frac{\partial \boldsymbol{\Lambda}}{\partial x_1}, \frac{\partial \boldsymbol{\Lambda}}{\partial x_2}, \cdots$) can be derived by direct differentiation with respect to the equation satisfied by variable coefficients $\boldsymbol{\Lambda}(\mathbf{x},t)$ themselves. Although our perspective and concerned aspects are different from that of gPINNs, the implementation method is the same.

Ulteriorly, we improve the original gPINNs by means of the idea of transfer learning since the introduction of additional constraints on gradients probably gives rise to low efficiency.

$\bullet$ \textbf{Procedure}

By fully leveraging the combined advantages of gradient-enhanced effect and transfer learning, TL-gPINN is proposed.

The followings are the main steps involved:

Firstly, the traditional PINNs are constructed to obtain the data-driven variable coefficients after defining the following loss function 
\begin{equation}
MSE_{inverse}=MSE_u+MSE_f+MSE_{u_{in}}+MSE_{\boldsymbol{\Lambda}},
\end{equation}
where 
\begin{equation}
MSE_u=\frac{1}{N_u}\sum^{N_u}_{i=1}|\Widehat{u}(\mathbf{x}_u^i,t_u^i)-u^i|^2,
\end{equation}
\begin{equation}
MSE_{f}=\frac{1}{N_f}\sum^{N_f}_{i=1}|f(\mathbf{x}_{f}^i,t_{f}^i)|^2,
\end{equation}
\begin{equation}
MSE_{u_{in}}=\frac{1}{N_{u_{in}}}\sum^{N_{u_{in}}}_{i=1}|\Widehat{u}(\mathbf{x}_{in}^i,t_{in}^i)-u^i_{in}|^2,
\end{equation}
\begin{equation}
MSE_{\boldsymbol{\Lambda}}=\frac{1}{N_{\boldsymbol{\Lambda}}}\sum^{N_{\boldsymbol{\Lambda}}}_{i=1}|\Widehat{{\boldsymbol{\Lambda}}}(\mathbf{x}_{\boldsymbol{\Lambda}}^i,t_{\boldsymbol{\Lambda}}^i)-{\boldsymbol{\Lambda}}^i|^2,
\end{equation}
and the set $\{\mathbf{x}^i_{\boldsymbol{\Lambda}},t^i_{\boldsymbol{\Lambda}},{\boldsymbol{\Lambda}}^i\}^{N_{\boldsymbol{\Lambda}}}_{i=1}$ denotes the boundary training data of the variable coefficients $\boldsymbol{\Lambda}(\mathbf{x},t)$. In particular, the networks of the solution $u(\mathbf{x},t)$ and the variable coefficients $\boldsymbol{\Lambda}(\mathbf{x},t)$ are set to be separate and named as the trunk and branch networks respectively in order to eliminate mutual influence. What's more, with regard to the network of variable coefficients $\boldsymbol{\Lambda}(\mathbf{x},t)$, the width and depth of it are usually narrower and shallower than that of the solution $u(\mathbf{x},t)$ since the function expression of variable coefficient is simpler in general.

Then at the end of the iteration process, the weight matrixes and bias vectors of PINNs are saved to initialize the gradient-enhanced PINNs with the advantage of transfer learning. The mean squared error loss function of gPINNs is given by
\begin{equation}
MSE_{inverse}^g=MSE_u+MSE_f+MSE_{u_{in}}+MSE_{\boldsymbol{\Lambda}}+MSE_{g},
\end{equation}
where
\begin{align}
MSE_g=\frac{1}{N_g} \left( \sum^{N}_{j=1} \sum^{N_g}_{i=1}|\frac{\partial f}{\partial x_j}(\mathbf{x}_{g}^i,t_{g}^i)|^2 + \sum^{N_g}_{i=1}|\frac{\partial f}{\partial t}(\mathbf{x}_{g}^i,t_{g}^i)|^2 \right).	
\end{align}
Based on MSE criteria, the parameters of gPINNs are optimized and we finally obtain the data-driven variable coefficients $\boldsymbol{\Lambda}(\mathbf{x},t)$.

\begin{figure}
\centering
\includegraphics[width=18cm,height=7.5cm]{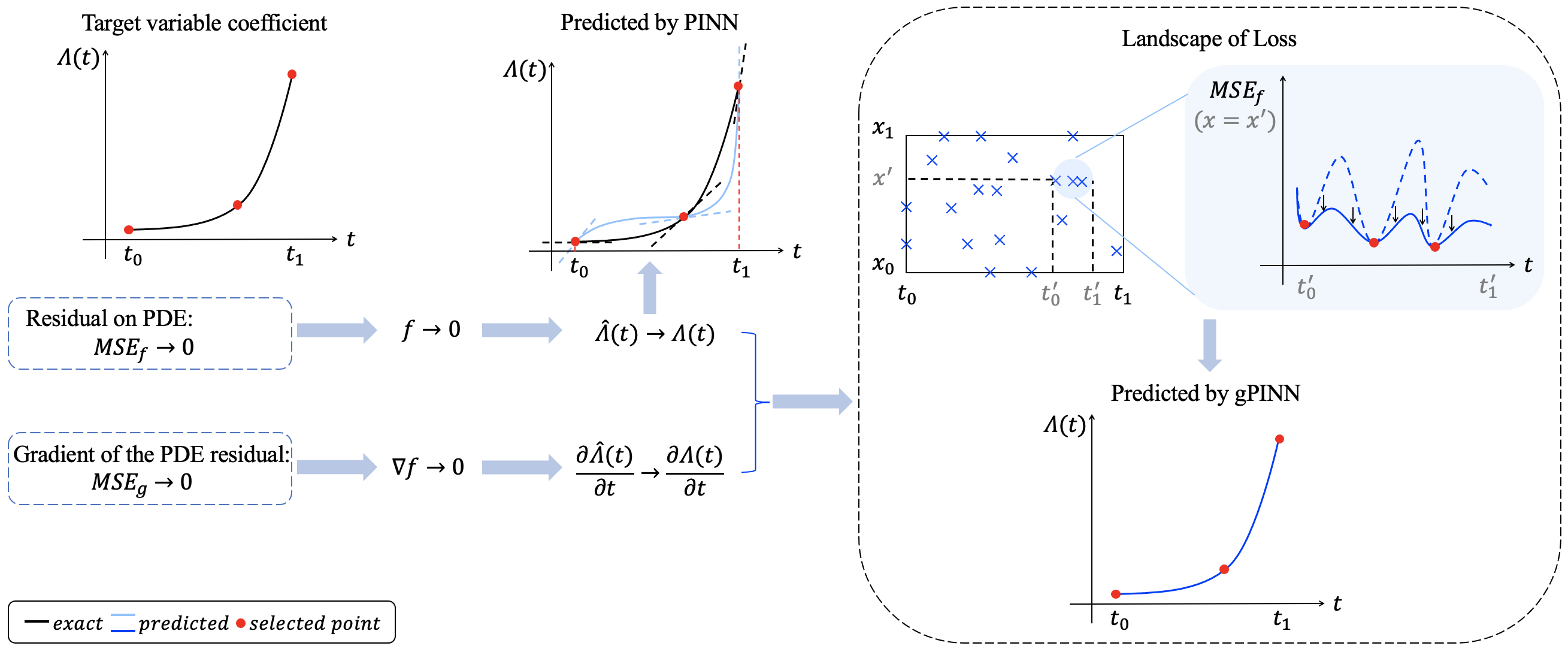}
\caption{(Color online) The effect of gPINN compared to PINN on the optimization of loss (PDE residual) and inference of time-varying variable coefficient $\Lambda(t)$.}
\label{fig2-1}
\end{figure}

To better understand our new angle of view for gPINNs, we take the inverse problem of identifying the time-varying variable coefficient $\Lambda(t)$ as an example to illustrate the effect of gradient information on the optimization of loss (PDE residual) and inference of time-varying variable coefficient. The corresponding sketch map is displayed in Fig. \ref{fig2-1}. Wherein, the value of the predicted variable coefficient will close to that of the exact one  at the selected point if the term of PDE residual $MSE_f$ is considered solely. The derivative $\frac{\partial \hat{\Lambda}(t)}{\partial t}$ of the predicted variable coefficient $\hat{\Lambda} (t)$ is constrained to approximate $\frac{\partial \Lambda(t)}{\partial t}$ at the given point due to the incorporation of the gradient term $MSE_g$, which reflects the equation information satisfied by the derivative of variable coefficient. Then the combined effect is remarkable and it enables the predicted $\hat{\Lambda} (t)$ to tend to the exact $\Lambda(t)$.

\begin{figure}
\centering
\includegraphics[width=18cm,height=14.5cm]{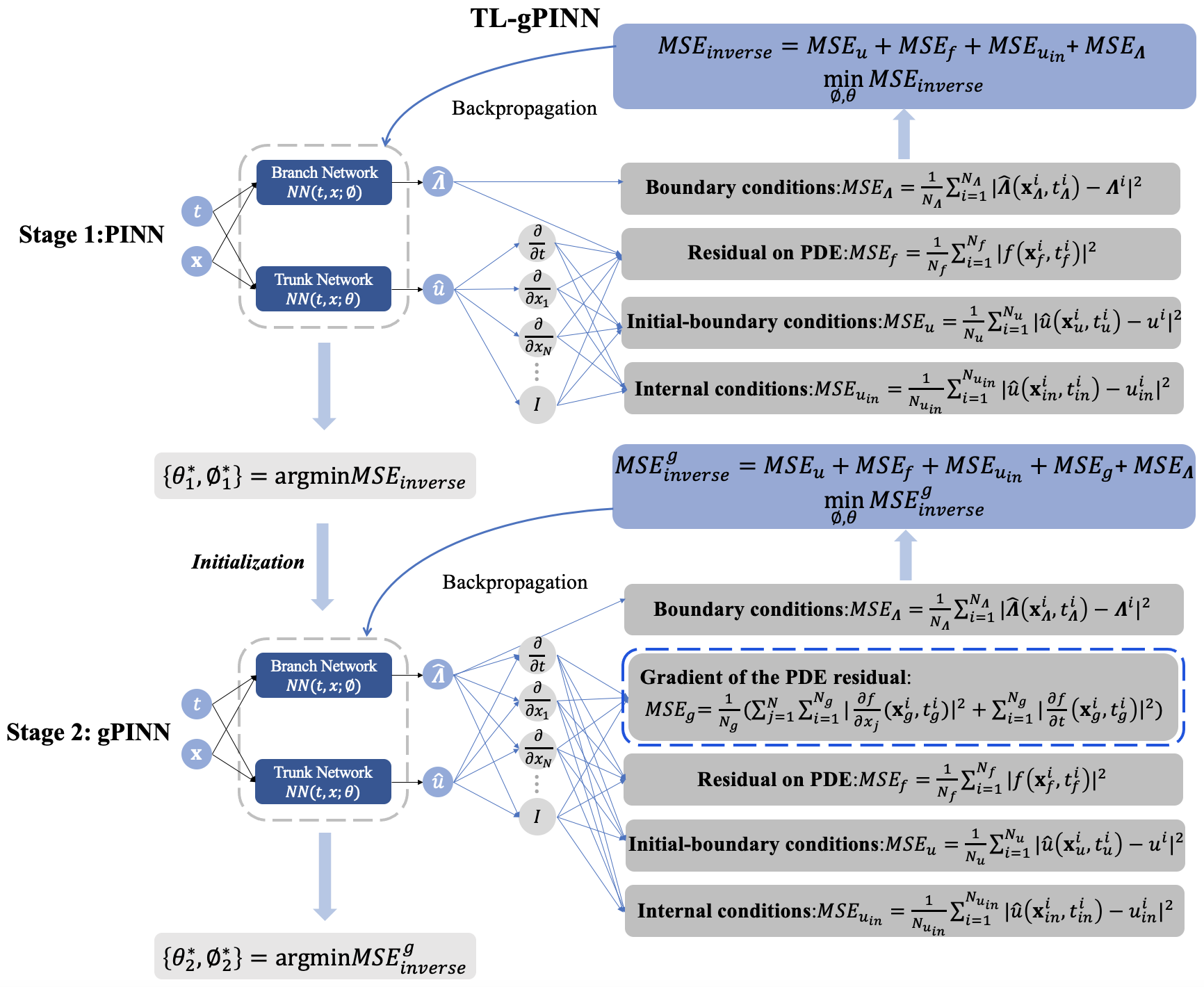}
\caption{(Color online) Schematic diagram of the TL-gPINN algorithm.}
\label{fig2-2}
\end{figure}

Our method uses a two-step optimization strategy and gradually increases the difficulty, resulting in better results than the direct one-step optimization, i.e., the original gPINN method. The process draft of the TL-gPINN method is sketched and shown in Fig. \ref{fig2-2}.

Partial derivatives of higher orders, of course, can be considered by adding the loss of $\frac{\partial^2 f}{\partial x_i \partial x_j}, \frac{\partial^2 f}{\partial x_i \partial t}$ and so on into the term $MSE_g$. However, excessive constraints may lead to high training costs and low efficiency, which is the reason why the transfer learning technique is introduced here.

The difficulty of the inverse problem lies in that the information about variable coefficients is insufficient and meanwhile, the properties or physical laws (if any) of variable coefficients remain to be discovered, which can be used to achieve higher accuracy. Therefore, we should make full use of existing information. For example, the equations satisfied by the partial derivatives of variable coefficients ($\frac{\partial \boldsymbol{\Lambda}}{\partial t}, \frac{\partial \boldsymbol{\Lambda}}{\partial x_1}, \frac{\partial \boldsymbol{\Lambda}}{\partial x_2}, \cdots$) can be derived by differentiating the equation satisfied by variable coefficients $\boldsymbol{\Lambda}(\mathbf{x},t)$ themselves. Thus, the gradient-enhanced PINN can be served as an effective tool to improve the accuracy of variable coefficients by fully utilizing the information of equations satisfied by the derivatives of variable coefficients. Further, TL-gPINNs can improve the accuracy and training efficiency of the original gPINNs.

All the codes in this article are based on Python 3.7 and Tensorflow 1.15, and the presented numerical experiments  are run on a DELL Precision 7920 Tower computer with 2.10 GHz 8-core Xeon Silver 4110 processor, 64 GB memory and 11 GB NVIDIA GeForce GTX 1080 Ti video card.

\section{Applications in the variable coefficient
nonlinear Schr\"{o}odinger equation}\label{example}

The nonlinear Schr\"{o}odinger (NLS) equation, one of the most classical equations in integrable systems, is commonly used in the field of optical fiber communication to describe the propagation of optical solitons \cite{Agrawal2000,Kivshar2003}. When it comes to inhomogeneous optical fibers, it is believed that the variable coefficient Schr\"{o}odinger equation is more accurate and realistic than the standard one since variable coefficients can reflect the inhomogeneities of media and nonuniformities of boundaries \cite{Gupta2018}. The research on variable coefficient NLS-type models has achieved very fruitful results \cite{Papaioannou1996,HaoPRE2004,Meng2007,Karpman2001,Ohta1994}, like the groundbreaking work of Serkin et al. \cite{Serkin2001}. Meanwhile, solutions of the variable coefficient NLS-type equations are also obtained by using powerful means, such as the Hirota bilinear method \cite{PRE2008,vcNLS}, the Darboux transformation \cite{Hao2004}, the Riemann–Hilbert approach \cite{RHvcNLS2021} and so on.

In this part, we discuss the mathematical model which can be used to describe the optical fiber system or the Rossby waves \cite{Rossby}, i.e., the variable coefficient nonlinear Schr\"{o}odinger (vcNLS) equation 
\begin{align}
{\rm{i}} A_t +\alpha(t) A_{xx}+\beta(t) A+\gamma(t) |A|^2 A=0,
\end{align}
where the variable coefficient $\alpha(t)$ denotes the dispersion effect and $\gamma(t)$ denotes the Kerr nonlinearity. When considering the inhomogeneities, the varying dispersion and Kerr nonlinearity are of practical importance in the optical-fiber transmission system.

Under the assumption that the amplitude $A(x,t)$ has the transformation
\begin{align}
A(x,t)=	{\rm{e}}^{{\rm{i}} \int \beta(t) dt} \frac{h(x,t)}{g(x,t)},
\end{align}
the one-soliton solution can be derived by the Hirota method \cite{vcNLS}
\begin{align}
A(x,t)=	{\rm{e}}^{{\rm{i}} \int \beta(t) dt} \frac{{\rm{e}}^{\theta}}{1+\frac{\gamma(t)}{2 \alpha(t) (k+k^{\ast})^2} {\rm{e}}^{\theta+ \theta^{\ast}}},	
\end{align}
where 
\begin{align}
&\phi(t)={\rm{i}} \int \alpha(t) k^2 dt,\\
&\theta = k x +\phi(t)+\eta,\\	
&\theta^{\ast} = k^{\ast} x +\phi(t)+\eta.
\end{align}
Here, $k$ is a complex constant and $\eta$ is a real constant.

The PINN, gPINN and TL-gPINN methods are applied to infer the unknown time-varying variable coefficients of the vcNLS equation. To avoid repetition, hyper-parameters of neural networks used for each case are listed in Table \ref{table3-0} and other parameters are selected as $k=1+{\rm{i}}, \eta=0$.

\begin{table}[htbp]
\caption{Hyper-parameters used for each case}
\label{table3-0}  
\centering
\begin{tabular}{c|cccccc}
\bottomrule
\multirow{2}{*}{Section} & \multicolumn{2}{c|}{Variable coefficients}                 & \multicolumn{2}{c|}{Trunk network}               & \multicolumn{2}{c}{Branch network} \\ \cline{2-7} 
  & \multicolumn{1}{c|}{Fixed} & \multicolumn{1}{c|}{Inferred} & \multicolumn{1}{c|}{Depth} & \multicolumn{1}{c|}{Width} & \multicolumn{1}{c|}{Depth}                                & Width                               \\ \hline
3.1.1                    &  $\alpha(t)=\frac{t}{2}, \beta(t) = \frac{t}{5}$                           &  $\gamma(t)=t$                           & 8                          & 40                         & 4                                                         & 30                                  \\
3.1.2                    &  $\alpha(t)=\frac{t^2}{2}, \beta(t) = \frac{t}{5}$                           &  $\gamma(t)=t^2$                           & 8                          & 40                         & 4                                                         & 30                                  \\
3.1.3                    &  $\alpha(t)=\sin(t),\beta(t)=\frac{t}{5}$                           &  $\gamma(t)=\sin(t)$                           & 8                          & 40                         & 4                                                         & 30                                  \\
3.1.4                    &  $\alpha(t)=\tanh(t), \beta(t) = \frac{t}{5}$                           &  $\gamma(t)=\tanh(t)$                           & 8                          & 40                         & 4                                                         & 30                                  \\
3.1.5                    &  $\alpha(t)=\frac{1}{2(1+t^2)},\beta(t)=\frac{t}{5}$                           &  $\gamma(t)=\frac{1}{1+t^2}$                           & 8                          & 40                         & 4                                                         & 30                                  \\
3.2.1                    &  $\alpha(t)=\sin(t)$                           &  $\beta(t)=\frac{t}{5},\gamma(t)=\sin(t)$                           & 8                          & 40                         & 4                                                         & 30                                  \\
3.2.2(1)                 &  -           & $\alpha(t)=\frac{t}{2},\beta(t)=\frac{t}{5},\gamma(t)=t$                              & 8                          & 40                         & 2                                                         & 10                                  \\
3.2.2(2)                 &  -           &  $\alpha(t)=\frac{1}{2(1+t^2)},\beta(t)=\frac{t}{5},\gamma(t)=\frac{1}{1+t^2}$                            & 8                          & 40                         & 4                                                         & 10                                  \\ \toprule
\end{tabular}
\end{table}

\subsection{Data-driven discovery of single variable coefficient}\label{NLS}
\quad

The general aim is to utilize TL-gPINNs to solve the inverse problem for the discovery of function coefficient $\gamma(t)$, and systematically compare the performance of three methods (PINNs, gPINNs and TL-gPINNs) under the circumstances that the other two variable coefficients $\alpha(t)$ and $\beta(t)$ are already known.

Several types of time-varying variable coefficients $\gamma(t)$ in common use are provided, such as linear, quadratic, sine, hyperbolic tangent and fractional functions: $\gamma(t)=t, t^2, \sin(t), \tanh(t), \frac{1}{1+t^2}$ separately.

\subsubsection{Data-driven discovery of  linear variable coefficient $\gamma(t)$}\label{3.1}
\quad

In this part, we take $\alpha(t)=\frac{t}{2}, \beta(t) = \frac{t}{5}$ and choose $[x_0,x_1]=[-4,4], [t_0,t_1]=[-4,4]$ as the training region. 

In consideration of the complexity of the structure of complex-valued solution $A(x,t)$, we decompose it into real part $u(x,t)$ and imaginary part $v(x,t)$, i.e., $A(x,t)=u(x,t)+{\rm{i}} v(x,t)$. After substituting it into the governing equation
\begin{equation}
f:={\rm{i}} A_t +\alpha(t) A_{xx}+\beta(t) A+\gamma(t) |A|^2 A=0,
\end{equation}
we have
\begin{equation}
f_u: =-v_t+\alpha(t) u_{xx} +\beta(t) u +\gamma(t) (u^2+v^2)u,
\end{equation}
\begin{equation}
f_v: =u_t+\alpha(t) v_{xx} +\beta(t) v +\gamma(t) (u^2+v^2)v.
\end{equation}
Define the loss function of PINNs for the inverse problem as follows:
\begin{equation}
MSE_{inverse}=MSE_A+MSE_f+MSE_{A_{in}}+MSE_{\gamma},
\end{equation}
where
\begin{equation}
MSE_A=MSE_u+MSE_v,\quad MSE_f=MSE_{f_u}+MSE_{f_v},\quad MSE_{A_{in}}=MSE_{u_{in}}+MSE_{v_{in}},	
\end{equation}
\begin{equation}
MSE_u=\frac{1}{N_A}\sum^{N_A}_{i=1}|\Widehat{u}(x_A^i,t_A^i)-u^i|^2,
\end{equation}
\begin{equation}
MSE_v=\frac{1}{N_A}\sum^{N_A}_{i=1}|\Widehat{v}(x_A^i,t_A^i)-v^i|^2,
\end{equation}
\begin{equation}
MSE_{f_u}=\frac{1}{N_f}\sum^{N_f}_{i=1}|f_u(x_{f}^i,t_{f}^i)|^2,
\end{equation}
\begin{equation}
MSE_{f_v}=\frac{1}{N_f}\sum^{N_f}_{i=1}|f_v(x_{f}^i,t_{f}^i)|^2,
\end{equation}
\begin{equation}
MSE_{u_{in}}=\frac{1}{N_{A_{in}}}\sum^{N_{A_{in}}}_{i=1}|\Widehat{u}(x_{in}^i,t_{in}^i)-u^i_{in}|^2,
\end{equation}
\begin{equation}
MSE_{v_{in}}=\frac{1}{N_{A_{in}}}\sum^{N_{A_{in}}}_{i=1}|\Widehat{v}(x_{in}^i,t_{in}^i)-v^i_{in}|^2,
\end{equation}
\begin{equation}
MSE_{\gamma}=\frac{1}{N_{\gamma}}\sum^{N_{\gamma}}_{i=1}|\Widehat{{\gamma}}(t_{\gamma}^i)-{\gamma}^i|^2=|\Widehat{{\gamma}}(t_0)-{\gamma}^0|^2.
\end{equation}
Here, $\{x^i_A,t^i_A,u^i,v^i\}^{N_A}_{i=1}$ and $\{x^i_{in},t^i_{in},u^i,v^i\}^{N_{A_{in}}}_{i=1}$ denote the training dataset consisting of initial-boundary points and internal points separately. Correspondingly, $\{\Widehat{u}(x_A^i,t_A^i), \Widehat{v}(x_A^i,t_A^i)\}^{N_A}_{i=1}$ and $\{\Widehat{u}(x_{in}^i,t_{in}^i), \Widehat{v}(x_{in}^i,t_{in}^i)\}^{N_{A_{in}}}_{i=1}$ are the predicted values. In order to calculate $\{f_u(x_{f}^i,t_{f}^i),f_v(x_{f}^i,t_{f}^i)\}^{N_f}_{i=1}$, the derivatives of the networks $u$ and $v$ with respect to time $t$ and space $x$ are derived by automatic differentiation \cite{AD}. Considering that the unknown time-varying variable coefficient $\gamma(t)$ is independent of space $x$ and the objective $\gamma(t)$ takes the form of linear function, we take $N_{\gamma}=1$ and choose $\{ t_0, {\gamma}^0\}$ as the training data.

Similarly, after additionally embedding the term of gradient-enhanced information into the loss function of PINNs, the mean squared error function of gPINNs is given by
\begin{equation}\label{MSE-inv}
MSE_{inverse}^g=MSE_A+MSE_f+MSE_{A_{in}}+MSE_{\gamma}+MSE_{g},
\end{equation}
where
\begin{equation}
MSE_{g}=MSE_{g_u}+MSE_{g_v},	
\end{equation}
\begin{equation}
MSE_{g_u}=\frac{1}{N_g} \left( \sum^{N_g}_{i=1}|\frac{\partial f_u}{\partial t}(x_{g}^i,t_{g}^i)|^2 \right),	
\end{equation}
\begin{equation}
MSE_{g_v}=\frac{1}{N_g} \left( \sum^{N_g}_{i=1}|\frac{\partial f_v}{\partial t}(x_{g}^i,t_{g}^i)|^2 \right),	
\end{equation}
\begin{equation}
\frac{\partial f_u}{\partial t} = -v_{tt}+\alpha(t)_t u_{xx} +\alpha(t) u_{xxt} +\beta(t)_t u +\beta(t) u_t +\gamma(t)_t (u^2+v^2)u	+\gamma(t) (2uu_t+2vv_t)u +\gamma(t) (u^2+v^2)u_t,	
\end{equation}
\begin{equation}
\frac{\partial f_v}{\partial t} = u_{tt}+\alpha(t)_t v_{xx} +\alpha(t) v_{xxt} +\beta(t)_t v +\beta(t) v_t +\gamma(t)_t (u^2+v^2)v+\gamma(t) (2uu_t+2vv_t)v +\gamma(t) (u^2+v^2)v_t.	
\end{equation}
For the time-varying variable coefficient $\gamma(t)$, the gradient-enhanced effect of $t$ is solely considered here by adding mean squared errors involving the partial derivatives of the governing functions with respect to time ($\frac{\partial f_u}{\partial t}$ and $\frac{\partial f_v}{\partial t}$). Besides, the functions $f_u$ and $f_v$ only reflect the value of variable coefficient itself while $\frac{\partial f_u}{\partial t}$ and $\frac{\partial f_v}{\partial t}$ embody the extra derivative information of $\gamma(t)$, i.e., the information of equations satisfied by  $\gamma_t$.

With the aid of the MATLAB software, the spatial region $[-4,4]$ and the temporal region $[-4,4]$ are divided into $N_x=513$ and $N_t=201$ discrete equidistance points, respectively. Thus, the reference one-soliton solution 
\begin{align}
A(x,t)=\frac{\mathrm{e}^{\frac{{\rm{i}}}{10} t^2} \mathrm{e}^{(1+{\rm{i}})x-\frac{t^2}{2}}}{1+\frac{\mathrm{e}^{-t^2+2x}}{4}},	
\end{align}
is discretized into $513 \times 201$ data points in the given spatiotemporal domain. Then $N_A=200$ points are randomly selected from the initial-boundary dataset and $N_{A_{in}}=2000$ points from interior point set. By means of the Latin hypercube sampling method \cite{Latin}, $N_f=N_g=40000$ collocation points are also sampled.

The neural network of the complex valued solution $A(x,t)$ (called as the trunk network) consists of one input layer, 7 hidden layers with 40 neurons per hidden layer and one output layer. The output layer has 2 neurons to learn the real part $u(x,t)$ and imaginary part $v(x,t)$. Given that the function expression of variable coefficient is simpler than that of the solution, we construct the branch network consisting of one input layer, 3 hidden layers as well as one output layer with one neuron to obtain the data-driven variable coefficient $\gamma(t)$ and each hidden layer has 30 neurons. The linear activation function is used in the branch network while $tanh$ function is selected as the activation function in the trunk network. Weights of the neural networks are initialized with Xavier initialization \cite{Xavier}. In addition, we apply L-BFGS algorithm \cite{LBFGS} to minimize the value of the loss function by optimizing the parameters of the neural networks.

To evaluate the performance of three methods (PINNs, gPINNs, TL-gPINNs), we calculate the absolute error and relative error of $\gamma(t)$: the mean absolute error ($MAE$) and relative $\mathbb{L}_2$ error ($RE$) of the variable coefficient $\gamma(t)$
\begin{align}
MAE_{\gamma}=\frac{1}{N_t^\prime} \sum^{N_t^\prime-1}_{k=0} |\Widehat{\gamma}(t_0+k\frac{t_1-t_0}{N_t^\prime-1})-\gamma^{k}|,\\
RE_{\gamma}=\frac{\sqrt{ \sum^{N_t^\prime-1}_{k=0} |\Widehat{\gamma}(t_0+k\frac{t_1-t_0}{N_t^\prime-1})-\gamma^{ k}|^2}}{\sqrt{ \sum^{N_t^\prime-1}_{k=0} |\gamma^{ k}|^2}},
\end{align}
after choosing the corresponding parameter as $N_t^\prime = 500$.

Firstly, the original PINNs is applied. Then, we save the weight matrixes and bias vectors of PINNs at the end of the iteration process to initialize corresponding parameters of gPINNs. After 1862.3727 seconds, the relative $\mathbb{L}_2$ errors of the real part $u$, the imaginary part $v$ and the modulus $|A|$ are 1.331004e-03, 1.407320e-03 and 9.441619e-04 separately. Besides, the mean absolute error and relative $\mathbb{L}_2$ error of the variable coefficient $\gamma(t)$ are: 1.915842e-05 and 9.511436e-06. Obviously, the training by gPINNs is based on training results of PINNs instead of training from scratch, which helps to accelerate convergence to the approximate optimal solution and variable coefficient.

\begin{figure}[htbp]
\centering
\includegraphics[width=7.5cm,height=5cm]{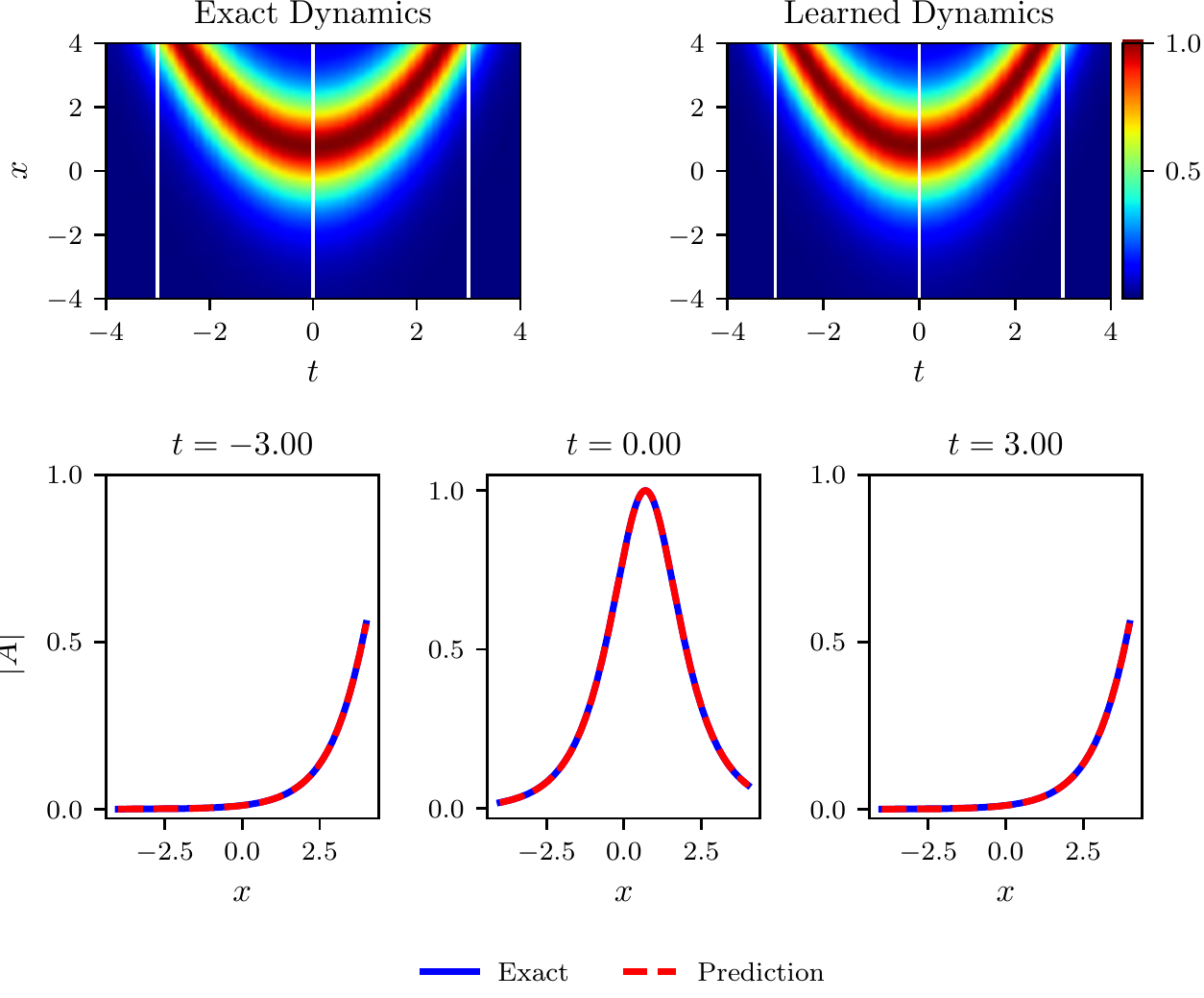}
$a$
\includegraphics[width=7.5cm,height=5cm]{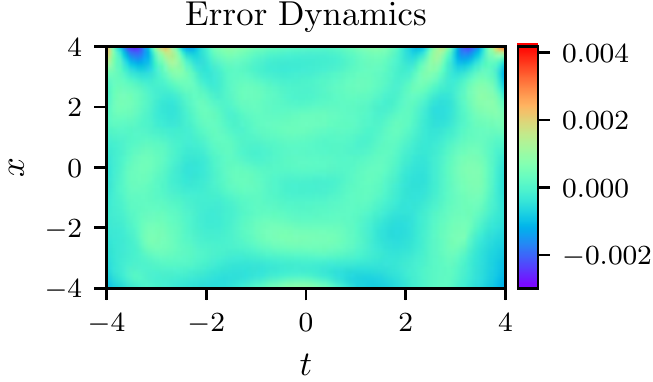}
$b$
\caption{(Color online) One-soliton solution $A(x,t)$ of the vcNLS equation by TL-gPINNs: (a) The density diagrams and comparison between the predicted solutions and exact solutions at the three temporal snapshots of $|A(x,t)|$; (b) The error density diagram of $|A(x,t)|$.}
\label{fig3-1}
\end{figure}

\begin{figure}[htbp]
\centering
\includegraphics[width=7.5cm,height=5cm]{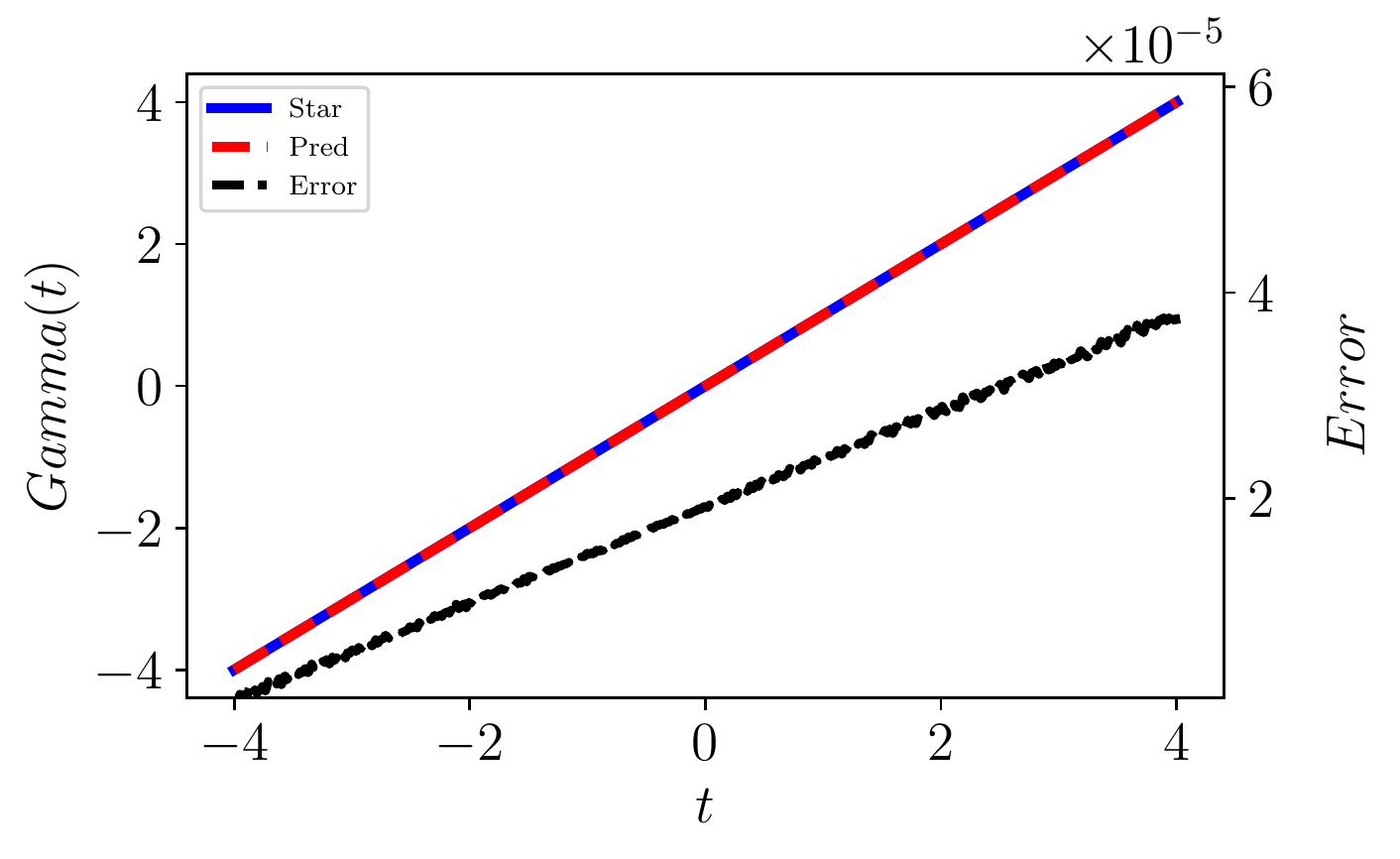}
$a$
\includegraphics[width=7.5cm,height=5cm]{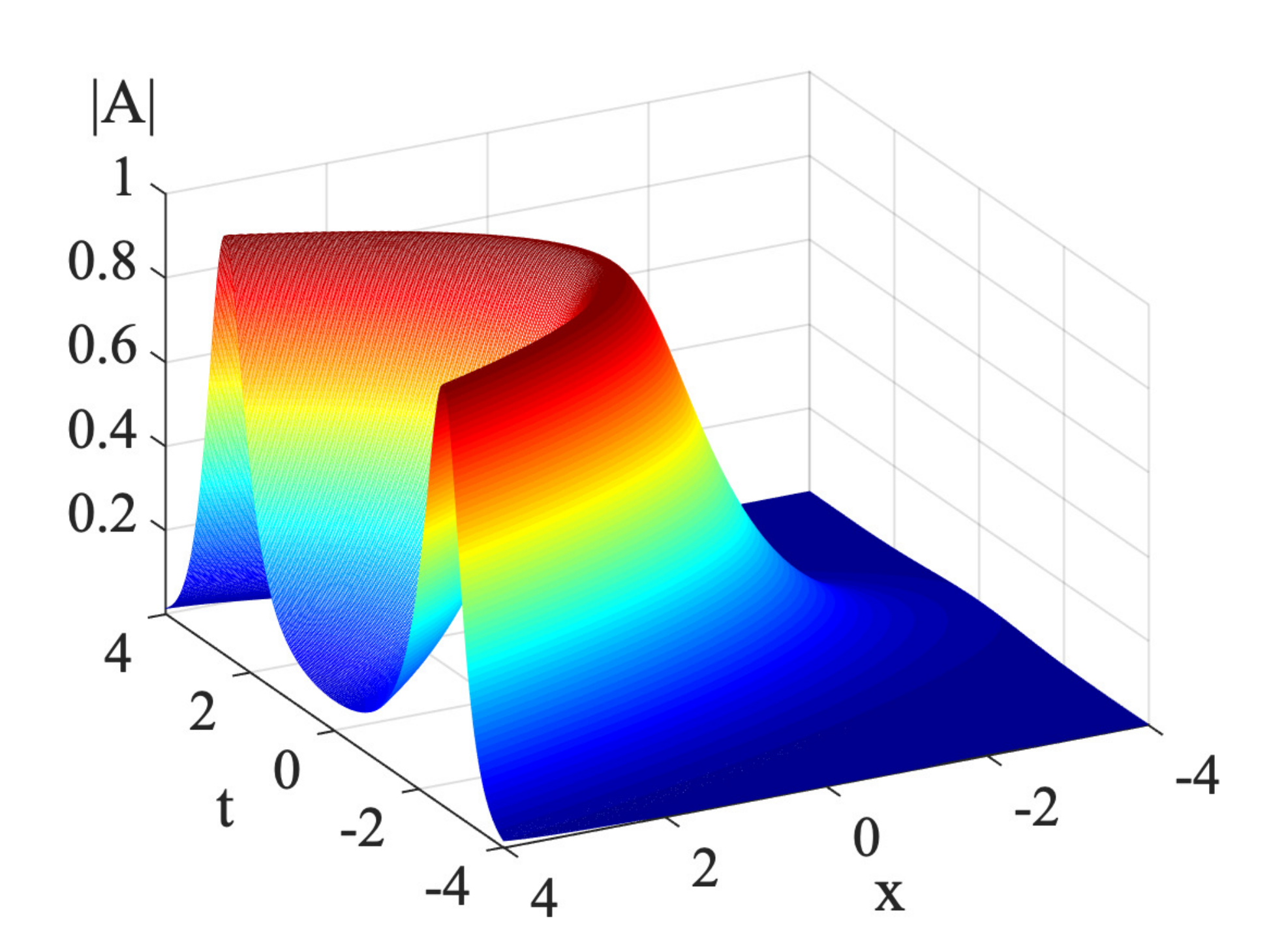}
$b$
\caption{(Color online) Results of function discovery for  the vcNLS equation by TL-gPINNs: (a) The absolute error and comparison between the predicted and exact variable coefficient $\gamma(t)$; (b) The three-dimensional plot of the data-driven one-soliton solution $|A(x,t)|$.}
\label{fig3-2}
\end{figure}

Ultimately, the unknown variable coefficient $\gamma(t)$ is learned simultaneously with the one-soliton solution $A(x,t)$ by TL-gPINNs. Density diagrams of the data-driven one-soliton solution, comparison between the predicted solution and exact solution as well as the evolution plots are shown in Fig. \ref{fig3-1}. It implies there is little difference between  the exact solution and the predicted one. Fig. \ref{fig3-2} (a) is a double coordinate plot, where the solid blue line and the dashed red line corresponding to the left coordinate axis represent the exact and predicted variable coefficient $\gamma(t)$ respectively while the curve of the absolute error is drawn with black dotted line corresponding to the right one. Obviously, the curve of the absolute error exhibits a characteristic of linear variation, and the error is close to 0 at the initial moment since the data of $\gamma(t)$ at $t_0=-4$ is provided and the linear activation function is selected in the branch network.  The predicted 3D plot of the soliton solution with a parabolic shape for the vcNLS equation are shown in Fig. \ref{fig3-2} (b). From the above figures, it can be intuitively seen that the experimental results of $\gamma(t)$ are in good agreement with the theoretical ones.

For intuitive comparison of the prediction accuracy of different methods, the error reduction rate ($ERR$) can be obtained according to the mean absolute error and relative $\mathbb{L}_2$ error achieved by PINNs and gPINNs (TL-gPINNs)
\begin{align}
&ERR_1=\frac{MAE_{\gamma}^{PINNs}-MAE_{\gamma}^{new}}{MSE_{\gamma}^{PINNs}},\\
&ERR_2=\frac{RE_{\gamma}^{PINNs}-RE_{\gamma}^{new}}{RE_{\gamma}^{PINNs}}	,
\end{align}
where 'new' can be replaced by 'gPINNs' and 'TL-gPINNs'.

Finally, the contrast in respect of efficiency and accuracy is presented in Table \ref{table3-1}, including the elapsed time, mean absolute error, relative $\mathbb{L}_2$ error and error reduction rates.

\begin{table}[htbp]
\caption{Performance comparison of three methods: the elapsed time, mean absolute errors and relative $\mathbb{L}_2$ errors of the linear variable coefficient $\gamma(t)$ as well as error reduction rates.}
\label{table3-1}  
\centering
\begin{tabular}{c|c|c|c}
\bottomrule
\diagbox{\textbf{\textbf{Results}}}{\textbf{Method}}  & PINNs & gPINNs & TL-gPINNs \\ \hline
Elapsed time (s)      & 302.1032         & 2579.5674   &  1862.3727        \\ \hline
$MAE_{\gamma}$ & 3.438509e-05 & 2.405417e-05 & 1.915842e-05\\ \hline
$RE_{\gamma}$ & 1.732523e-05 & 1.199509e-05 & 9.511436e-06\\ \hline
$ERR_1$ & - & 30.04\% & 44.28\%  \\ \hline
$ERR_2$ & - & 30.77\% & 45.10\%         \\ \toprule
\end{tabular}
\end{table}

\subsubsection{Data-driven discovery of  quadratic variable coefficient $\gamma(t)$}\label{3.1.2}
\quad

Here, we fix $\alpha(t)=\frac{t^2}{2}, \beta(t) = \frac{t}{5}$ and our objective function is  $\gamma(t)=t^2$ based on the dataset of the corresponding solution $A(x,t)$ for the variable coefficient nonlinear Schr\"{o}odinger equation:
\begin{align}
A(x,t)=\frac{\mathrm{e}^{\frac{{\rm{i}}}{10} t^2} \mathrm{e}^{(1+{\rm{i}})x-\frac{t^3}{3}}}{1+\frac{\mathrm{e}^{2x-\frac{2 t^3}{3}}}{4}}.	
\end{align}

Since the functional form of the target variable coefficient $\gamma(t)$ is quadratic and no longer linear as in Sec. \ref{3.1}, we add sampling data of it and change the term measuring the difference between the predicted values and the true values of $\gamma(t)$ into
\begin{equation}
MSE_{\gamma}=\frac{1}{2}\left(|\Widehat{{\gamma}}(t_0)-{\gamma}^0|^2 + |\Widehat{{\gamma}}(t_1)-{\gamma}^1|^2 \right).
\end{equation}
But apart from that, the loss functions of PINNs and gPINNs (TL-gPINNs) are consistent with the previous subsection.

Obviously, $N_{\gamma}=2$ here and then the training region is selected as $[x_0,x_1]\times [t_0,t_1]=[-4,4]\times [-2,2]$. After exploiting the same data discretization method, we divide the spatial region $[x_0,x_1]=[-4,4]$ into $N_x=513$ discrete equidistance points and the time region $[t_0,t_1]=[-2,2]$ into $N_t=201$ discrete equidistance points. The initial-boundary dataset ($N_A=200$) and the internal point set ($N_{A_{in}}=2000$) are sampled randomly from $513 \times 201$ data points of the solution $A(x,t)$ and extract $N_f=N_g=40000$ collocation points via the Latin hypercube sampling method.

We firstly initialize weights of PINNs with Xavier initialization. A 7-hidden-layer feedforward neural network with 40 neurons per hidden layer and a 3-hidden-layer feedforward neural network with 30 neurons per hidden layer are constructed to learn the one soliton solution and the variable coefficient $\gamma(t)$ of the vcNLS equation, respectively. We use the hyperbolic tangent (tanh) activation function to add nonlinear factors into neural networks. At the end of the iteration process, the parameter data of PINNs is stored and then we use the saved data to fine-tune gPINNs with the same structure by changing the loss function into \eqref{MSE-inv}.

\begin{figure}[htbp]
\centering
\includegraphics[width=7.5cm,height=5cm]{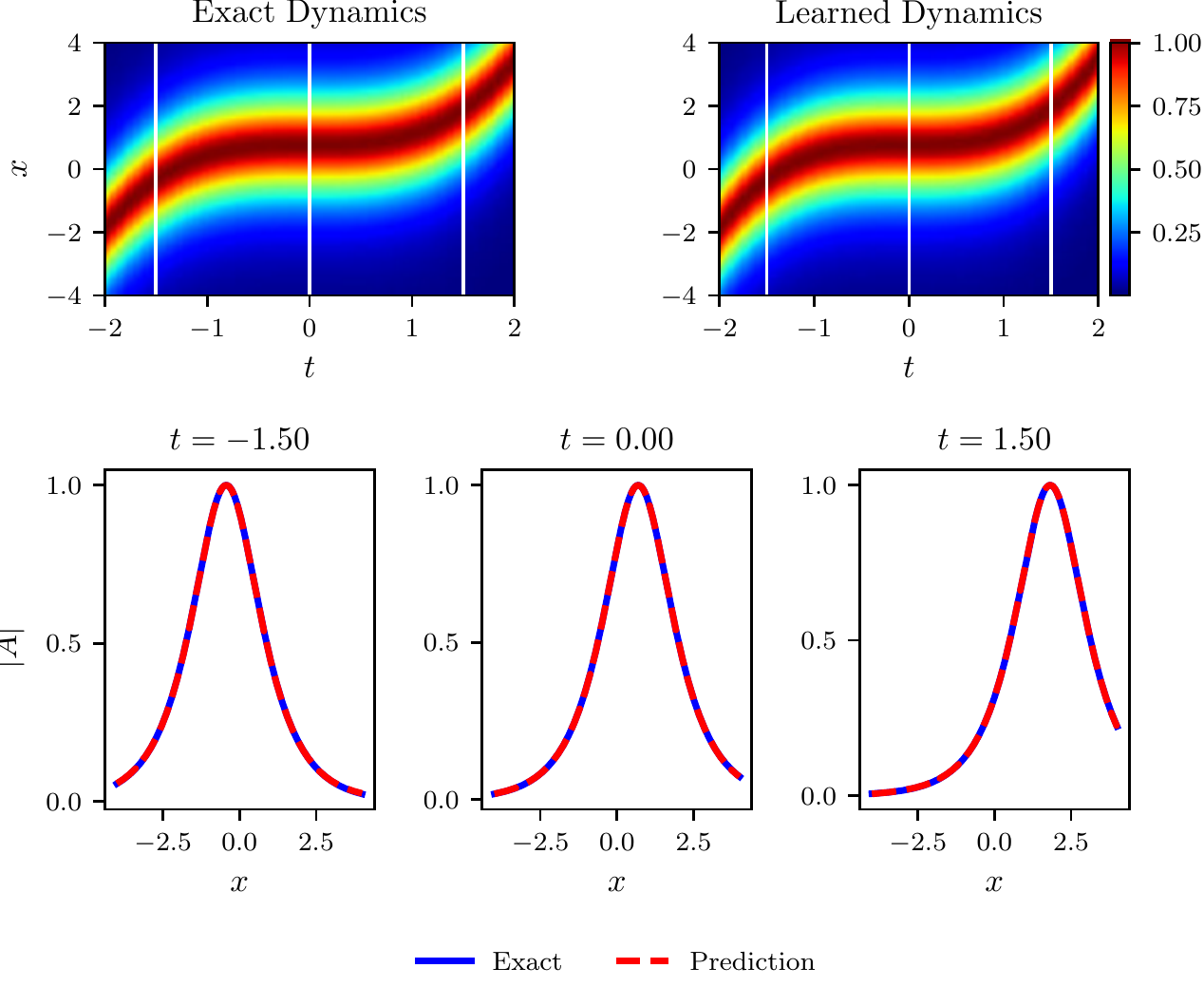}
$a$
\includegraphics[width=7.5cm,height=5cm]{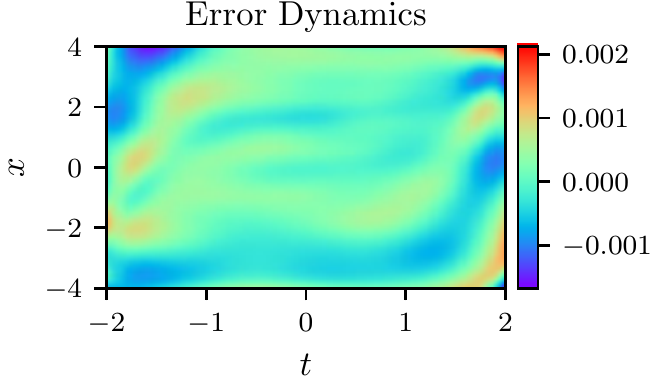}
$b$
\caption{(Color online) One-soliton solution $A(x,t)$ of the vcNLS equation by TL-gPINNs: (a) The density diagrams and comparison between the predicted solutions and exact solutions at the three temporal snapshots of $|A(x,t)|$; (b) The error density diagram of $|A(x,t)|$.}
\label{fig3-3}
\end{figure}

\begin{figure}[htbp]
\centering
\includegraphics[width=7.5cm,height=5cm]{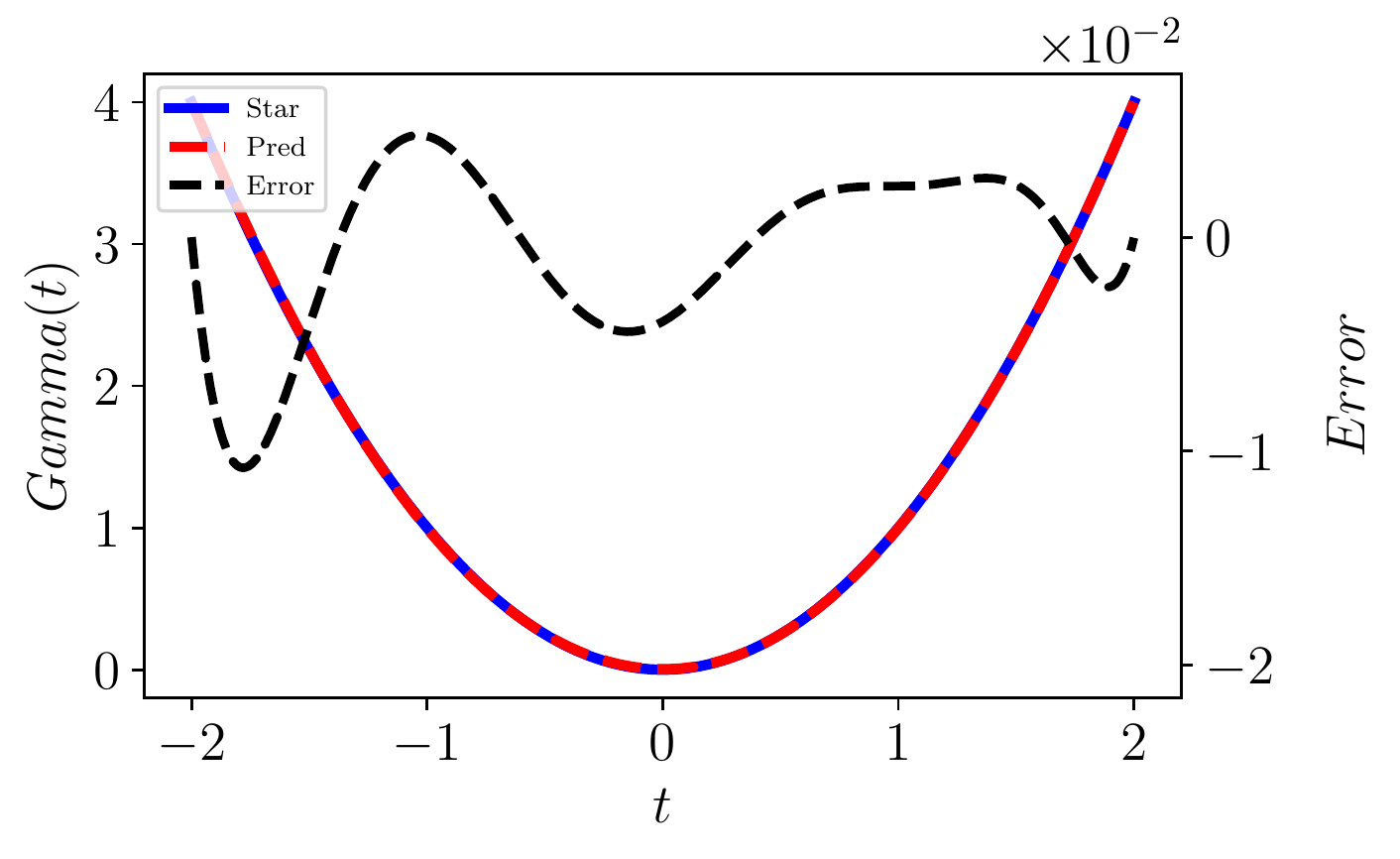}
$a$
\includegraphics[width=7.5cm,height=5cm]{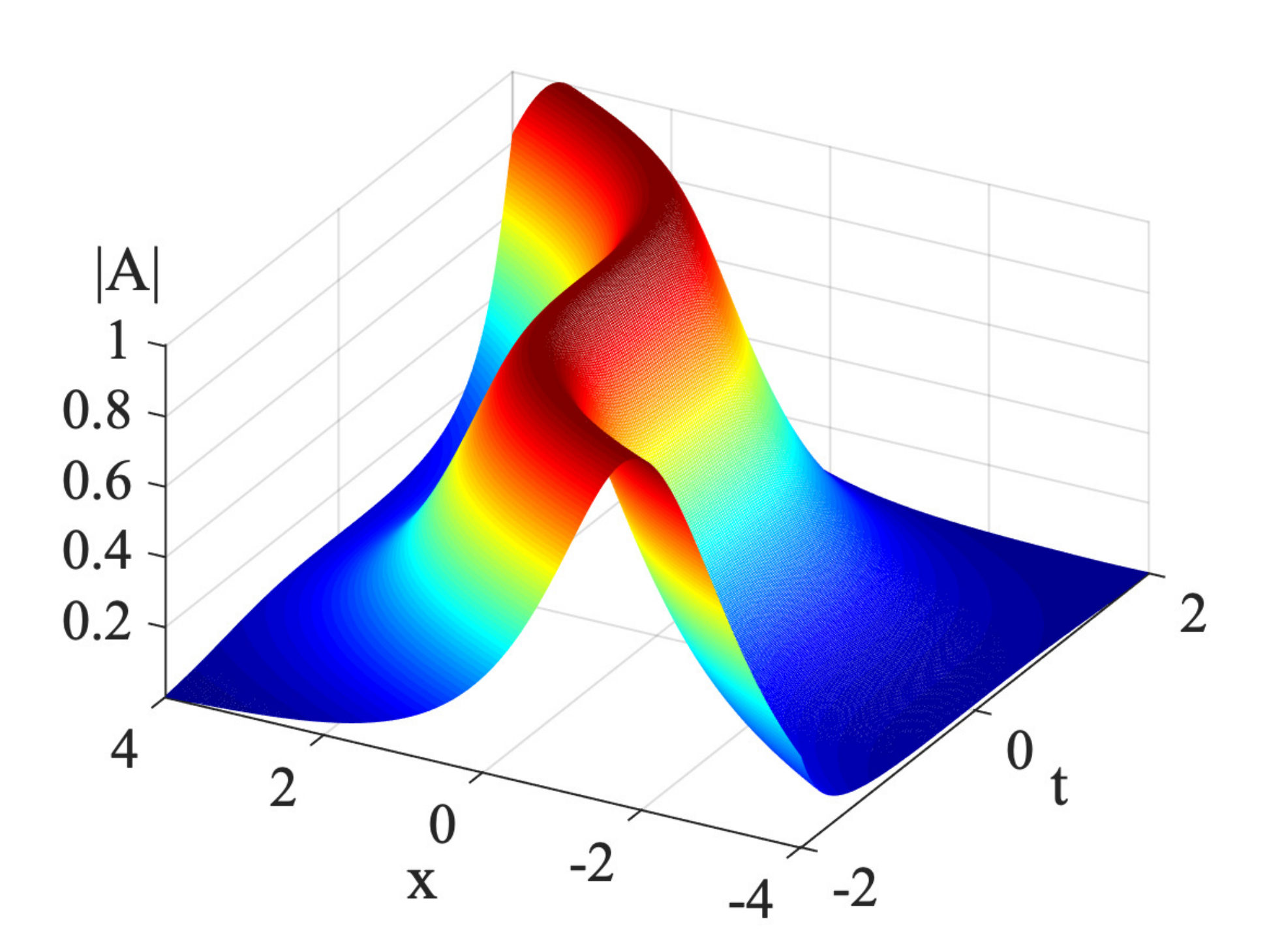}
$b$
\caption{(Color online) Results of function discovery for  the vcNLS equation by TL-gPINNs: (a) The absolute error and comparison between the predicted and exact variable coefficient $\gamma(t)$; (b) The three-dimensional plot of the data-driven one-soliton solution $|A(x,t)|$.}
\label{fig3-4}
\end{figure}

In about 2293.3249 seconds, the data-driven solution of the vcNLS equation is obtained by gPINNs based on transfer learning (TL-gPINNs) and the relative $\mathbb{L}_2$ errors of the real part $u$, the imaginary part $v$ and the modulus $|A|$ are 1.336860e-03, 1.452912e-03 and 8.587186e-04. Simultaneously, the variable coefficient $\gamma(t)$ is successfully inferred with the mean absolute error of 3.100830e-03 and relative $\mathbb{L}_2$ error of 2.163681e-03. Fig. \ref{fig3-3} displays the exact, learned and error density diagrams as well as evolution plots of one-soliton solution at different time points $t=-1.5, 0, 1.5$. In Fig. \ref{fig3-4}, the curve plots of the predicted and the exact variable coefficient $\gamma(t)$, the absolute error and the predicted 3D graph of the cubic soliton solution for the vcNLS equation are plotted. As can be seen from these diagrams and performance comparison of three methods shown in Table \ref{table3-2}, TL-gPINN is capable of correctly identifying the unknown variable coefficient $\gamma(t)$ and learning the cubic soliton solution with very high accuracy while gPINN doesn't work as expected.

\begin{table}[htbp]
\caption{Performance comparison of three methods: the elapsed time, mean absolute errors and relative $\mathbb{L}_2$ errors of the quadratic variable coefficient $\gamma(t)$ as well as error reduction rates.}
\label{table3-2}  
\centering
\begin{tabular}{c|c|c|c}
\bottomrule
\diagbox{\textbf{\textbf{Results}}}{\textbf{Method}}  & PINNs & gPINNs & TL-gPINNs \\ \hline
Elapsed time (s)      &  789.8311        & 3133.7038    & 2293.3249          \\ \hline
$MAE_{\gamma}$ &  4.211790e-03    & 1.072530e-02     &  3.100830e-03  \\ \hline
$RE_{\gamma}$ & 3.003559e-03     &  7.299052e-03    & 2.163681e-03    \\ \hline
$ERR_1$ & - & -154.65\%      & 26.38\%     \\ \hline
$ERR_2$ & - & -143.01\%      &  27.96\%        \\ \toprule
\end{tabular}
\end{table}

\subsubsection{Data-driven discovery of  sine variable coefficient $\gamma(t)$}
\quad

After fixing $\alpha(t)=\sin(t), \beta(t) = \frac{t}{5}$, we aim to infer the unknown $\gamma(t)$ in the variable coefficient nonlinear Schr\"{o}odinger equation based on the solution data:
\begin{align}
A(x,t)=\frac{\mathrm{e}^{\frac{{\rm{i}}}{10} t^2} \mathrm{e}^{(1+{\rm{i}})x +2 \cos(t)}}{1+\frac{\mathrm{e}^{2x+4\cos(t)}}{8}}.	
\end{align}

For simplicity, we confine our sampling and training in a rectangular region $(x,t) \in [-4,4]\times [-5,5]$. To generate a dataset for this example, we choose $N_A=200$ points from the initial-boundary dataset and $N_{A_{in}}=2000$ points from interior point set at random after equidistant discretization. In addition, we employ the Latin hypercube sampling method to select $N_f=N_g=40000$ collocation points.

Similarly, we also establish the fully-connected PINNs with Xavier initialization at first and proceed by adopting gPINNs with the advantage of transfer learning. The structure of networks, including the width and depth, activation function, definition of loss functions as well as the optimization algorithm, is the same as the previous subsection.

\begin{figure}[htbp]
\centering
\includegraphics[width=7.5cm,height=5cm]{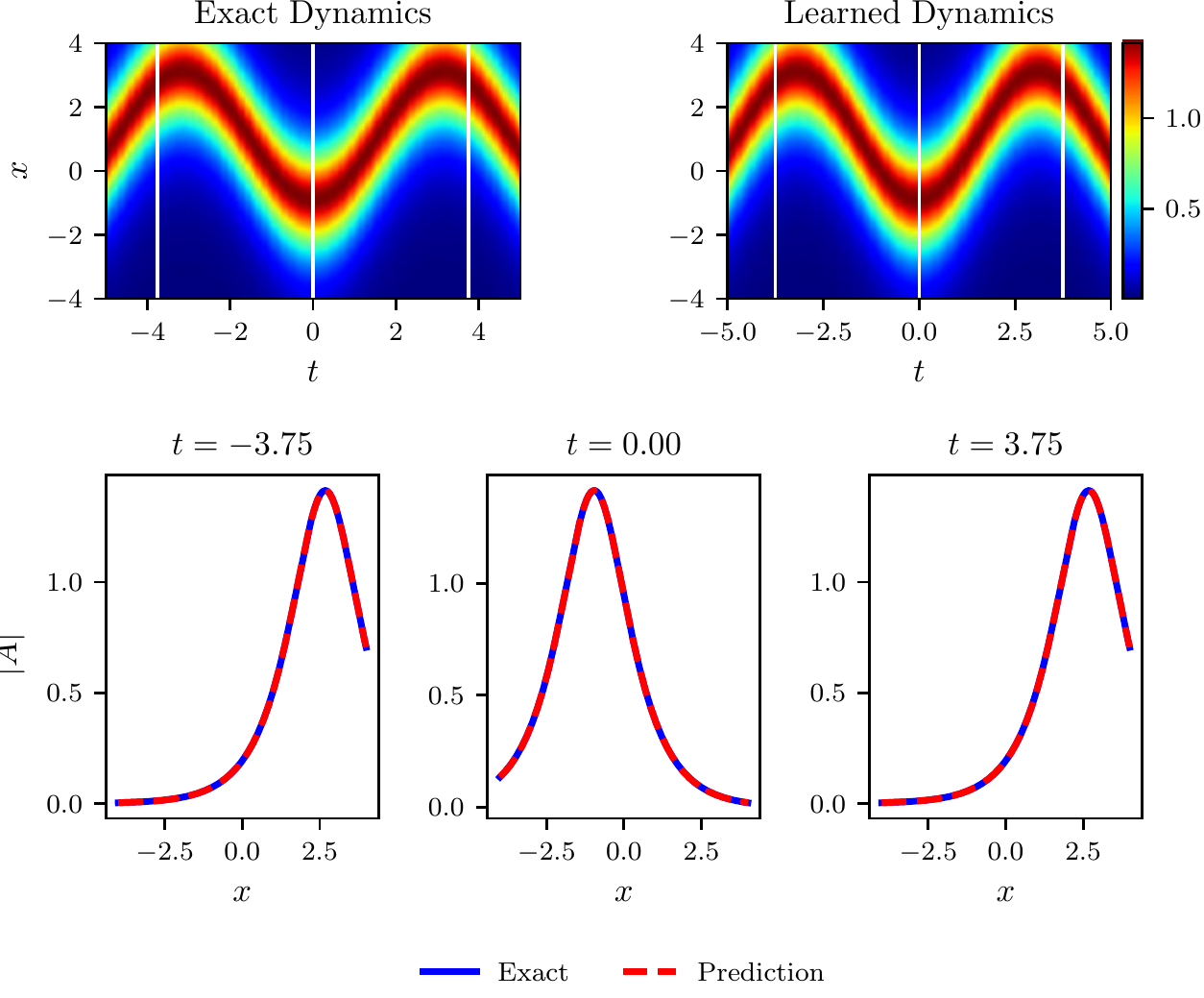}
$a$
\includegraphics[width=7.5cm,height=5cm]{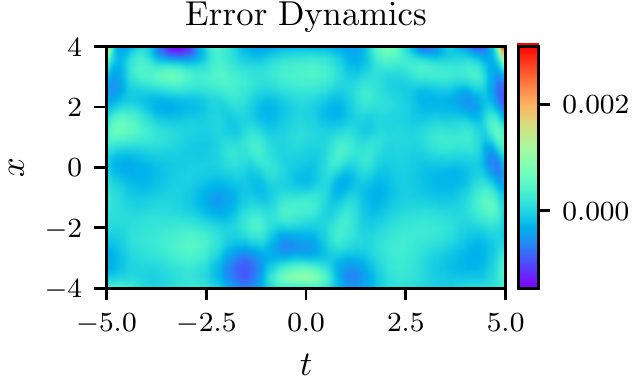}
$b$
\caption{(Color online) One-soliton solution $A(x,t)$ of the vcNLS equation by TL-gPINNs: (a) The density diagrams and comparison between the predicted solutions and exact solutions at the three temporal snapshots of $|A(x,t)|$; (b) The error density diagram of $|A(x,t)|$.}
\label{fig3-5}
\end{figure}

\begin{figure}[htbp]
\centering
\includegraphics[width=7.5cm,height=5cm]{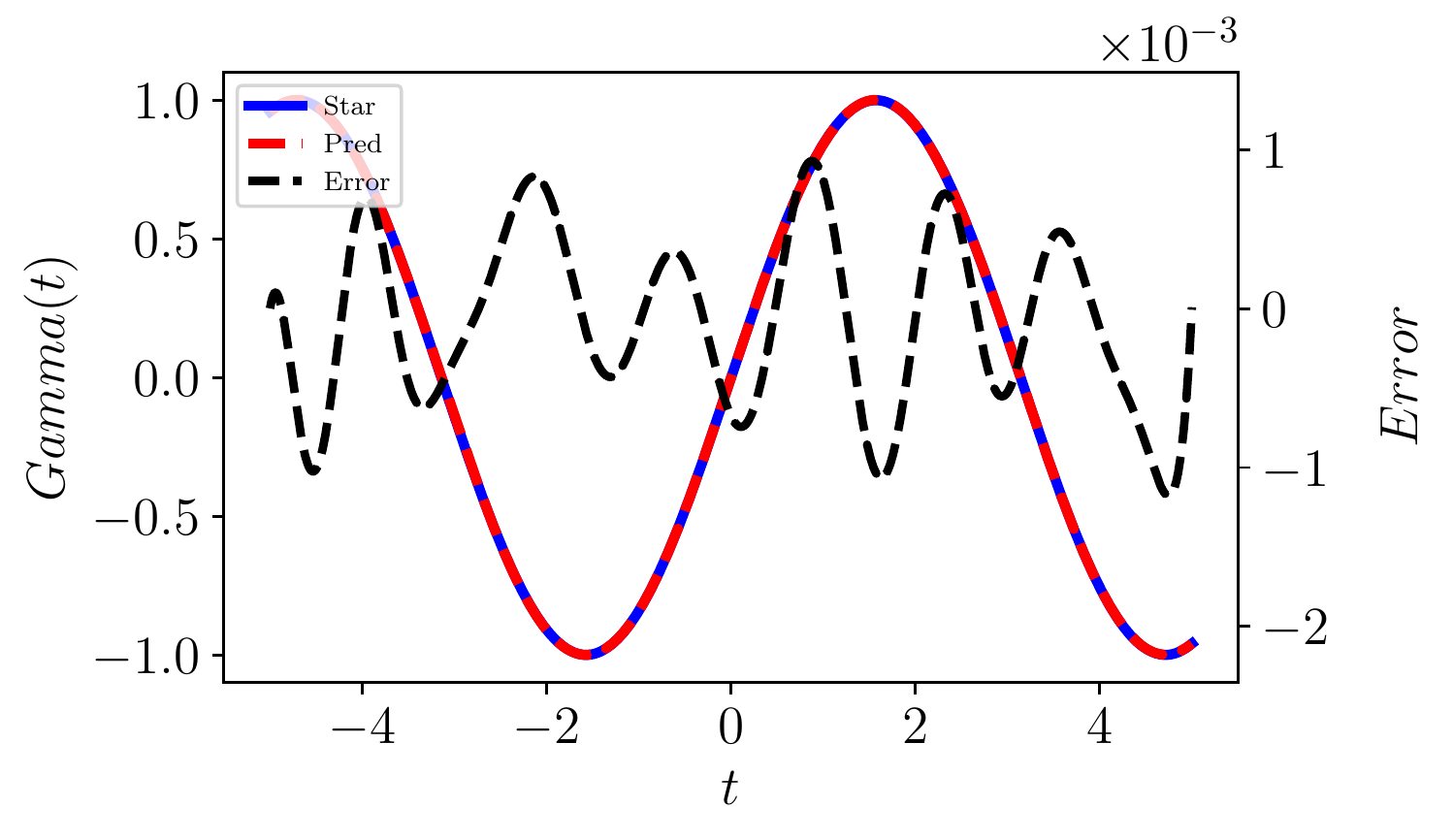}
$a$
\includegraphics[width=7.5cm,height=5cm]{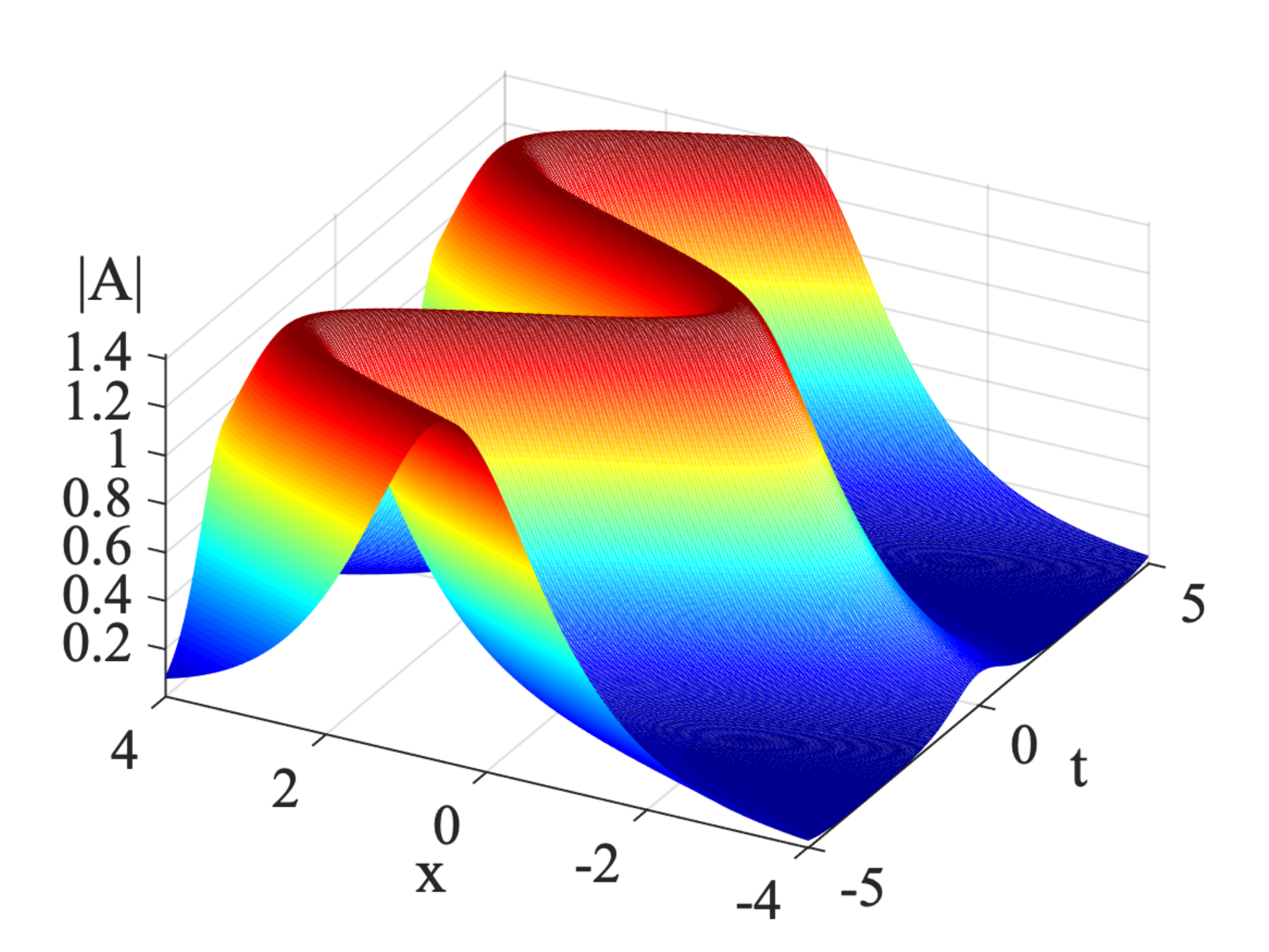}
$b$
\caption{(Color online) Results of function discovery for  the vcNLS equation by TL-gPINNs: (a) The absolute error and comparison between the predicted and exact variable coefficient $\gamma(t)$; (b) The three-dimensional plot of the data-driven one-soliton solution $|A(x,t)|$.}
\label{fig3-6}
\end{figure}

Dynamic behaviors of the soliton solution $A(x,t)$ and variable coefficient $\gamma(t)$ inferred by TL-gPINNs are shown in Fig.\ref{fig3-5} and Fig.\ref{fig3-6}, which contain comparison between the predicted solutions and exact ones, the three-dimensional plots of predicted $A(x,t)$ and the curve graph of the variable coefficient $\gamma(t)$. The absolute error curve of $\gamma(t)$ is drawn with black dotted line corresponding to the right coordinate axis in Fig.\ref{fig3-6} (a). An empirical inference is given that the high-frequency oscillation of absolute error is caused by the periodic oscillation and the change in concavity and convexity of the variable coefficient. It can be observed that we obtain a periodical soliton solution as sine or cosine function and the predicted  variable coefficient is well fitted with the exact one with absolute error less than $2 \times 10^{-3}$. In addition, based on the results in Table \ref{table3-3}, it illustrates that both the mean absolute error and relative $\mathbb{L}_2$ error of the variable coefficient $\gamma(t)$ achieved by TL-gPINNs reach the level of $10^{-4}$, about one order of magnitude lower than those by PINNs.

\begin{table}[htbp]
\caption{Performance comparison of three methods: the elapsed time, mean absolute errors and relative $\mathbb{L}_2$ errors of the sine variable coefficient $\gamma(t)$ as well as error reduction rates.}
\label{table3-3}  
\centering
\begin{tabular}{c|c|c|c}
\bottomrule
\diagbox{\textbf{\textbf{Results}}}{\textbf{Method}}  & PINNs & gPINNs & TL-gPINNs \\ \hline
Elapsed time (s)      & 729.1269         & 5283.3492    & 4272.9441          \\ \hline
$MAE_{\gamma}$ &  1.463990e-03    & 7.123562e-04     &  4.664226e-04 \\ \hline
$RE_{\gamma}$ & 2.703498e-03     & 1.363048e-03     &  7.559607e-04 \\ \hline
$ERR_1$ & - &  51.34\%     & 68.14\%     \\ \hline
$ERR_2$ & - &  49.58\%      &  72.04\%        \\ \toprule
\end{tabular}
\end{table}

\subsubsection{Data-driven discovery of  hyperbolic tangent variable coefficient $\gamma(t)$}
\quad

Given $\alpha(t)=\tanh(t), \beta(t) = \frac{t}{5}$, our goal is to identify the unknown variable parameter $\gamma(t)$ from the vcNLS equation with remarkable accuracy.

After utilizing the same generation and sampling method of training data as above, we acquire the training set consists of $N_A=200$ initial-boundary points,  $N_{A_{in}}=2000$ internal points and a random selection of $N_f=N_g=40000$ collocation points in the given spatiotemporal domain $[x_0,x_1]\times [t_0,t_1]=[-2,4]\times [-5,5]$ where the corresponding soliton solution is 
\begin{align}
A(x,t)=\frac{\mathrm{e}^{\frac{{\rm{i}}}{10} t^2} \mathrm{e}^{(1+{\rm{i}})x -2 \ln(\cosh(t))}}{1+\frac{\mathrm{e}^{2x-4\ln(\cosh(t))}}{8}}	.
\end{align}

The first step is to construct the conventional PINNs. The architecture of multi-out neural networks consists of one input layer, 7 hidden layers with 40 neurons per hidden layer and one output layer with 2 neurons to learn the real part $u(x,t)$ and imaginary part $v(x,t)$ of the soliton solution. A 3-hidden-layer feedforward neural network with 30 neurons per hidden layer is employed to infer the variable parameter $\gamma(t)$. This process can be regarded as the pre-training of the gPINNs, which helps accelerate the convergence of training. Next, we initialize gPINNs with the saved weights of PINNs. The activation function and optimization algorithm used here are the $tanh$ function and L-BFGS optimizer respectively.

\begin{figure}[htbp]
\centering
\includegraphics[width=7.5cm,height=5cm]{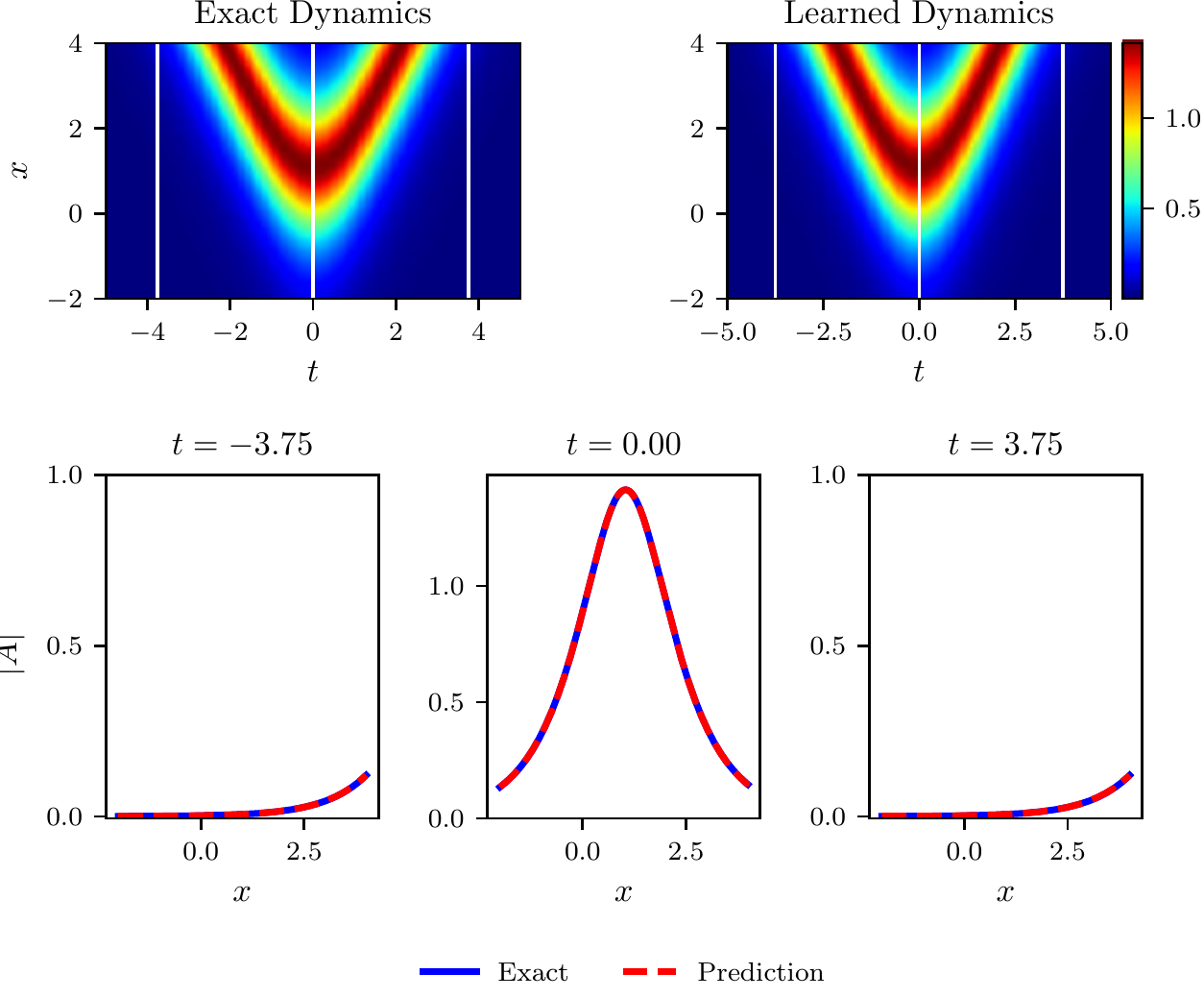}
$a$
\includegraphics[width=7.5cm,height=5cm]{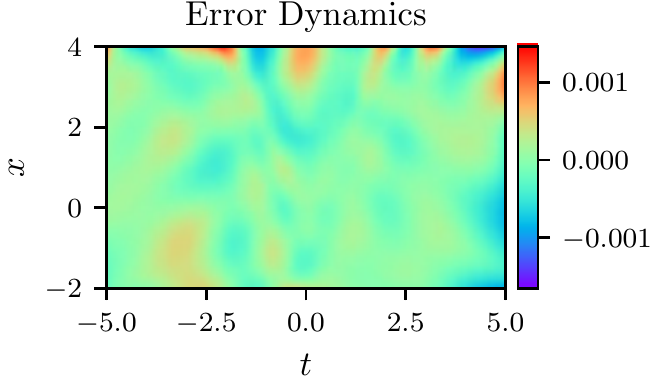}
$b$
\caption{(Color online) One-soliton solution $A(x,t)$ of the vcNLS equation by TL-gPINNs: (a) The density diagrams and comparison between the predicted solutions and exact solutions at the three temporal snapshots of $|A(x,t)|$; (b) The error density diagram of $|A(x,t)|$.}
\label{fig3-7}
\end{figure}

\begin{figure}[htbp]
\centering
\includegraphics[width=7.5cm,height=5cm]{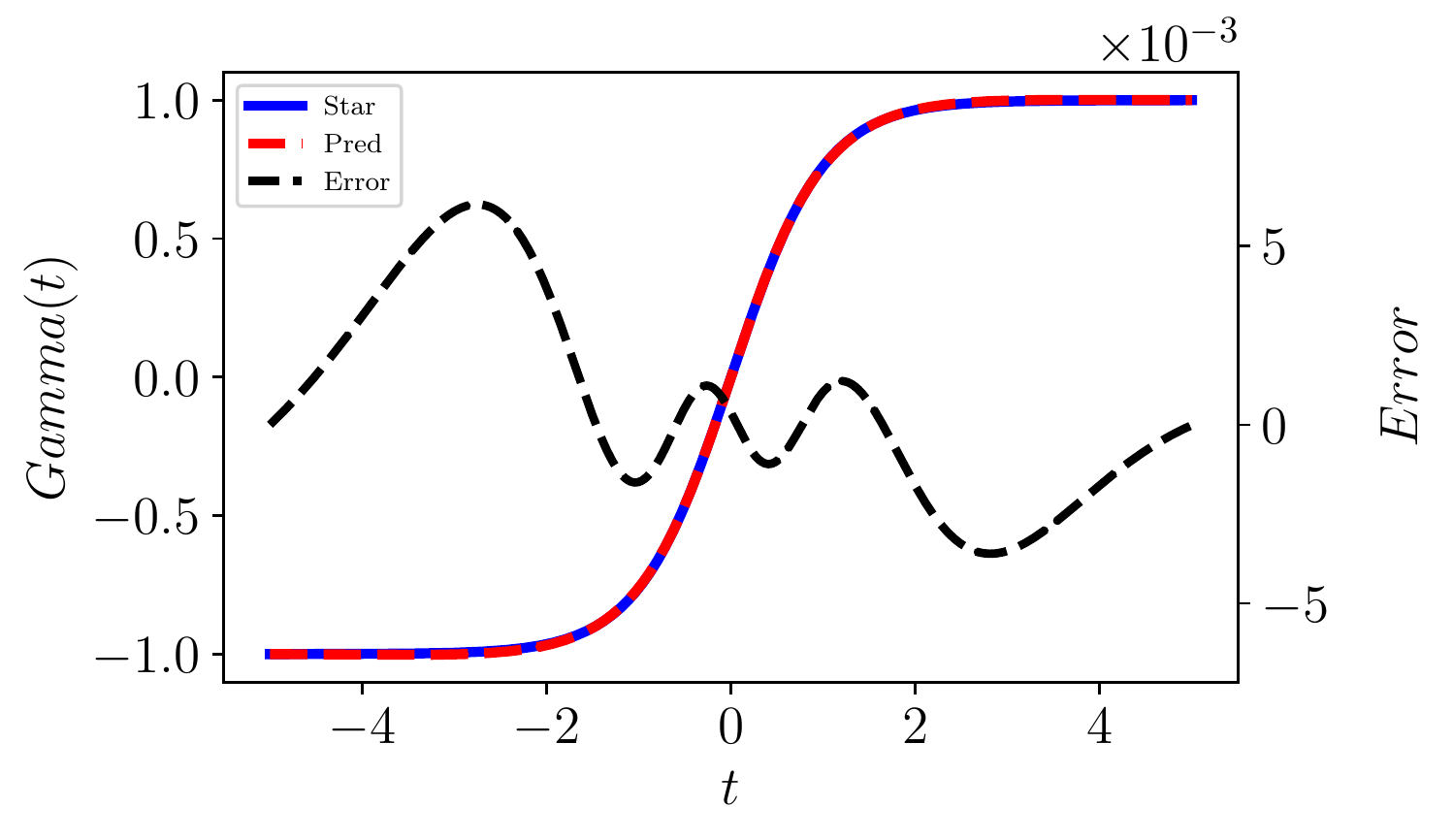}
$a$
\includegraphics[width=7.5cm,height=5cm]{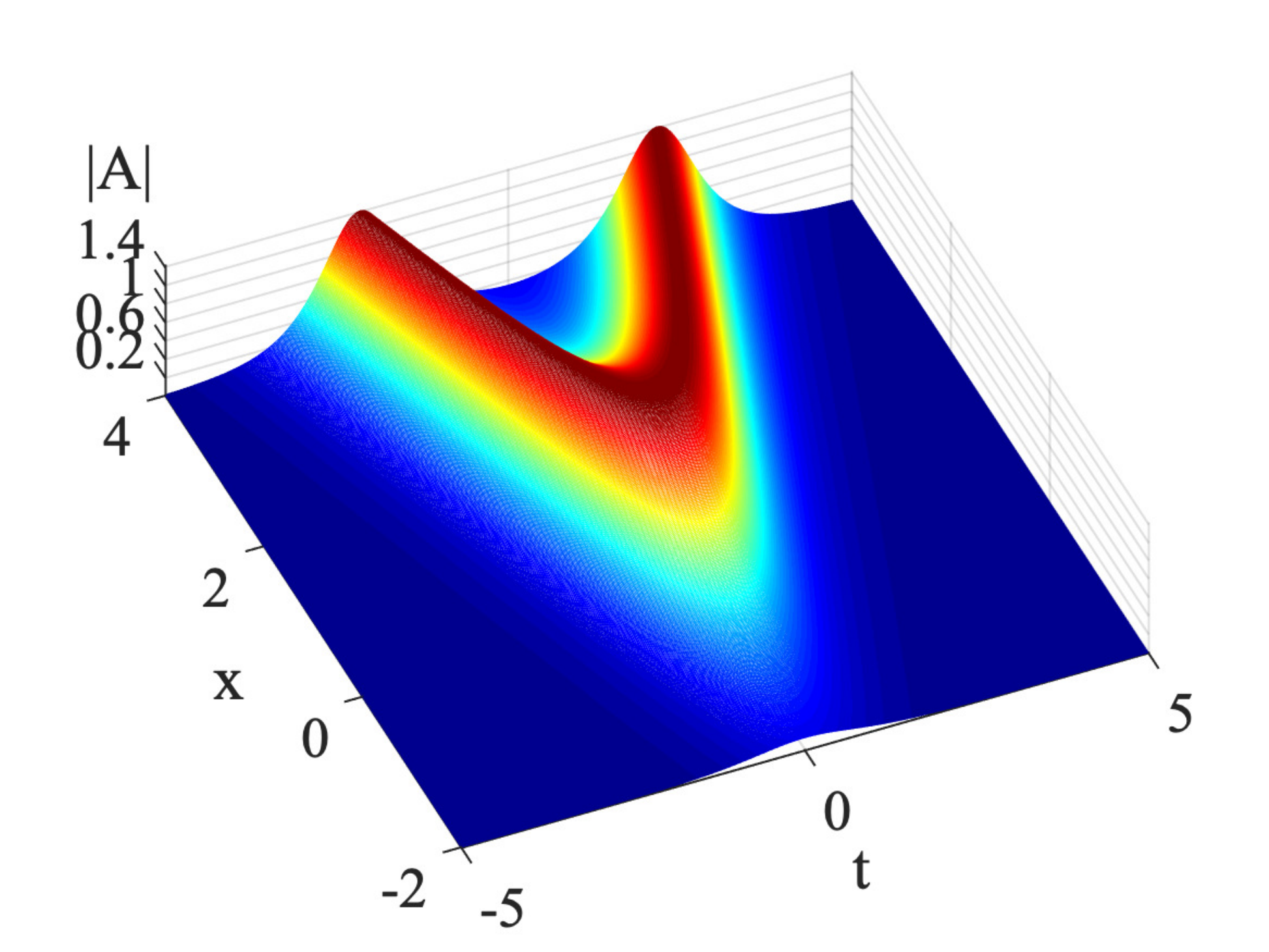}
$b$
\caption{(Color online) Results of function discovery for the vcNLS equation by TL-gPINNs: (a) The absolute error and comparison between the predicted and exact variable coefficient $\gamma(t)$; (b) The three-dimensional plot of the data-driven one-soliton solution $|A(x,t)|$.}
\label{fig3-8}
\end{figure}

By leveraging TL-gPINNs, the data-driven soliton solution $A(x,t)$ and variable coefficient $\gamma(t)$ are plotted in Fig. \ref{fig3-7} and Fig. \ref{fig3-8}. For the double coordinate plot in Fig. \ref{fig3-8} (a), the black dotted line corresponding to the right coordinate axis represents the absolute error curve, which exhibits a certain degree of symmetry since the variable coefficient itself is centrosymmetric. Empirically speaking, the error will increase accordingly when the value of the function to be learned is large or changes greatly. However, the error is relatively small during the period with high slopes, i.e. $t \in [-2,2]$. Presumably it's because the introduction of gradient information into the loss function is conducive to learn the features of variable coefficient where the slope is relatively large. We observe that this V-shaped soliton and the variable coefficient with the function type of hyperbolic tangent are both accurately inferred. Furthermore, Table \ref{table3-4} gives a brief overview of accuracy and efficiency of three methods.

\begin{table}[htbp]
\caption{Performance comparison of three methods: the elapsed time, mean absolute errors and relative $\mathbb{L}_2$ errors of the hyperbolic tangent variable coefficient $\gamma(t)$ as well as error reduction rates.}
\label{table3-4}  
\centering
\begin{tabular}{c|c|c|c}
\bottomrule
\diagbox{\textbf{\textbf{Results}}}{\textbf{Method}}  & PINNs & gPINNs & TL-gPINNs \\ \hline
Elapsed time (s)      & 223.075         & 1621.832    &  1014.9962         \\ \hline
$MAE_{\gamma}$ & 6.833609e-03     & 2.909324e-03     &  2.209660e-03  \\ \hline
$RE_{\gamma}$ & 9.421331e-03     & 3.897784e-03     &  3.178629e-03   \\ \hline
$ERR_1$ & - &  57.43\%     & 67.66\%     \\ \hline
$ERR_2$ & - &   58.63\%     &  66.26\%        \\ \toprule
\end{tabular}
\end{table}

\subsubsection{Data-driven discovery of  fractional variable coefficient $\gamma(t)$}\label{3.1.5}
\quad

When $\alpha(t)$, $\beta(t)$ are respectively fixed as $\frac{1}{2(1+t^2)}$, $\frac{t}{5}$ and the training of this case is confined in a rectangular region $(x,t) \in [-4,5]\times [-5,5]$, the target here is to infer the unknown variable coefficient $\gamma(t)$ on the basis of the dataset of the corresponding solution
\begin{align}
A(x,t)=\frac{\mathrm{e}^{\frac{{\rm{i}}}{10} t^2 }\mathrm{e}^{(1+{\rm{i}})x -\arctan(t)}}{1+\frac{(2t^2+2) \mathrm{e}^{2x-2 \arctan(t)}}{8(t^2+1)}}.
\end{align}

Since there are large amounts of descriptions of the sampling method and network structure above, we won't reiterate them here to avoid repetition. All details are the same as the previous subsection.

\begin{figure}[htbp]
\centering
\includegraphics[width=7.5cm,height=5cm]{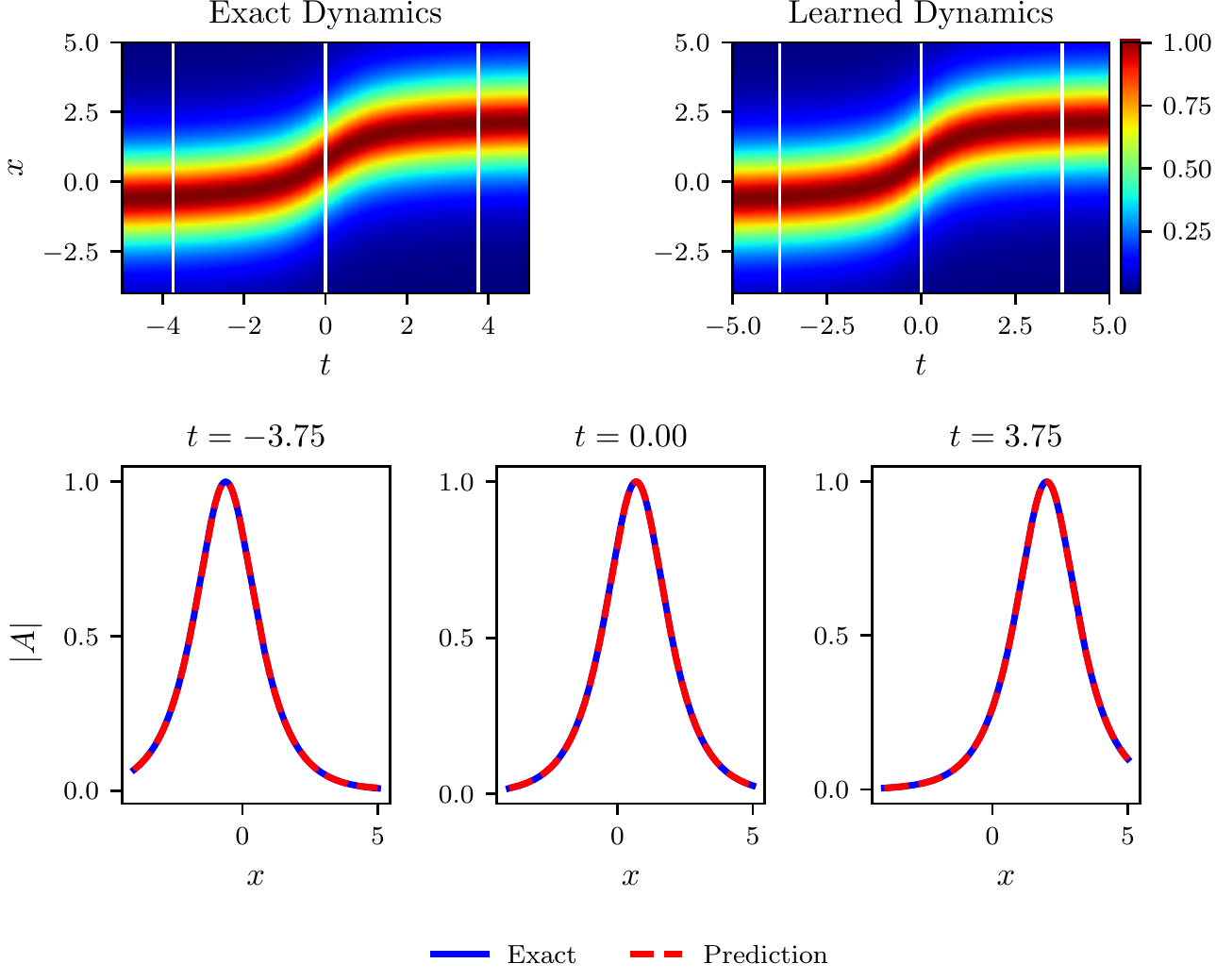}
$a$
\includegraphics[width=7.5cm,height=5cm]{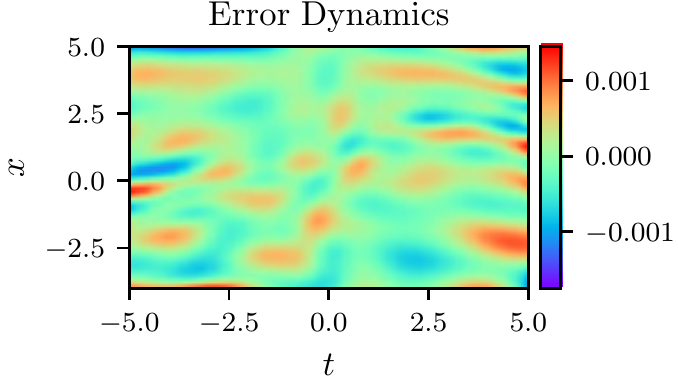}
$b$
\caption{(Color online) One-soliton solution $A(x,t)$ of the vcNLS equation by TL-gPINNs: (a) The density diagrams and comparison between the predicted solutions and exact solutions at the three temporal snapshots of $|A(x,t)|$; (b) The error density diagram of $|A(x,t)|$.}
\label{fig3-9}
\end{figure}

\begin{figure}[htbp]
\centering
\includegraphics[width=7.5cm,height=5cm]{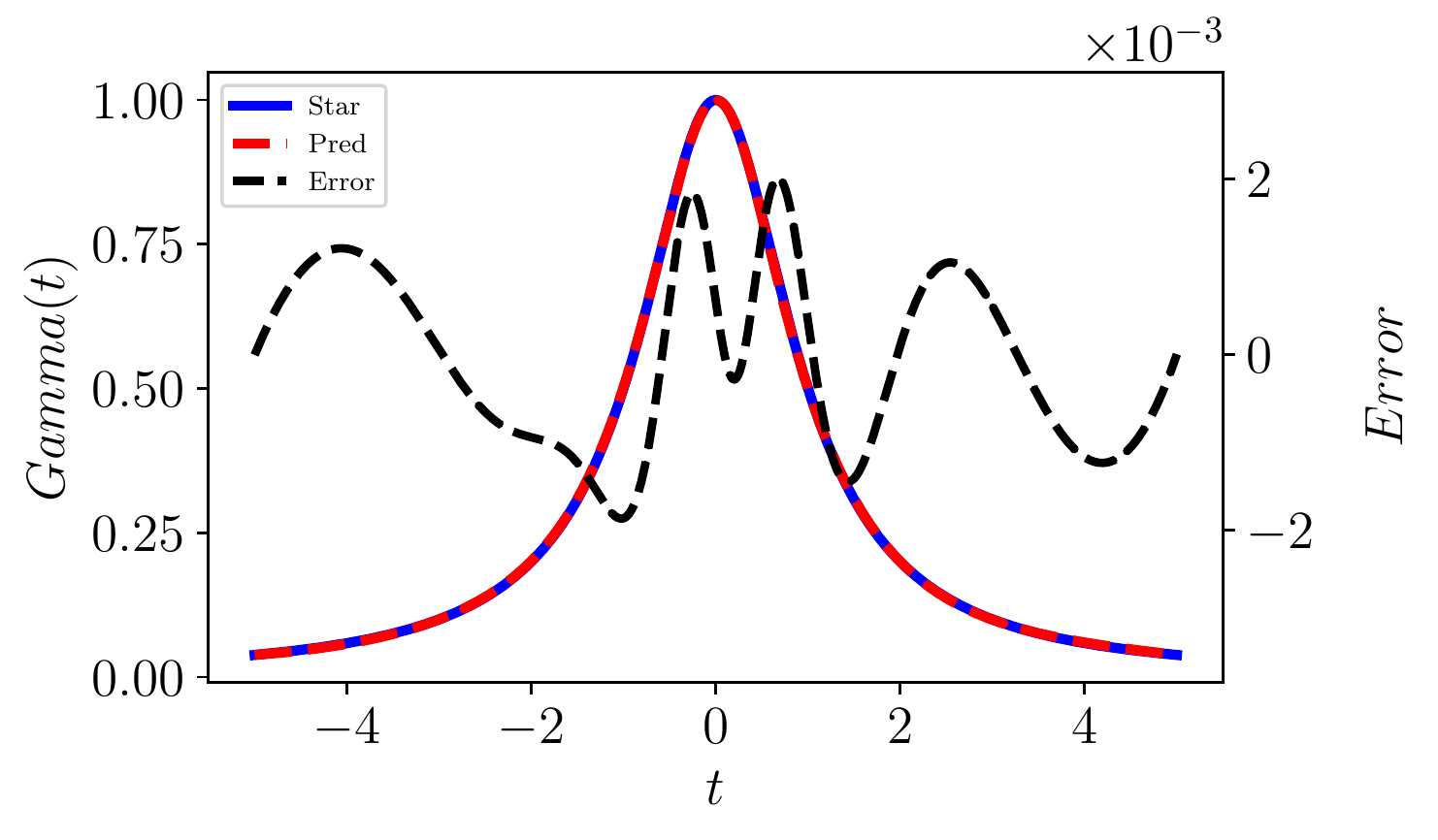}
$a$
\includegraphics[width=7.5cm,height=5cm]{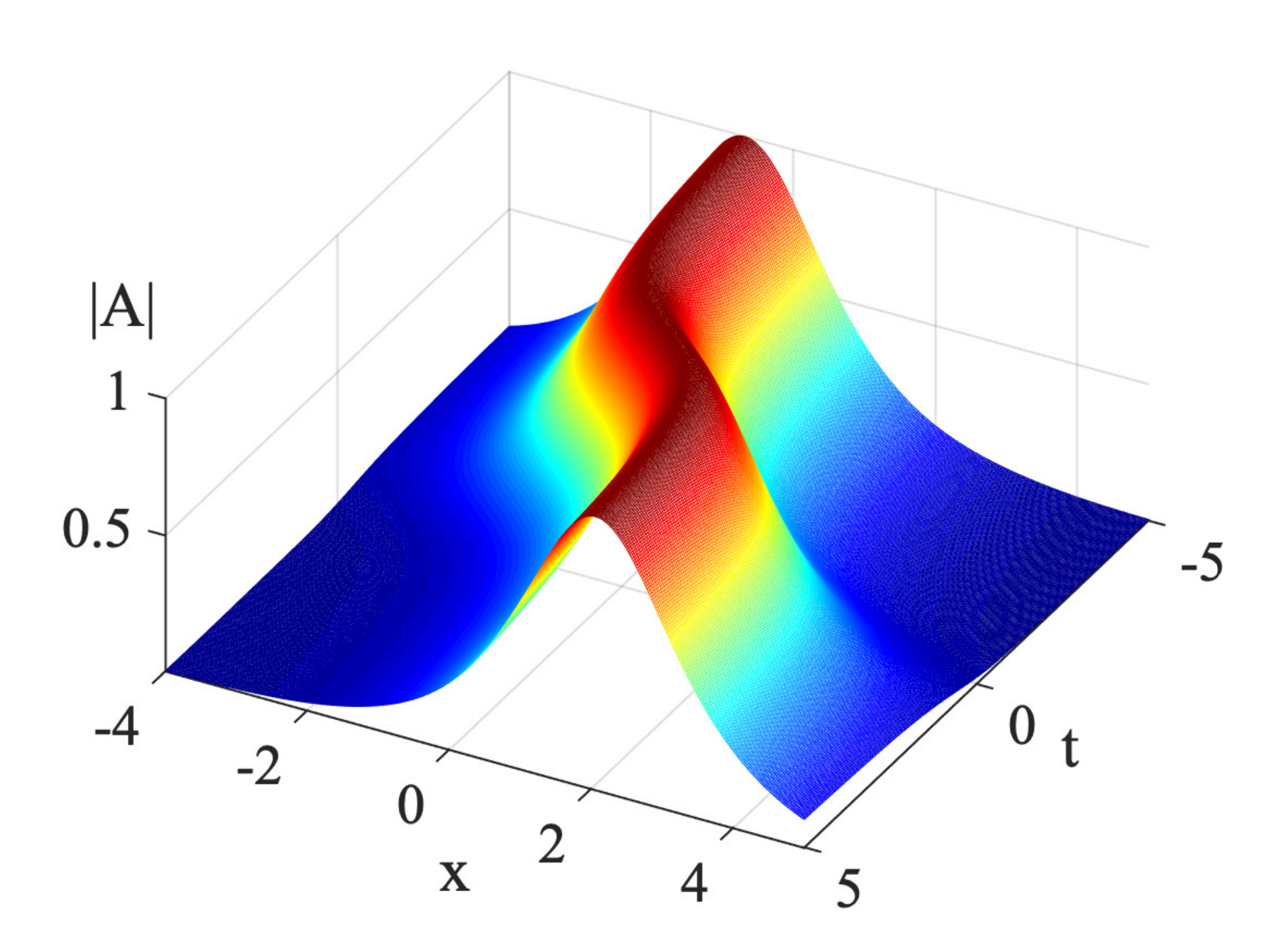}
$b$
\caption{(Color online) Results of function discovery for  the vcNLS equation by TL-gPINNs: (a) The absolute error and comparison between the predicted and exact variable coefficient $\gamma(t)$; (b) The three-dimensional plot of the data-driven one-soliton solution $|A(x,t)|$.}
\label{fig3-10}
\end{figure}

Table \ref{table3-5} summarizes the results of our experiment and compares the performance of PINNs, TL-gPINNs and gPINNs. A more detailed assessment of the predicted soliton solution $A(x,t)$ and variable coefficient $\gamma(t)$ by leveraging TL-gPINNs is presented in Fig.\ref{fig3-9} and Fig.\ref{fig3-10}. Specifically, the comparison between the exact and the predicted solutions at different time points $t=-3.75, 0, 3.75$ as well as that between the predicted and exact variable coefficient $\gamma(t)$ is also displayed. A rule of thumb is that the error is large when the value of variable coefficient $\gamma(t)$ is large or $\gamma(t)$ changes sharply. The change of the absolute error curve plotted with black dashed line shown in Fig.\ref{fig3-10} (a) is in good agreement with this experiential conclusion to a certain extent. In addition, TL-gPINN is capable of accurately capturing the intricate nonlinear behaviors of the vcNLS equation, including the dynamic behaviors of the solution and the Kerr nonlinearity $\gamma(t)$.

\begin{table}[htbp]
\caption{Performance comparison of three methods: the elapsed time, mean absolute errors and relative $\mathbb{L}_2$ errors of the fractional variable coefficient $\gamma(t)$ as well as error reduction rates.}
\label{table3-5}  
\centering
\begin{tabular}{c|c|c|c}
\bottomrule
\diagbox{\textbf{\textbf{Results}}}{\textbf{Method}}  & PINNs & TL-gPINNs & gPINNs \\ \hline
Elapsed time (s)      & 155.3069         & 612.4961     & 327.6574           \\ \hline
$MAE_{\gamma}$ & 1.504344e-03    & 1.662907e-03     &  8.788681e-04 \\ \hline
$RE_{\gamma}$ & 5.177898e-03     &  5.403644e-03   & 2.536767e-03   \\ \hline
$ERR_1$ & - & -10.54\%      &  41.58\%     \\ \hline
$ERR_2$ & - &  -4.36\%     & 51.01\%        \\ \toprule
\end{tabular}
\end{table}

\subsection{Data-driven discovery of multiple variable coefficients}
\quad

We extend the research of data-driven discovery for single variable coefficient to that of multiple ones, and the hyper-parameters of which are given in outline in Table \ref{table3-0}. For each case discussed here, the L-BFGS algorithm is utilized to optimize loss functions.

\subsubsection{Data-driven discovery of two variable coefficients: linear $\beta(t)$ and sine $\gamma(t)$}
\quad

In this part, we use the TL-gPINNs to identify two unknown variable coefficients: linear $\beta(t)$ and sine $\gamma(t)$ when the other variable coefficient ($\alpha (t)=\sin(t)$) is fixed and the training dataset consisting of initial-boundary data $\{x^i_A,t^i_A,u^i,v^i\}^{N_A}_{i=1}$($N_A=200$) and internal data $\{x^i_{in},t^i_{in},u^i,v^i\}^{N_{A_{in}}}_{i=1}$(${N_{A_{in}}}=2000$) is randomly selected. Then the loss functions of PINNs and gPINNs are redefined 
\begin{equation}\label{PINN-two}
MSE_{inverse}=MSE_A+MSE_f+MSE_{A_{in}}+MSE_{\boldsymbol{\Lambda}},
\end{equation}
\begin{equation}\label{gPINN-two}
MSE_{inverse}^g=MSE_A+MSE_f+MSE_{A_{in}}+MSE_{\boldsymbol{\Lambda}}+MSE_{g},
\end{equation}
where
\begin{equation}
MSE_{\boldsymbol{\Lambda}}=MSE_{\beta}+MSE_{\gamma},	
\end{equation}
\begin{equation}
MSE_{\beta}=|\Widehat{{\beta}}(t_0)-{\beta}^0|^2,
\end{equation}
\begin{equation}
MSE_{\gamma}=\frac{1}{2}\left(|\Widehat{{\gamma}}(t_0)-{\gamma}^0|^2 + |\Widehat{{\gamma}}(t_1)-{\gamma}^1|^2 \right).
\end{equation}
The depth and width of neural networks for inferring the solution and variable coefficients are listed in Table \ref{table3-0}.

\begin{figure}[htbp]
\centering
\includegraphics[width=7.5cm,height=5cm]{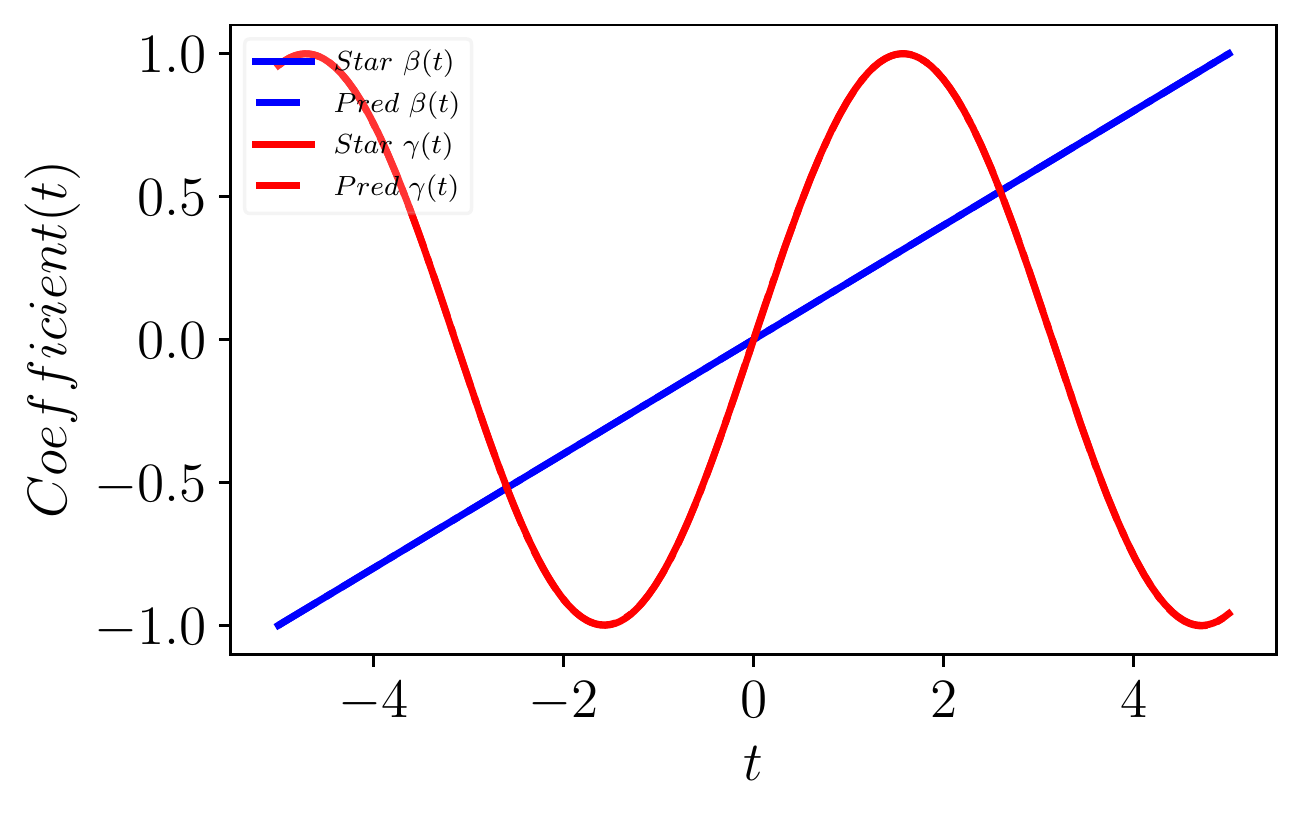}
$a$
\includegraphics[width=7.5cm,height=5cm]{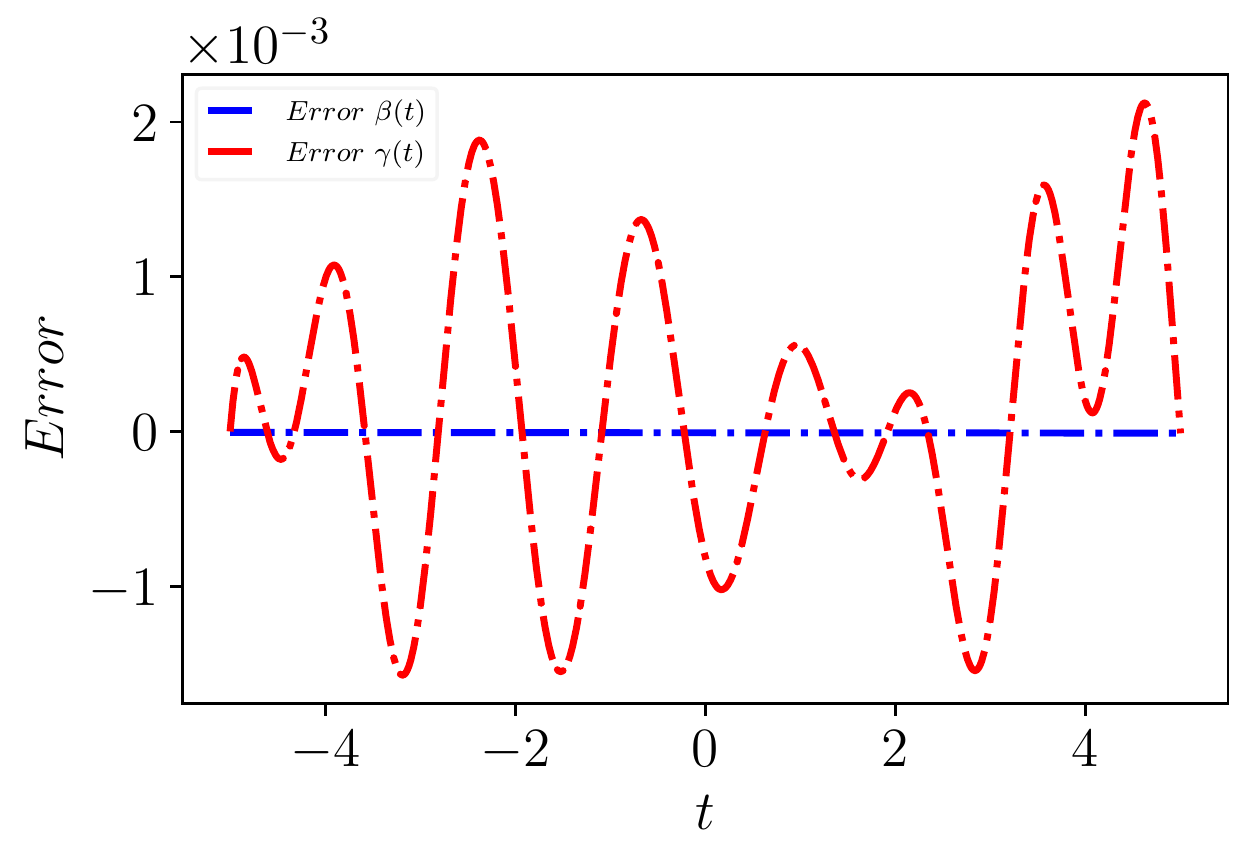}
$b$
\caption{(Color online) Results of function discovery for the vcNLS equation by TL-gPINNs: (a) Comparison between the predicted and exact variable coefficients $\beta(t)$ and $\gamma(t)$; (b) The absolute errors.}
\label{fig3-11}
\end{figure}

By employing the TL-gPINN method, the data-driven soliton solution and variable coefficients for the vcNLS equation are successfully simulated. Comparison between the predicted and exact variable coefficients $\beta(t)$ and $\gamma(t)$ as well as the corresponding absolute errors is displayed in Fig. \ref{fig3-11}. It can be seen that the absolute error of linear $\beta(t)$ is negligible compared with that of nonlinear $\gamma(t)$, which exhibits the feature of high-frequency oscillation due to the periodic oscillation and the change in concavity and convexity of the variable coefficient. Table \ref{table3-6} gives a brief overview of the method performance.

\begin{table}[htbp]
\caption{Performance comparison of three methods: mean absolute errors and relative $\mathbb{L}_2$ errors of the variable coefficients $\beta(t)$ and $\gamma(t)$ as well as error reduction rates.}
\label{table3-6}   
\centering
\begin{tabular}{cc|ccc}
\bottomrule
\multicolumn{2}{c|}{\multirow{2}{*}{Results}} & \multicolumn{3}{c}{Method}                                           \\ \cline{3-5} 
\multicolumn{2}{c|}{}                             & \multicolumn{1}{c|}{PINNs} & \multicolumn{1}{c|}{gPINNs} & TL-gPINNs \\ \hline
\multicolumn{1}{c|}{\multirow{2}{*}{$\beta(t)$}}   & $MAE_{\beta}$ ($ERR_1$)  & 1.246162e-05                           & 1.861258e-05 \textbf{(-49.36\%)}                                &  8.680616e-06 \textbf{(30.34\%)}    \\ \cline{2-2}
\multicolumn{1}{c|}{}                     & $RE_{\beta}$ ($ERR_2$)   & 2.323916e-05                           &  3.888068e-05 \textbf{(-67.31\%)}                            & 1.517259e-05 \textbf{(30.34\%)}      \\ \cline{1-2}
\multicolumn{1}{c|}{\multirow{2}{*}{$\gamma(t)$}}   & $MAE_{\gamma}$ ($ERR_1$)   & 1.412442e-03                           &  1.128851e-03 \textbf{(-67.31\%)}                            & 7.726441e-04 \textbf{(34.71\%)}      \\ \cline{2-2}
\multicolumn{1}{c|}{}                     & $RE_{\gamma}$ ($ERR_2$)   & 2.712897e-03                           &  2.220811e-03 \textbf{(18.14\%)}                             &   1.309694e-03 \textbf{(51.72\%)}   \\ \toprule
\end{tabular}
\end{table}

\subsubsection{Data-driven discovery of three variable coefficients}\label{3.2.2}
\quad

Note that all variable coefficients of the vcNLS equation are unknown here.

$\bullet$ \textbf{Linear $\alpha(t)$, $\beta(t)$ and $\gamma(t)$}

For the identification of three linear variable coefficients, the term embodying the training data in the loss functions in Eq. \eqref{PINN-two} and \eqref{gPINN-two} need to be modified
\begin{equation}
MSE_{\boldsymbol{\Lambda}}=MSE_{\alpha}+MSE_{\beta}+MSE_{\gamma},	
\end{equation}
\begin{equation}
MSE_{\alpha}=|\Widehat{{\alpha}}(t_0)-{\alpha}^0|^2,
\end{equation}
\begin{equation}
MSE_{\beta}=|\Widehat{{\beta}}(t_0)-{\beta}^0|^2,
\end{equation}
\begin{equation}
MSE_{\gamma}=|\Widehat{{\gamma}}(t_0)-{\gamma}^0|^2.
\end{equation}

With the aid of the same generation and sampling method above, we obtain the training data (size: $N_A=200, N_{A_{in}}=2000$) in the given spatiotemporal region $[-4,4]\times[-4,4]$, where the corresponding soliton solution is 
\begin{align}
A(x,t)=\frac{\mathrm{e}^{\frac{{\rm{i}}}{10} t^2} \mathrm{e}^{(1+{\rm{i}})x -2 \ln(\cosh(t))}}{1+\frac{\mathrm{e}^{2x-4\ln(\cosh(t))}}{8}}.	
\end{align}
The linear and $\tanh$ activation functions are adopted to infer the variable coefficients and soliton solution separately.

\begin{figure}[htbp]
\centering
\includegraphics[width=7.5cm,height=5cm]{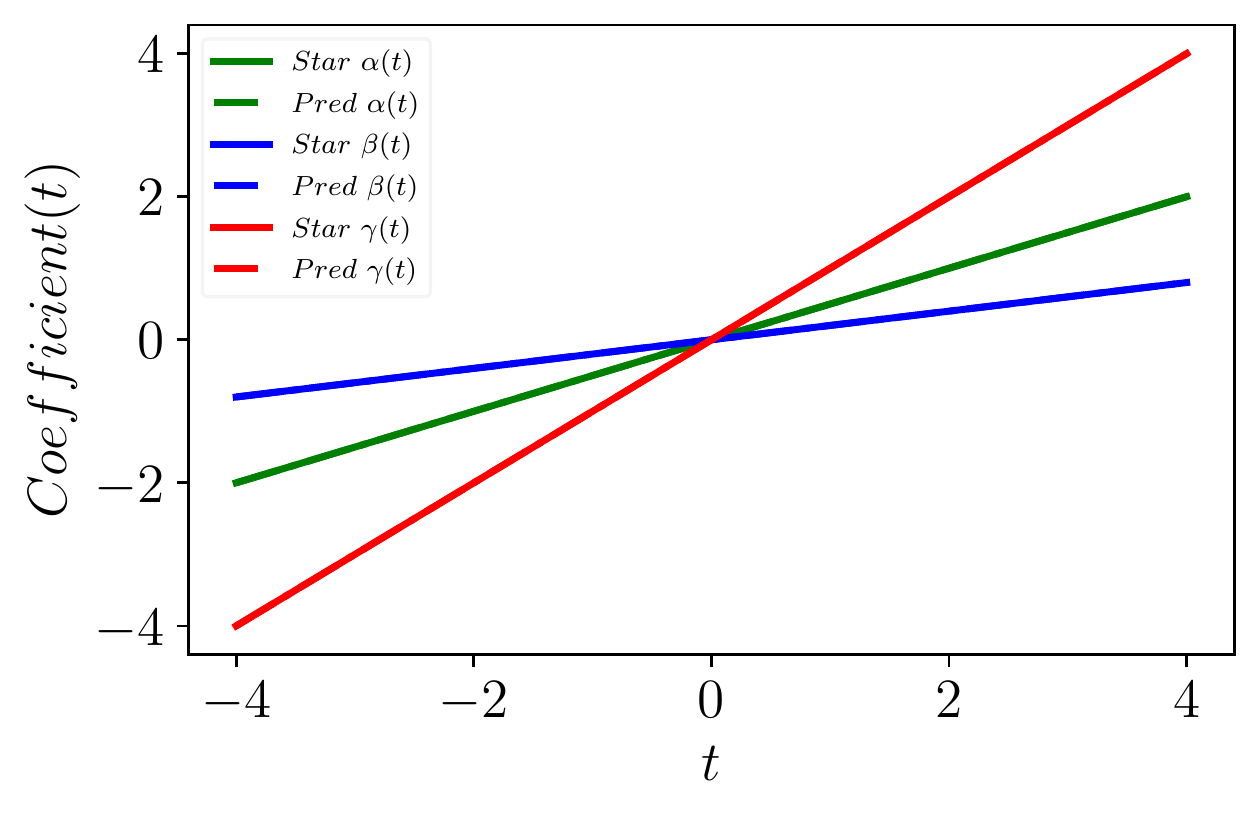}
$a$
\includegraphics[width=7.5cm,height=5cm]{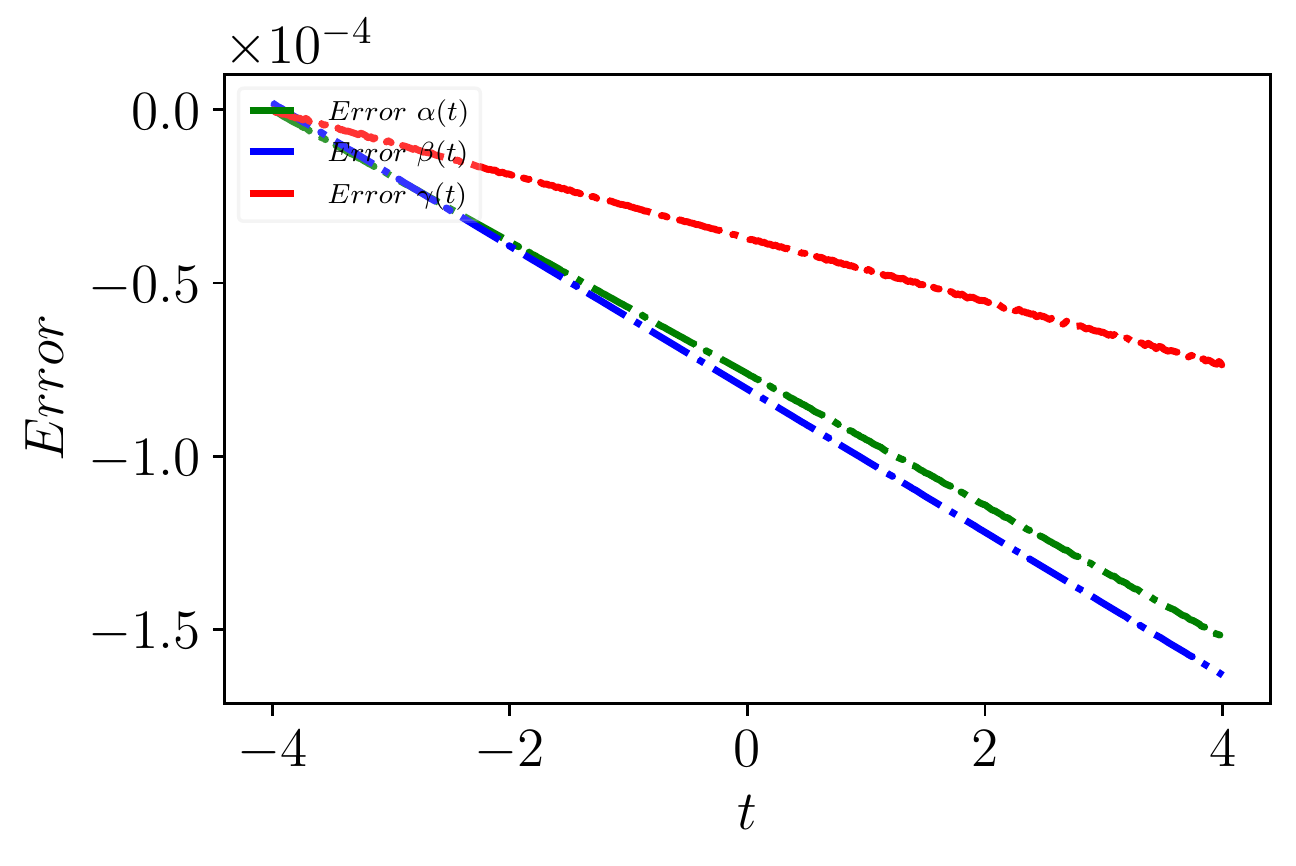}
$b$
\caption{(Color online) Results of function discovery for the vcNLS equation by TL-gPINNs: (a) Comparison between the predicted and exact variable coefficients $\alpha(t)$, $\beta(t)$ and $\gamma(t)$; (b) The absolute errors.}
\label{fig3-12}
\end{figure}

Finally, Fig. \ref{fig3-12} shows the curve plots of the predicted and the exact variable coefficients as well as absolute errors obtained by TL-gPINN, and Table \ref{table3-7} summarizes the detailed results of three methods in the term of prediction accuracy. The change of absolute error curves here is similar to that in Sec. \ref{3.1}.

\begin{table}[htbp]
\caption{Performance comparison of three methods:  mean absolute errors and relative $\mathbb{L}_2$ errors of three linear variable coefficients as well as error reduction rates.}
\label{table3-7}  
\centering
\begin{tabular}{cc|ccc}
\bottomrule
\multicolumn{2}{c|}{\multirow{2}{*}{Results}} & \multicolumn{3}{c}{Method}                                           \\ \cline{3-5} 
\multicolumn{2}{c|}{}                             & \multicolumn{1}{c|}{PINNs} & \multicolumn{1}{c|}{gPINNs} & TL-gPINNs \\ \hline
\multicolumn{1}{c|}{\multirow{2}{*}{$\alpha(t)$}}   & $MAE_{\alpha}$ ($ERR_1$)  & 7.617165e-05                           & 1.224185e-04 \textbf{(-60.71\%)}                                &  6.102315e-05 \textbf{(19.89\%)}   \\ \cline{2-2}
\multicolumn{1}{c|}{}                     & $RE_{\alpha}$ ($ERR_2$)   &  7.604356e-05                          &  1.225268e-04 \textbf{(-61.13\%)}                            &  6.186360e-05 \textbf{(18.65\%)}   \\ \cline{1-2}
\multicolumn{1}{c|}{\multirow{2}{*}{$\beta(t)$}}   & $MAE_{\beta}$ ($ERR_1$)  & 2.193591e-04                           & 5.647830e-04 \textbf{(-157.47\%)}                              &  8.066458e-05 \textbf{(63.23\%)}      \\ \cline{2-2}
\multicolumn{1}{c|}{}                     & $RE_{\beta}$ ($ERR_2$)   & 5.484737e-04                           & 1.412454e-03 \textbf{(-157.52\%)}                               &  2.025123e-04 \textbf{(63.08\%)}    \\ \cline{1-2}
\multicolumn{1}{c|}{\multirow{2}{*}{$\gamma(t)$}}   & $MAE_{\gamma}$ ($ERR_1$)   & 1.086061e-04                            & 4.343718e-04 \textbf{(-299.95\%)}                                 &  3.701118e-05 \textbf{(65.92\%)}     \\ \cline{2-2}
\multicolumn{1}{c|}{}                     & $RE_{\gamma}$ ($ERR_2$)   & 5.447867e-05                           & 2.165039e-04 \textbf{(-297.41\%)}                               &  1.842847e-05 \textbf{(66.17\%)}    \\ \toprule
\end{tabular}
\end{table}

$\bullet$ \textbf{Linear $\beta(t)$, fractional $\alpha(t)$ and $\gamma(t)$}

Based on the initial-boundary data of the soliton solution
\begin{align}
A(x,t)=\frac{\mathrm{e}^{\frac{{\rm{i}}}{10} t^2 }\mathrm{e}^{(1+{\rm{i}})x -\arctan(t)}}{1+\frac{(2t^2+2) \mathrm{e}^{2x-2 \arctan(t)}}{8(t^2+1)}},
\end{align}
corresponding to $\alpha(t)=\frac{1}{2(1+t^2)},\beta(t)=\frac{t}{5},\gamma(t)=\frac{1}{1+t^2}$, we utilize the TL-gPINNs to infer these three unknown variable coefficients. Here, the loss term of nonlinear variable coefficients should be changed into
\begin{equation}
MSE_{vc}=MSE_{\alpha}+MSE_{\beta}+MSE_{\gamma},	
\end{equation}
where
\begin{equation}
MSE_{\beta}=|\Widehat{{\beta}}(t_0)-{\beta}^0|^2,
\end{equation}
\begin{equation}
MSE_{\alpha}=\frac{1}{2}\left(|\Widehat{{\alpha}}(t_0)-{\alpha}^0|^2 + |\Widehat{{\alpha}}(t_1)-{\alpha}^1|^2 \right),
\end{equation}
\begin{equation}
MSE_{\gamma}=\frac{1}{2}\left(|\Widehat{{\gamma}}(t_0)-{\gamma}^0|^2 + |\Widehat{{\gamma}}(t_1)-{\gamma}^1|^2 \right).
\end{equation}

\begin{figure}[htbp]
\centering
\includegraphics[width=7.5cm,height=5cm]{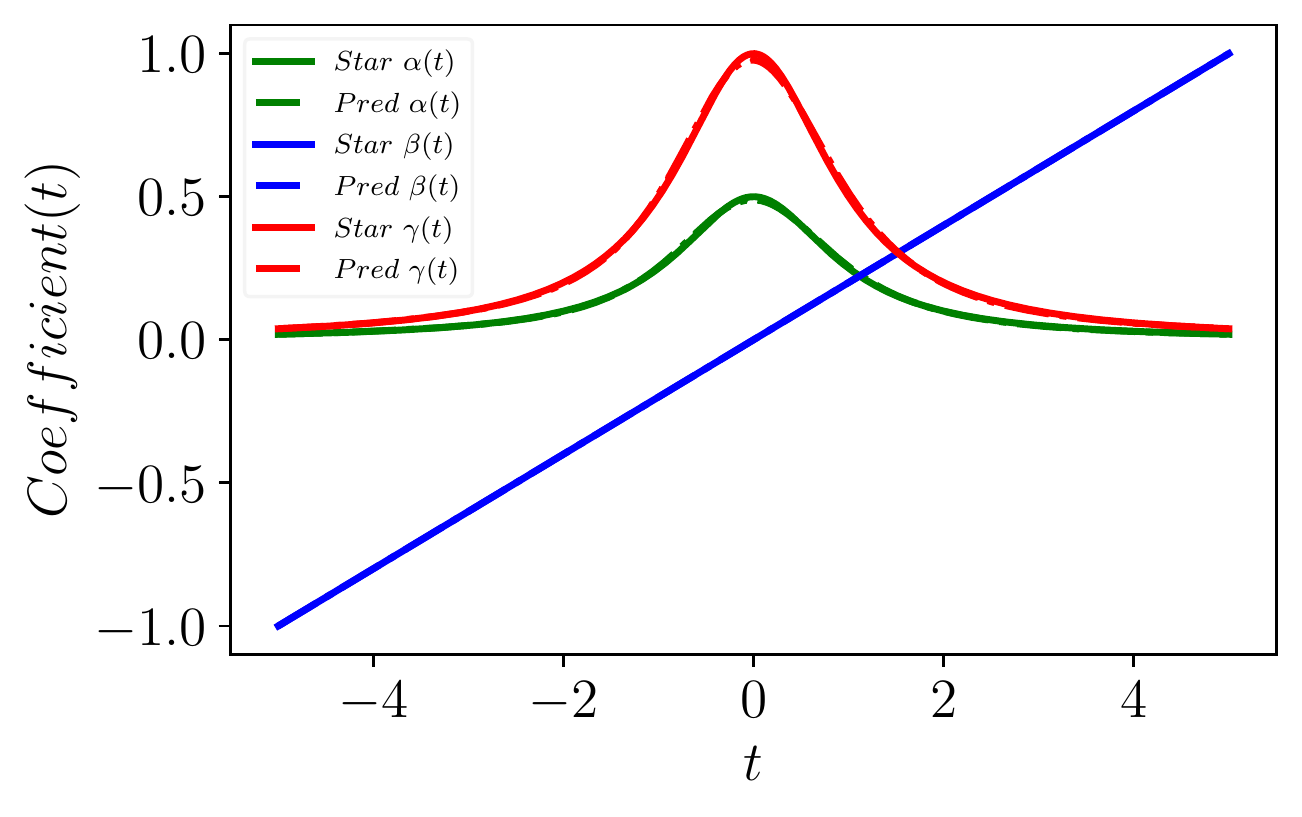}
$a$
\includegraphics[width=7.5cm,height=5cm]{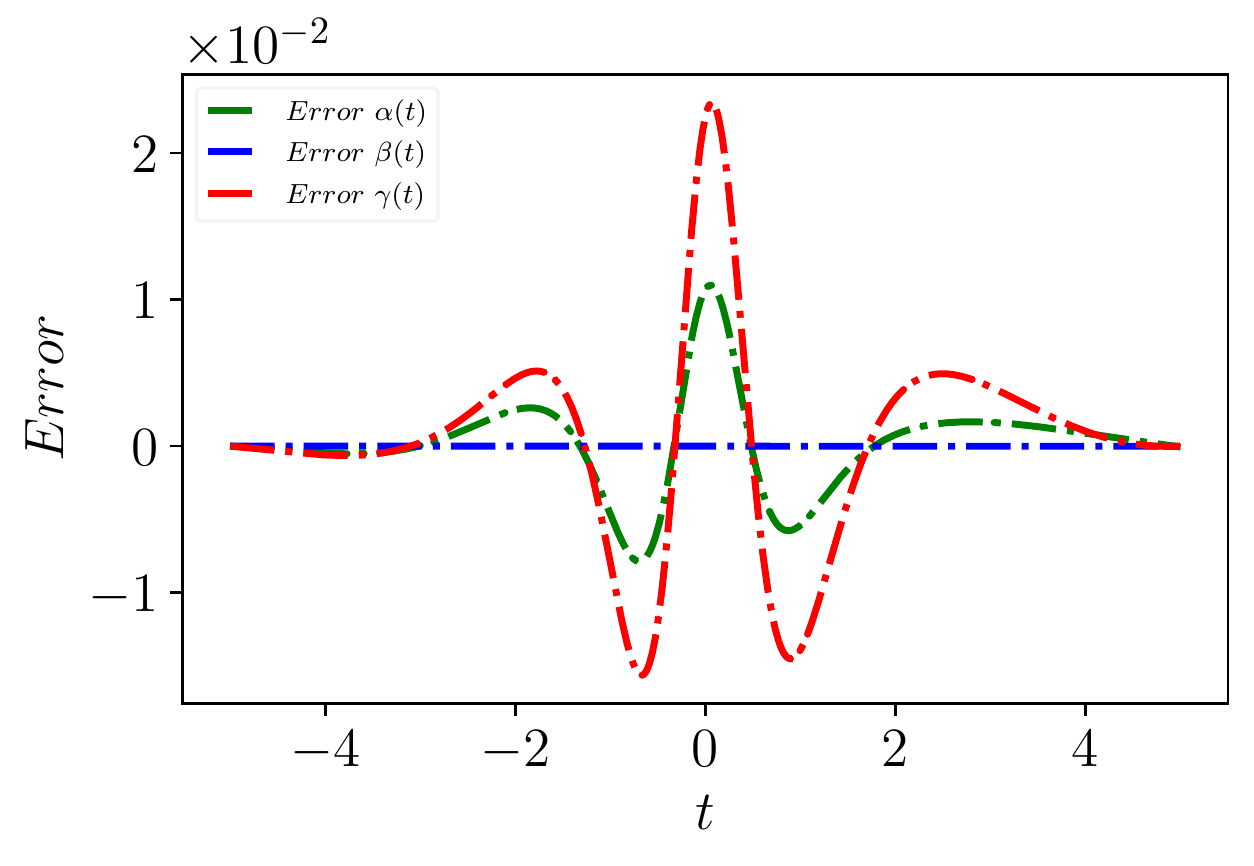}
$b$
\caption{(Color online) Results of function discovery for  the vcNLS equation by TL-gPINNs: (a) Comparison between the predicted and exact variable coefficients $\alpha(t)$, $\beta(t)$ and $\gamma(t)$; (b) The absolute errors.}
\label{fig3-13}
\end{figure}

Results of function discovery for the vcNLS equation, i.e. comparisons between the predicted and exact variable coefficients $\alpha(t)$, $\beta(t)$ and $\gamma(t)$ as well as their respective absolute errors are presented in Fig. \ref{fig3-13}. Similarly, the absolute error of linear variable coefficient $\beta(t)$ is negligible compared with those of nonlinear ones. The variable coefficients $\alpha(t)$ and $\gamma(t)$ in the form of fractional polynomials also basically meets the rule of thumb mentioned in Sec. \ref{3.1.5}. However, the phenomenon of multiple intersections between error curves and more specific feature analysis remain to be further explored in future work. Besides, the running time of PINNs, TL-gPINNs and gPINNs is: 242.0248, 320.302 and 892.13 seconds, respectively. The performance comparison of these three methods is shown in Table \ref{table3-8}.

\begin{table}[htbp]
\caption{Performance comparison of three methods: mean absolute errors and relative $\mathbb{L}_2$ errors of three variable coefficients as well as error reduction rates.}
\label{table3-8}  
\centering
\begin{tabular}{cc|ccc}
\bottomrule
\multicolumn{2}{c|}{\multirow{2}{*}{Results}} & \multicolumn{3}{c}{Method}                                           \\ \cline{3-5} 
\multicolumn{2}{c|}{}                             & \multicolumn{1}{c|}{PINNs} & \multicolumn{1}{c|}{gPINNs} & TL-gPINNs \\ \hline
\multicolumn{1}{c|}{\multirow{2}{*}{$\alpha(t)$}}   & $MAE_{\alpha}$($ERR_1$)  & 2.536442e-03                           & 2.372040e-03 \textbf{(6.48\%)}                               & 2.069355e-03 \textbf{(18.42\%)}    \\ \cline{2-2}
\multicolumn{1}{c|}{}                     & $RE_{\alpha}$($ERR_2$)   & 1.855795e-02                           & 1.840729e-02 \textbf{(0.81\%)}                             & 1.616184e-02 \textbf{(12.91\%)}      \\ \cline{1-2}
\multicolumn{1}{c|}{\multirow{2}{*}{$\beta(t)$}}   & $MAE_{\beta}$($ERR_1$)  & 3.263991e-05                           & 2.116216e-05 \textbf{(35.16\%)}                                & 9.873565e-06 \textbf{(69.75\%)}    \\ \cline{2-2}
\multicolumn{1}{c|}{}                     & $RE_{\beta}$($ERR_2$)   &5.729922e-05                            &3.738890e-05 \textbf{(34.75\%)}                                & 1.751263e-05 \textbf{(69.44\%)}      \\ \cline{1-2}
\multicolumn{1}{c|}{\multirow{2}{*}{$\gamma(t)$}}   & $MAE_{\gamma}$($ERR_1$)   & 5.610852e-03                           & 5.163366e-03 \textbf{(7.98\%)}                               &  4.453668e-03 \textbf{(20.62\%)}     \\ \cline{2-2}
\multicolumn{1}{c|}{}                     & $RE_{\gamma}$($ERR_2$)   & 2.201791e-02                           & 2.034598e-02 \textbf{(7.59\%)}                              &  1.753237e-02 \textbf{(20.37\%)}    \\ \toprule
\end{tabular}
\end{table}

\subsection{Result analysis}
\quad

According to the performance comparison of three methods (PINNs, TL-gPINNs and gPINNs) presented in Table. \ref{table3-1} - Table. \ref{table3-8},  TL-gPINNs possess the notable performance of high accuracy whether in identifying single variable coefficient or in inferring multiple ones compared with the other two methods. Meanwhile, TL-gPINNs can accelerate convergence of iteration and reduce calculation time since the technique of transfer learning helps to mitigate the problem of inefficiency caused by extra loss terms of the gradient.

The reason why gPINN doesn't perform up to expectations here may be that the solution $A(x,t)$ for the variable coefficient nonlinear Schr\"{o}odinger equation is complex-valued and each constraint function in neural networks should be decomposed into two parts: the real and imaginary parts. Thus, the loss function itself consists of many constraint terms even without regard to the gradient restriction. When solving the multi-objective optimization problems, the local optimum that it ultimately converges to is obtained based on the competitive relationship between various objectives. Therefore, the result may not necessarily be better even if more constraints are imposed. Evidently, the experiments show that gPINN has lower prediction accuracy than PINN even at the cost of sacrificing efficiency especially in Case \ref{3.1.2} (shown in Table \ref{table3-2}), Case \ref{3.1.5} (Table \ref{table3-5}) and Case \ref{3.2.2} (Table \ref{table3-7}). The advantage of the TL-gPINN method lies in that gPINN inherits the saved weight matrixes and bias vectors of PINN at the end of the iteration process as the initialization parameters, and thus the subsequent training of gPINN is based on that of PINN by leveraging the transfer learning technique instead of training from scratch. Consequently, TL-gPINN is steadier on precision promotion compared to gPINN, a method which has been proved to be efficient in improving the accuracy of PINN \cite{gPINN}.

\begin{figure}[htbp]
\centering
\includegraphics[width=8cm,height=5.5cm]{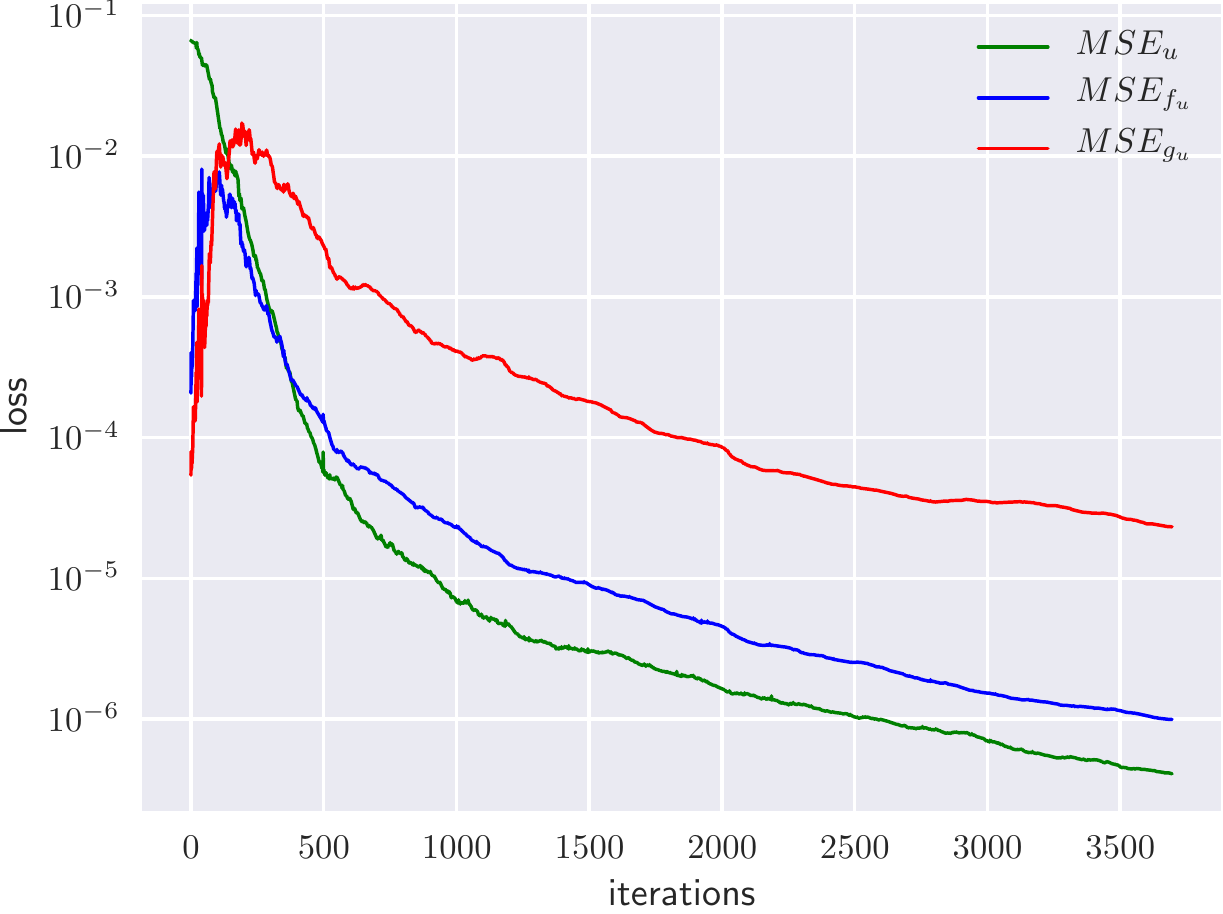}
$a$
\includegraphics[width=8cm,height=5.5cm]{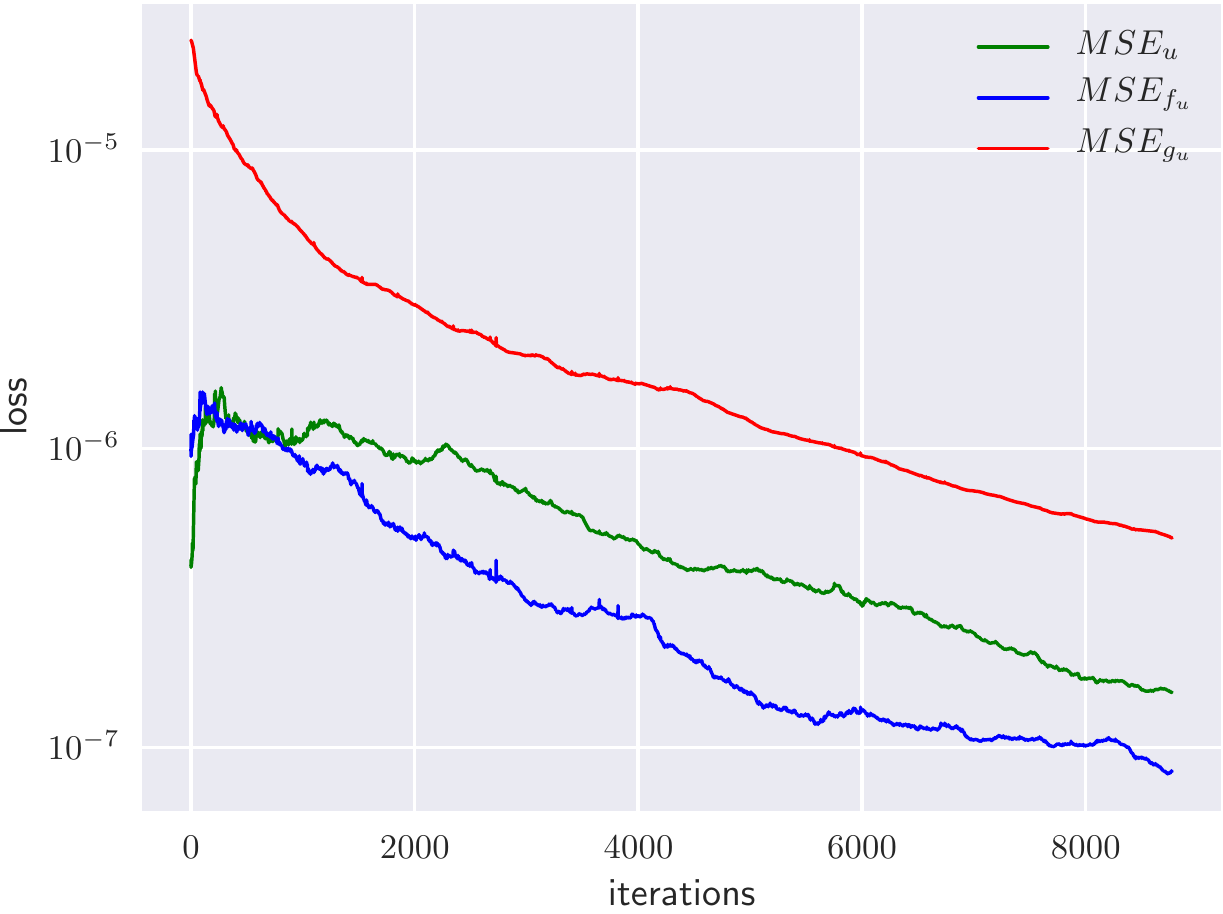}
$b$
\caption{(Color online) Evolution of the loss functions in inferring linear variable coefficient $\gamma(t)$ for the vcNLS equation by two methods: (a) PINN; (b) TL-gPINN.}
\label{fig3-14}
\end{figure}

What's more, the loss curve figures of inferring linear variable coefficient $\gamma(t)$ in Sec. \ref{3.1} are plotted in Fig. \ref{fig3-14} for the sake of more intuitive analysis. Here, only loss functions corresponding to the real part, i.e. $MSE_u, MSE_{f_u}$ and $MSE_{g_u}$, are considered and counterparts of the imaginary part ($MSE_v, MSE_{f_v}$ and $MSE_{g_v}$) change approximately in the same way. The values of each loss term at the beginning and end of iterations are listed in Table \ref{table3-9}. 

\begin{table}[htbp]
\caption{Results of losses at the beginning and end of iteration in inferring linear variable coefficient $\gamma(t)$ for the vcNLS equation by three methods.}
\label{table3-9}  
\centering
\begin{tabular}{c|cccccc}
\bottomrule
\multirow{3}{*}{\diagbox{\textbf{\textbf{Results}}}{\textbf{Method}}} & \multicolumn{2}{c|}{PINNs}                                                                                                                                              & \multicolumn{2}{c|}{gPINNs}                                                                                                                                          & \multicolumn{2}{c}{TL-gPINNs}                                                                                                                         \\ \cline{2-7} 
& \multicolumn{1}{c|}{\begin{tabular}[c]{@{}c@{}}The zeroth\\ iteration\end{tabular}} & \multicolumn{1}{c|}{\begin{tabular}[c]{@{}c@{}}The last\\ iteration\end{tabular}} & \multicolumn{1}{c|}{\begin{tabular}[c]{@{}c@{}}The zeroth\\ iteration\end{tabular}} & \multicolumn{1}{c|}{\begin{tabular}[c]{@{}c@{}}The last\\ iteration\end{tabular}} & \multicolumn{1}{c|}{\begin{tabular}[c]{@{}c@{}}The zeroth\\ iteration\end{tabular}} & \begin{tabular}[c]{@{}c@{}}The last\\ iteration\end{tabular} \\ \hline
$MSE_u$  & 6.672856e-02                                                                                    &  4.078804e-07                                                                                   &  6.672856e-02                                                                                   &  1.650002e-07 &  4.078804e-07                                                                                                                                                  &  1.530691e-07 \\
$MSE_v$  & 8.681140e-02                                                                                   & 4.347213e-07                                                                                  & 8.681140e-02                                                                                     & 2.090931e-07 & 4.347213e-07                                                                                                                                                  & 2.113367e-07    \\
$MSE_{f_u}$  & 2.073652e-04                                                                                   & 9.956128e-07                                                                                   & 2.073652e-04                                                                                     & 1.071937e-07 & 9.956128e-07                                                                                                                                                 & 8.313516e-08  \\
$MSE_{f_v}$  & 4.143952e-04                                                                                    &  1.020429e-06                                                                                 &  4.143952e-04                                                                                   &  1.637525e-07 &  1.020429e-06                                                                                                                                                  &  1.001443e-07   \\
$MSE_{g_u}$  & 5.406206e-05                                                                                    & 2.326252e-05                                                                                  & 5.406206e-05                                                                                     & 6.869229e-07 & 2.326252e-05                                                                                                                                                    & 5.017756e-07  \\
$MSE_{g_v}$ & 8.511276e-05                                                                                    &  3.057861e-05                                                                                 &  8.511276e-05                                                                                    &  7.017902e-07 &  3.057861e-05                                                                                                                                                  &  4.487380e-07   \\
$MSE_{\gamma}$  & 16.521566                                                                                    & 3.637979e-12                                                                                  & 16.521566                                                                                     & 2.273737e-13 & 3.637979e-12                                                                                                                                                  & 2.273737e-13   \\ \toprule
\end{tabular}
\end{table}

As we can see from Fig. \ref{fig3-14} (a),   $MSE_{g_u}$ fluctuated at a relatively high level and the value of the last iteration is almost the same as that at the beginning of the iteration while $MSE_u$ and $MSE_{f_u}$ decreased to 4.078804e-07 and 9.956128e-07 respectively during the training process of PINNs, where the loss of gradient $MSE_{g_u}$ has no contribution to optimization. In Fig. \ref{fig3-14} (b) and Table \ref{table3-9}, it is obvious to note that the values of loss terms of the gradients (i.e.,$MSE_{g_u}$ and $MSE_{g_v}$) are larger by several orders of magnitude than those of other loss terms in the zeroth iteration when the weight transfer is just completed. Specifically, the values of $MSE_u, MSE_v, MSE_{f_u}, MSE_{f_v}, MSE_{u_{in}}, MSE_{v_{in}}$ are approximately remain between 10e-07 and 10e-06, and that of $MSE_{\gamma}$ maintains at 10e-11 to 10e-10 while the values of $MSE_{g_u}$ and $MSE_{g_v}$ are at a high level of 10e-5 to 10e-4. It reveals that there is still some deviation between the variable coefficients themselves and the ones learned by the PINN method from the perspective of gradients. In other words, the PINN method lacks sufficient attention to gradients and leads to inadequate optimization, which may be an underlying cause why the training of gPINNs can go on effectively after finishing the weight transfer. Then the values of $MSE_{g_u}$ dropped pretty steadily while $MSE_u$ and $MSE_{f_u}$ showed a downward trend after an initial ascent. Meanwhile, the process of their ascent happens to be that of the fastest descent of $MSE_{g_u}$, and we deduced that it may be a process of escaping from the local optimal point obtained by PINN, where the values of gradient loss are large although those of other loss terms are at a fairly low level.

With regard to efficiency, gPINNs significantly increase the time cost of training due to the introduction of additional gradient constraints while TL-gPINNs shorten the training time in contrast to the original gPINNs by taking full advantage of transfer learning.

In short, the TL-gPINN method achieves the highest prediction accuracy among the three methods whether in inferring unknown single variable coefficient or in identifying multiple ones. However, gPINN shows an unstable performance here and even performs no better than PINN in accuracy in some cases. For TL-gPINNs, the application of transfer learning technique can contribute to both higher efficiency and greater reliability than the original PINN. It outperforms the PINNs in accuracy and gPINNs in both accuracy and efficiency. Thereupon the TL-gPINN method is more superior compared with the PINN and gPINN here.

\section{Analysis and discussion}\label{analysis}

\subsection{Robustness analysis}\label{Robustness}
\quad

Numerical results presented in Sec. \ref{NLS} are based on noise-free training data, and here we carry out experiments when the training data was corrupted with noise to test the robustness of the TL-gPINNs. 

Specifically, the training data, including the initial-boundary data $\{x^i_A,t^i_A,u^i,v^i\}^{N_A}_{i=1}$,  internal data $\{x^i_{in},t^i_{in},u^i,v^i\}^{N_{A_{in}}}_{i=1}$ and the data $\{t^i_{{\gamma}},{\gamma}^i\}^{N_{\gamma}}_{i=1}$ of the variable coefficient $\gamma(t)$, is corrupted by four different noise levels: 0.5\%, 1\%, 3\% and 5\%.

\begin{table}[htbp]
\caption{Performance comparison of three methods in identifying variable coefficient $\gamma(t)$ for the vcNLS equation under different noise conditions.}
\label{table4-1}  
\centering
\begin{tabular}{ccc|ccccc}
\bottomrule
\multicolumn{3}{c|}{\multirow{2}{*}{Results}}                                                   & \multicolumn{5}{c}{Correct $\gamma(t)$} \\ \cline{4-8} 
\multicolumn{3}{c|}{}                                                                               & $t$   & $t^2$   & $\sin(t)$   & $\tanh(t)$   & $\frac{1}{1+t^2}$   \\ \hline
\multicolumn{1}{c|}{\multirow{6}{*}{0.5\% noise}} & \multicolumn{2}{c|}{$MAE_{\gamma}$}                              & 1.322670e-04    &  9.615011e-03   & 1.097929e-03 & 3.973503e-03    & 1.009619e-03    \\ \cline{2-3}
\multicolumn{1}{c|}{}                       & \multicolumn{2}{c|}{$RE_{\gamma}$}                               & 6.551441e-05    & 7.959068e-03    & 2.579889e-03    & 5.486678e-03    & 2.968212e-03    \\ \cline{2-3}
\multicolumn{1}{c|}{}                       & \multicolumn{1}{c|}{\multirow{2}{*}{$ERR_1$}} & TL-gPINNs & 0.00\%    & 3.16\%    & 49.86\%    & 51.81\%    & 5.06\%    \\ \cline{3-3}
\multicolumn{1}{c|}{}                       & \multicolumn{1}{c|}{}                     & gPINNs    & -73.90\%    & -90.00\%    & -26.63\%    & 39.10\%    & -10.64\%    \\ \cline{2-3}
\multicolumn{1}{c|}{}                       & \multicolumn{1}{c|}{\multirow{2}{*}{$ERR_2$}} & TL-gPINNs & 0.00\%    & 1.55\%    &39.72\%     &50.65\%     & 10.76\%    \\ \cline{3-3}
\multicolumn{1}{c|}{}                       & \multicolumn{1}{c|}{}                     & gPINNs    & -75.55\%    & -61.89\%    & -35.03\%    &39.97\%     & -16.03\%  \\ \hline
\multicolumn{1}{c|}{\multirow{6}{*}{1\% noise}} & \multicolumn{2}{c|}{$MAE_{\gamma}$}                              &3.393921e-04     & 1.990454e-02    & 2.040646e-03    & 5.745897e-03    & 1.529080e-03    \\ \cline{2-3}
\multicolumn{1}{c|}{}                       & \multicolumn{2}{c|}{$RE_{\gamma}$}                               &1.693800e-04     & 2.038240e-02     & 6.616817e-03    & 7.922166e-03    &  4.553724e-03   \\ \cline{2-3}
\multicolumn{1}{c|}{}                       & \multicolumn{1}{c|}{\multirow{2}{*}{$ERR_1$}} & TL-gPINNs &26.92\%     &3.80\%     & 28.23\%    & 9.67\%    &39.17\%     \\ \cline{3-3}
\multicolumn{1}{c|}{}                       & \multicolumn{1}{c|}{}                     & gPINNs    & 47.24\%    & 18.84\%    & -66.55\%    & 4.16\%    & 29.09\%    \\ \cline{2-3}
\multicolumn{1}{c|}{}                       & \multicolumn{1}{c|}{\multirow{2}{*}{$ERR_2$}} & TL-gPINNs &26.93\%     & -1.78\%    & 10.75\%    & 10.54\%    &42.98\%    \\ \cline{3-3}
\multicolumn{1}{c|}{}                       & \multicolumn{1}{c|}{}                     & gPINNs    & 47.31\%    & 8.33\%    & -37.25\%    & 2.59\%    & 36.48\%    \\ \hline
\multicolumn{1}{c|}{\multirow{6}{*}{3\% noise}} & \multicolumn{2}{c|}{$MAE_{\gamma}$}                              & 3.642581e-04    & 1.241316e-02    & 3.276918e-03    & 1.005220e-02    & 2.160257e-03    \\ \cline{2-3}
\multicolumn{1}{c|}{}                       & \multicolumn{2}{c|}{$RE_{\gamma}$}                               & 1.825354e-04    & 1.403889e-02    & 8.239684e-03    & 1.449078e-02    & 6.147466e-03    \\ \cline{2-3}
\multicolumn{1}{c|}{}                       & \multicolumn{1}{c|}{\multirow{2}{*}{$ERR_1$}} & TL-gPINNs & 48.18\%    &14.93\%     & 9.34\%    & 9.38\%    & 13.12\%    \\ \cline{3-3}
\multicolumn{1}{c|}{}                       & \multicolumn{1}{c|}{}                     & gPINNs    & 27.43\%    &-0.74\%     &-12.80\%     & 11.36\%    & 7.65\%    \\ \cline{2-3}
\multicolumn{1}{c|}{}                       & \multicolumn{1}{c|}{\multirow{2}{*}{$ERR_2$}} & TL-gPINNs & 48.11\%    & 0.09\%    & -2.99\%    &8.76\%     &33.17\%     \\ \cline{3-3}
\multicolumn{1}{c|}{}                       & \multicolumn{1}{c|}{}                     & gPINNs    & 27.17\%    & -2.20\%    & -5.96\%    & 10.90\%    & 1.89\%    \\ \hline
\multicolumn{1}{c|}{\multirow{6}{*}{5\% noise}} & \multicolumn{2}{c|}{$MAE_{\gamma}$}                              & 3.007159e-04    & 3.890978e-02    & 5.462772e-03    & 8.254448e-03    & 2.393983e-03   \\ \cline{2-3}
\multicolumn{1}{c|}{}                       & \multicolumn{2}{c|}{$RE_{\gamma}$}                               & 1.475556e-04    & 3.750857e-02    & 1.526540e-02    &  1.280577e-02   & 7.430092e-03    \\ \cline{2-3}
\multicolumn{1}{c|}{}                       & \multicolumn{1}{c|}{\multirow{2}{*}{$ERR_1$}} & TL-gPINNs & 32.08\%    & 5.61\%    & 25.55\%    & 11.41\%    & 26.18\%    \\ \cline{3-3}
\multicolumn{1}{c|}{}                       & \multicolumn{1}{c|}{}                     & gPINNs    & 23.02\%    &6.66\%     & 20.13\%    &-24.42\%     & -27.61\%    \\ \cline{2-3}
\multicolumn{1}{c|}{}                       & \multicolumn{1}{c|}{\multirow{2}{*}{$ERR_2$}} & TL-gPINNs & 32.65\%    & 0.68\%    & 9.56\%    &9.59\%     &37.61\%     \\ \cline{3-3}
\multicolumn{1}{c|}{}                       & \multicolumn{1}{c|}{}                     & gPINNs    &23.34\%     & -1.64\%    &7.05\%     &-23.97\%     & -7.41\%    \\ \toprule
\end{tabular}
\end{table}

Table \ref{table4-1} summarizes the results of the numerical experiments in the conditions of different noise levels and the indexes $MAE_{\gamma}$ and $RE_{\gamma}$ listed here are achieved by TL-gPINNs. The detailed results of PINNs and gPINNs are not provided here but shown in Table \ref{tableA-0} in Appendix A due to length limitations. Here, the reason why 0.00\% appears is that TL-gPINNs  converge rapidly after merely a few iterations, which means the local optimum obtained by PINNs also belongs to TL-gPINNs and then the training will not continue after initialization with saved weight data of PINNs.

According to the mean absolute error ($MAE_{\gamma}$) and relative $\mathbb{L}_2$ error ($RE_{\gamma}$) achieved by TL-gPINNs, different types of the variable coefficients $\gamma(t)$ can be identified accurately via this method. Evidently, the predictions of the unknown variable coefficient retain good robustness even when the training data was corrupted with different levels of noise. It also turns out that the accuracy of TL-gPINNs doesn't necessarily become worse with the increase of noise intensity, but may also increase in some cases.

Since the values of $ERR_1$ and $ERR_2$ indicate the degree of prediction accuracy improvement in the sense of the mean absolute error ($MAE_{\gamma}$) and relative $\mathbb{L}_2$ error ($RE_{\gamma}$) respectively, the results demonstrated that the ability of TL-gPINNs in precision promotion also remains robust to noise. We observe that the vast majority of experiments by TL-gPINNs have better performance than that of gPINNs in enhancing the accuracy of inferring the unknown variable coefficient and improving the generalization capability after assessing and comparing $ERR_1$ and $ERR_2$ of these two methods. Meanwhile, the higher efficiency of TL-gPINNs compared with the original gPINNs is a distinct advantage as well.

Based on the performance in Sec. \ref{NLS} and Sec. \ref{Robustness}, regardless of whether the training data is corrupted with noise or not, TL-gPINNs possess the ability to successfully infer the unknown variable coefficient $\gamma(t)$ with satisfactory accuracy. Taken overall, the TL-gPINNs meet the robustness and computational accuracy standards required in practice.

\subsection{Parametric sensitivity analysis}
\quad

The training results of neural networks are influenced by many factors, such as the architecture of neural networks and the size of training dataset. Thus, the parametric sensitivity analysis is conducted here to disclose the effect of these hyper-parameters on predictions of the single nonlinear variable coefficient $\gamma(t)$.

$\bullet$ \textbf{The architecture of neural networks}

With regard to the structure of fully-connected neural networks (FNN), the emphasis is put on the number of weighted layers (depth) and the number of neurons per hidden layer (width). Then we explore how the change of width and depth of the branch network for inferring the variable coefficient will affect the experimental results.

Meanwhile, we mainly investigate nonlinear variable coefficient $\gamma(t)$ mentioned in Sec. \ref{NLS}, which is more common in practice. For each form of the unknown nonlinear variable coefficient $\gamma(t)$, two hyper-parameters are changed: depth from 4 to 5 and width from 10 to 50 with step size 10.

\begin{figure}[htbp]
\centering
\includegraphics[width=5.5cm,height=3.5cm]{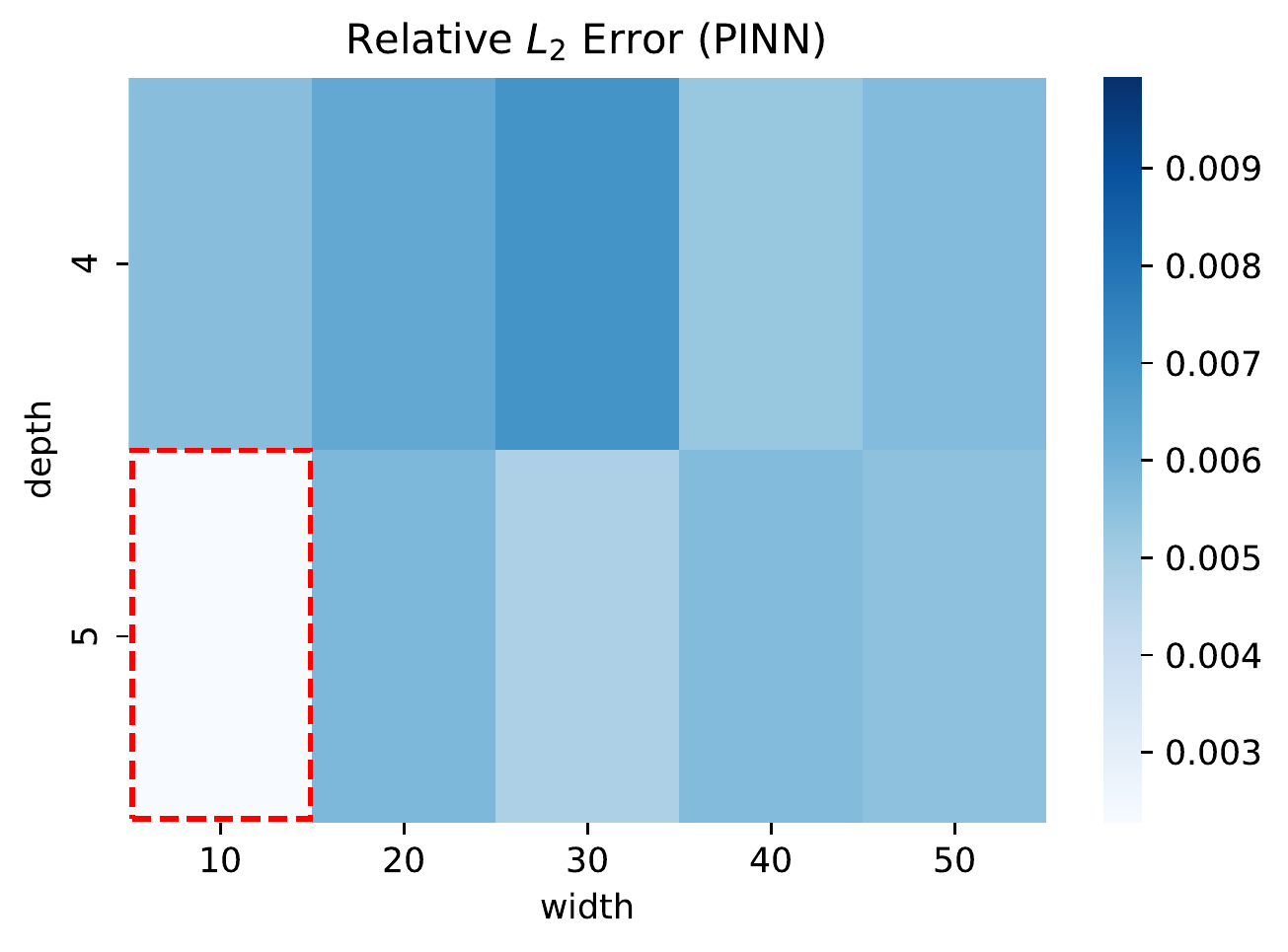}
\includegraphics[width=5.5cm,height=3.5cm]{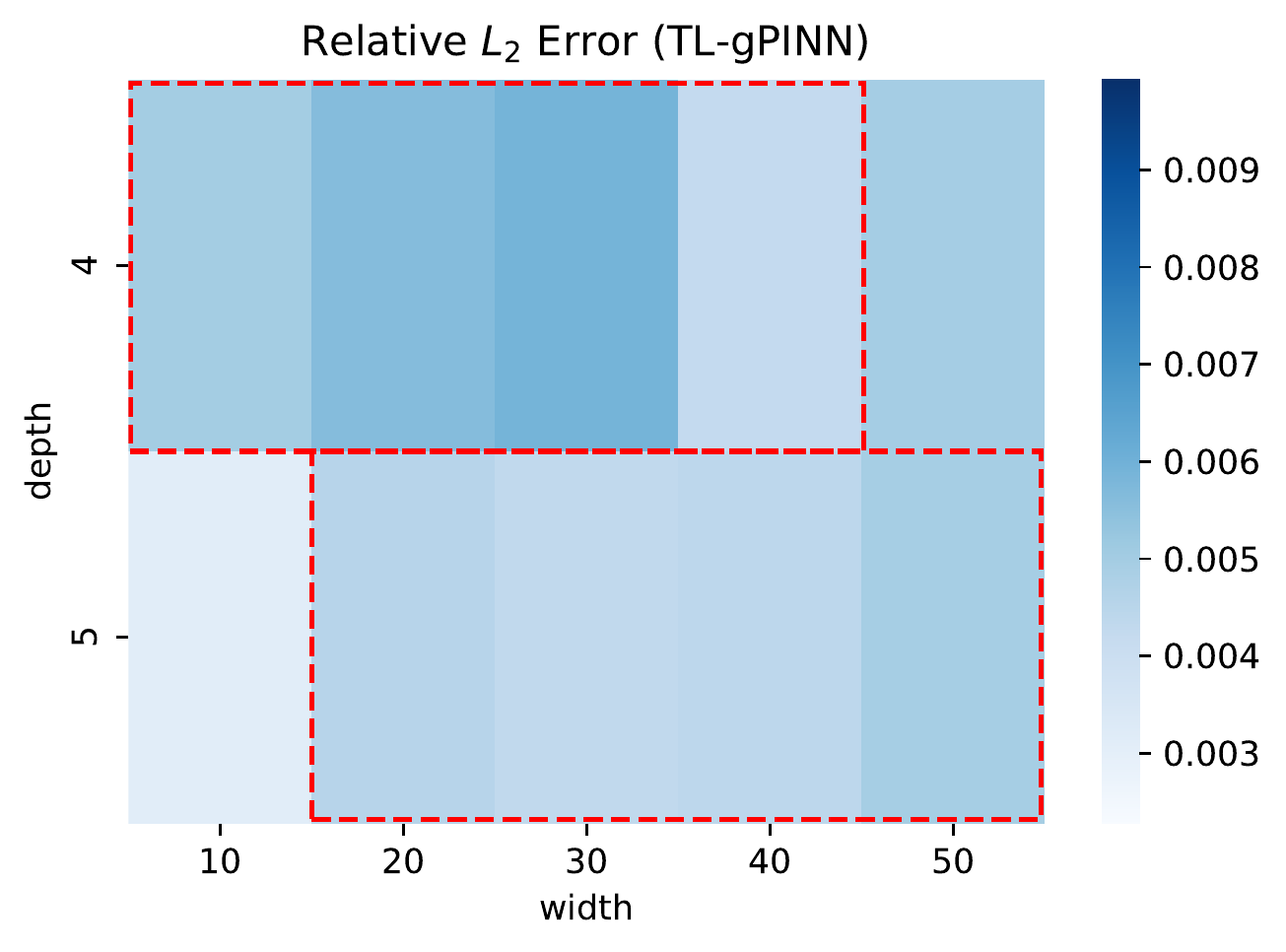}
\includegraphics[width=5.5cm,height=3.5cm]{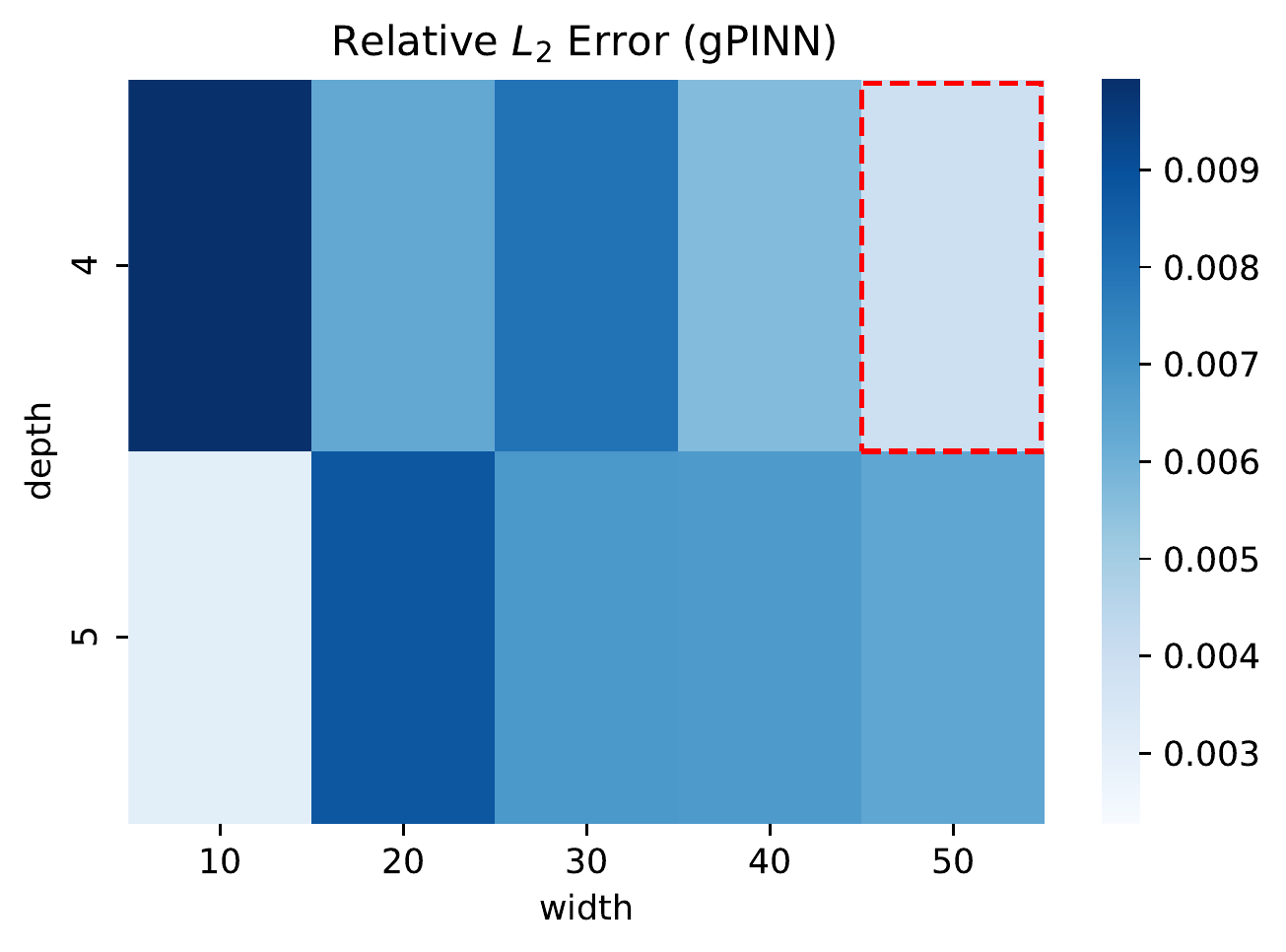}\\
$a$\\
\includegraphics[width=5.5cm,height=3.5cm]{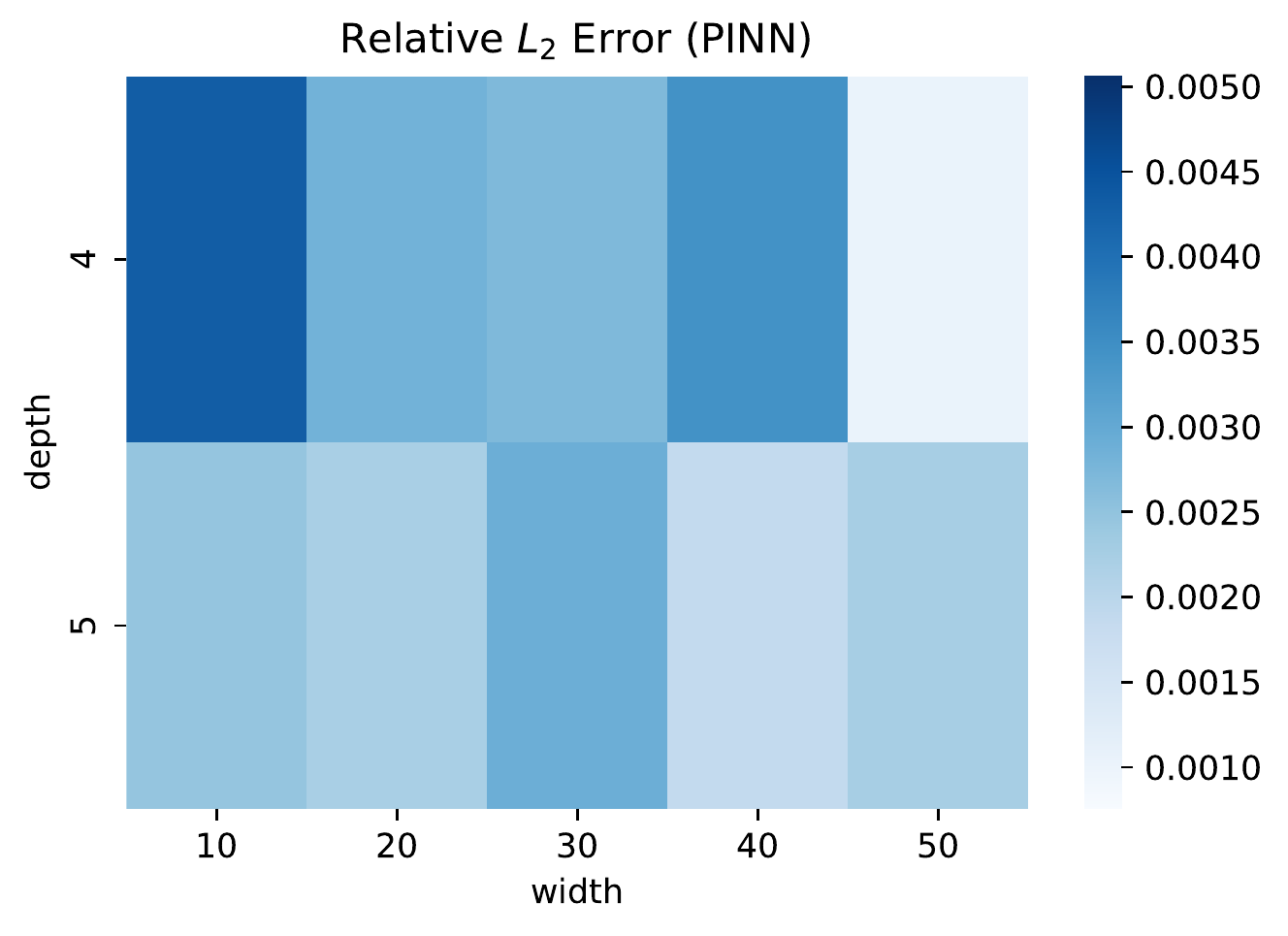}
\includegraphics[width=5.5cm,height=3.5cm]{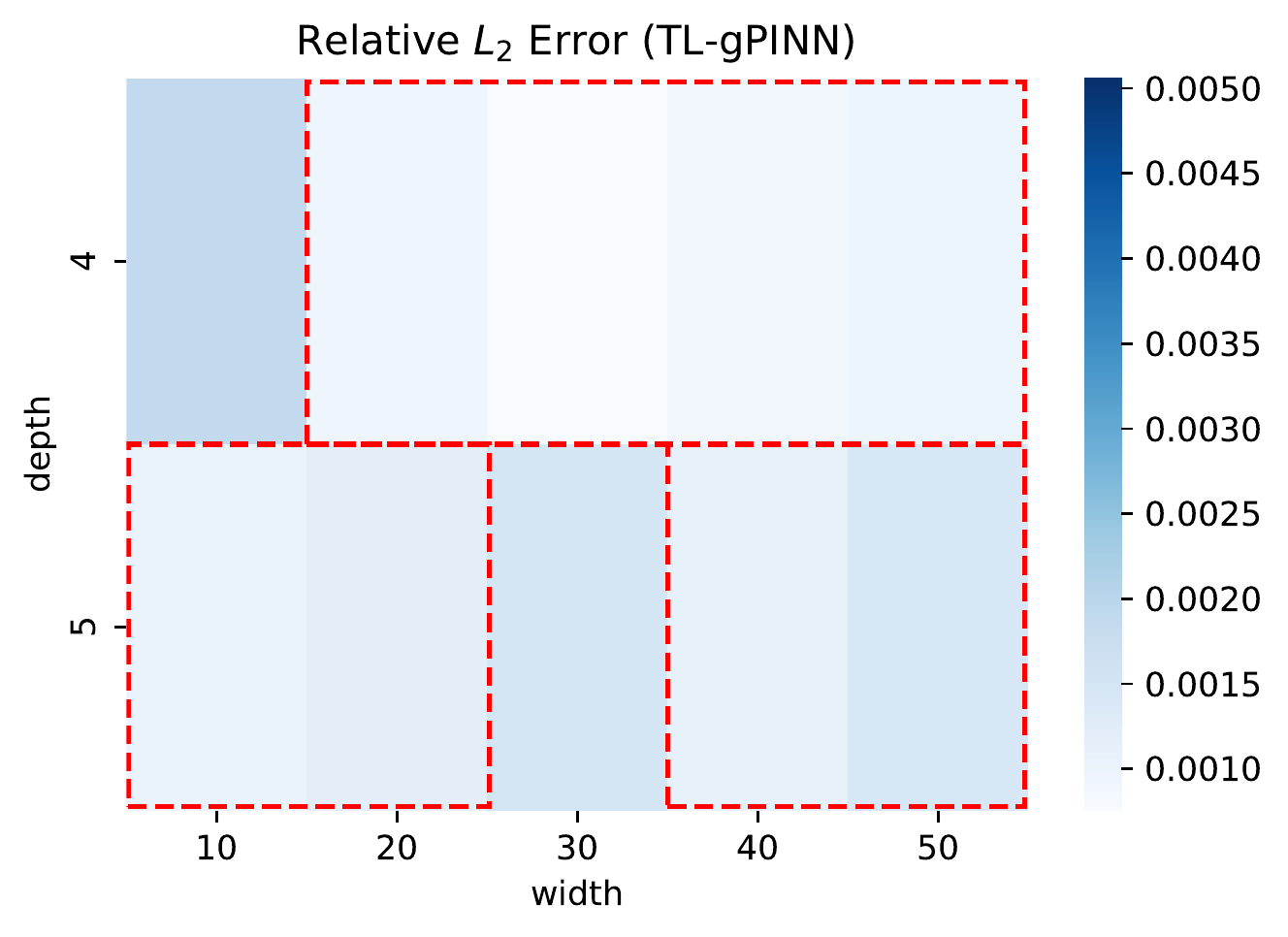}
\includegraphics[width=5.5cm,height=3.5cm]{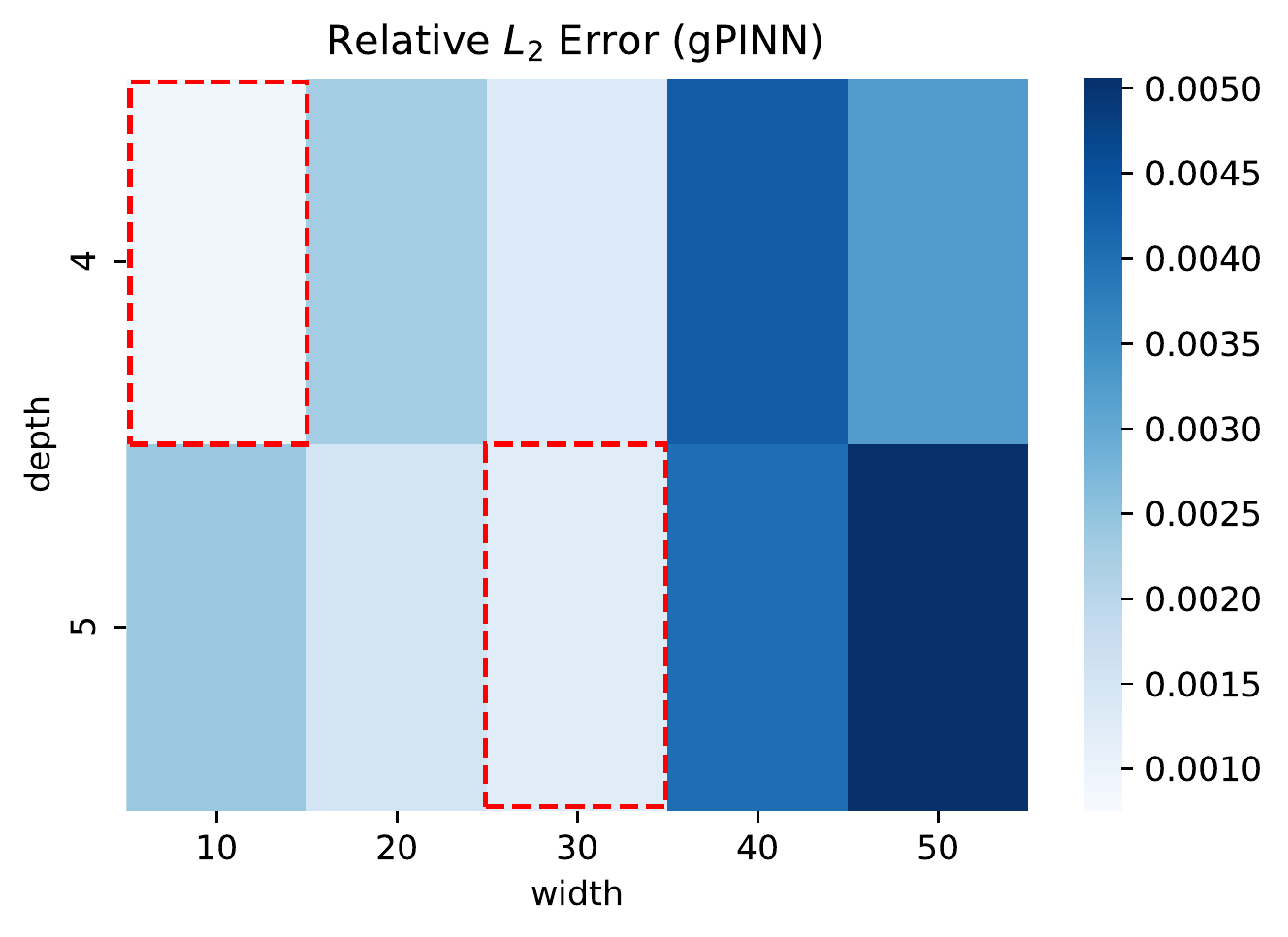}\\
$b$\\
\includegraphics[width=5.5cm,height=3.5cm]{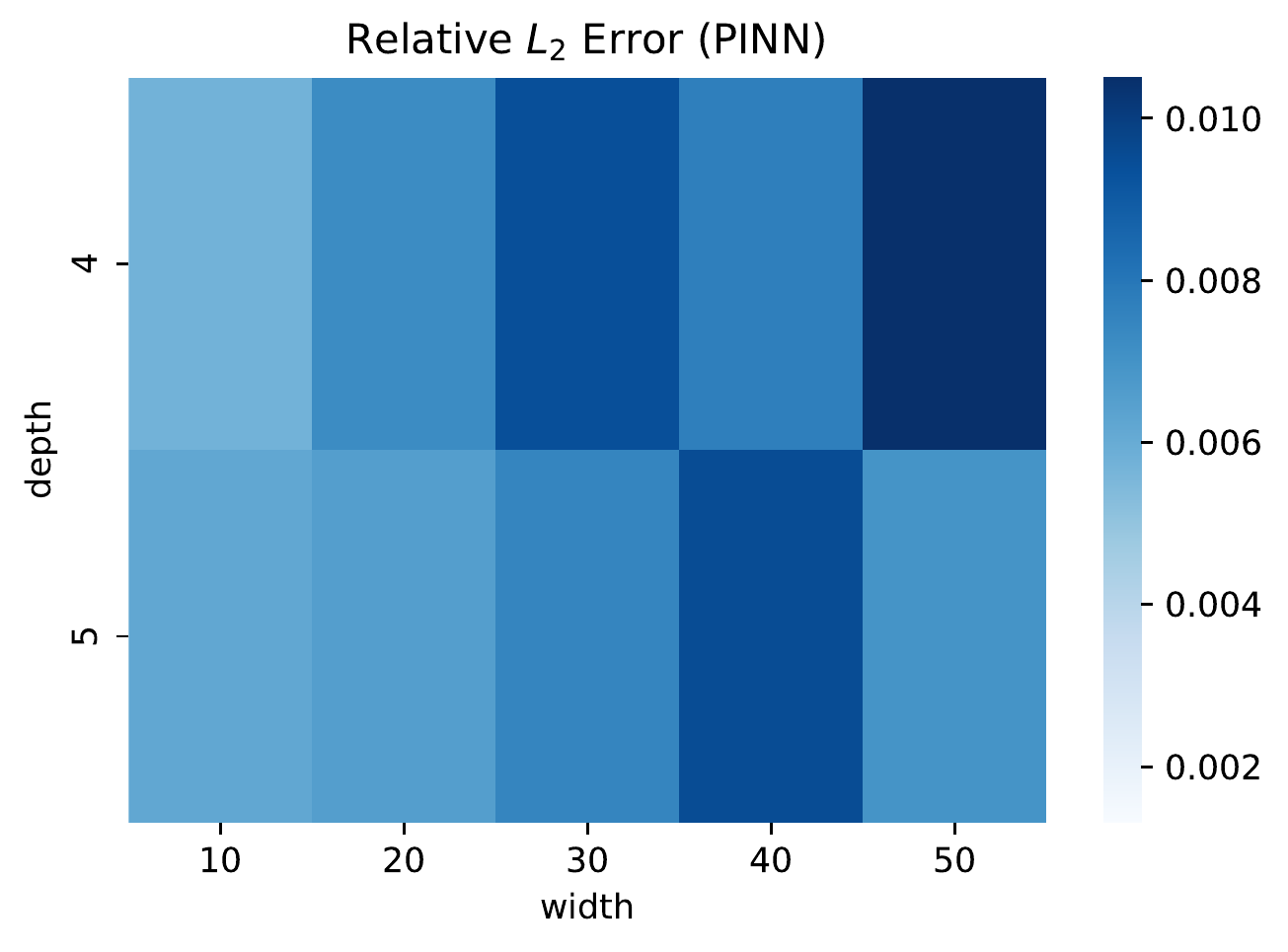}
\includegraphics[width=5.5cm,height=3.5cm]{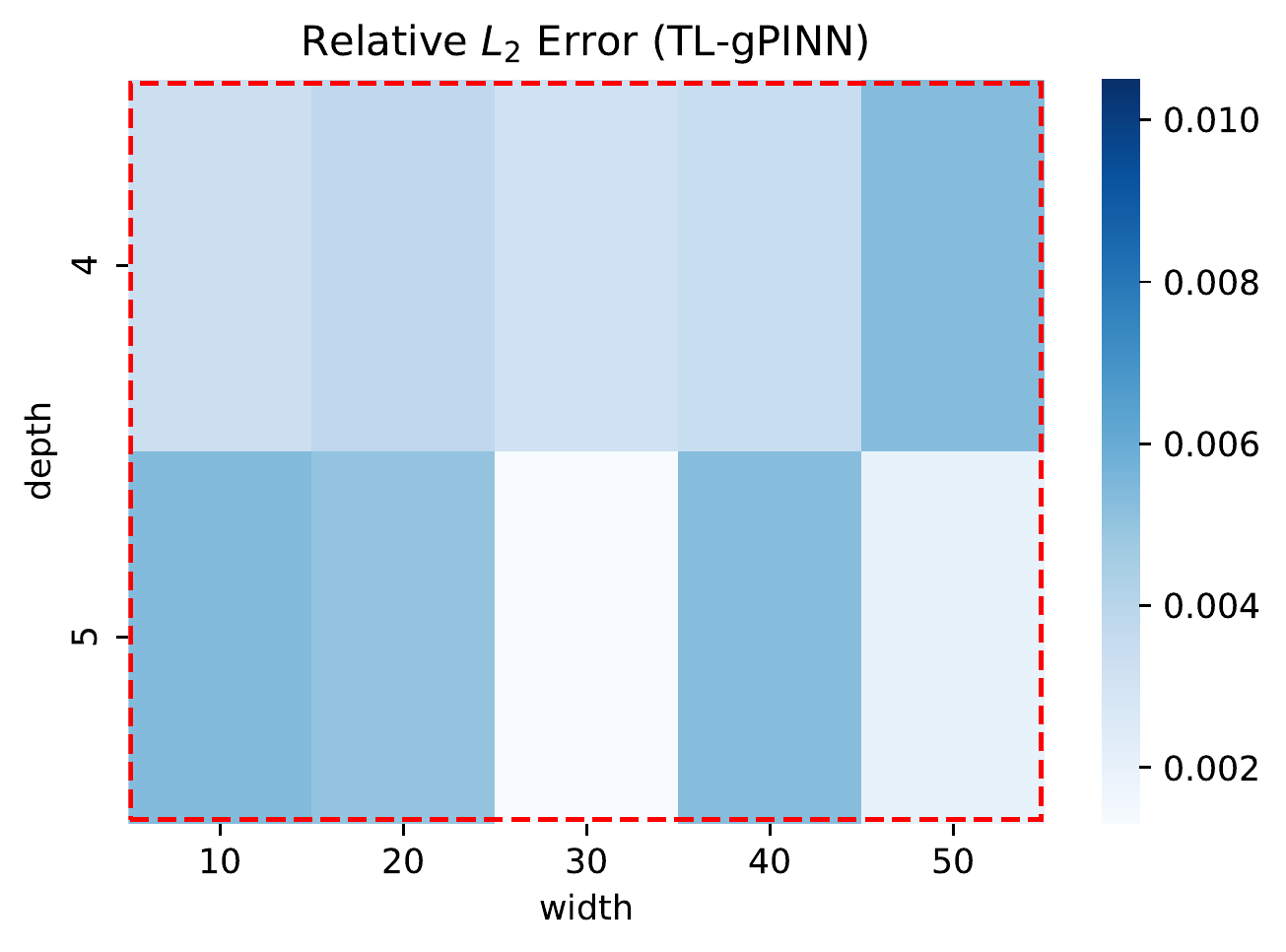}
\includegraphics[width=5.5cm,height=3.5cm]{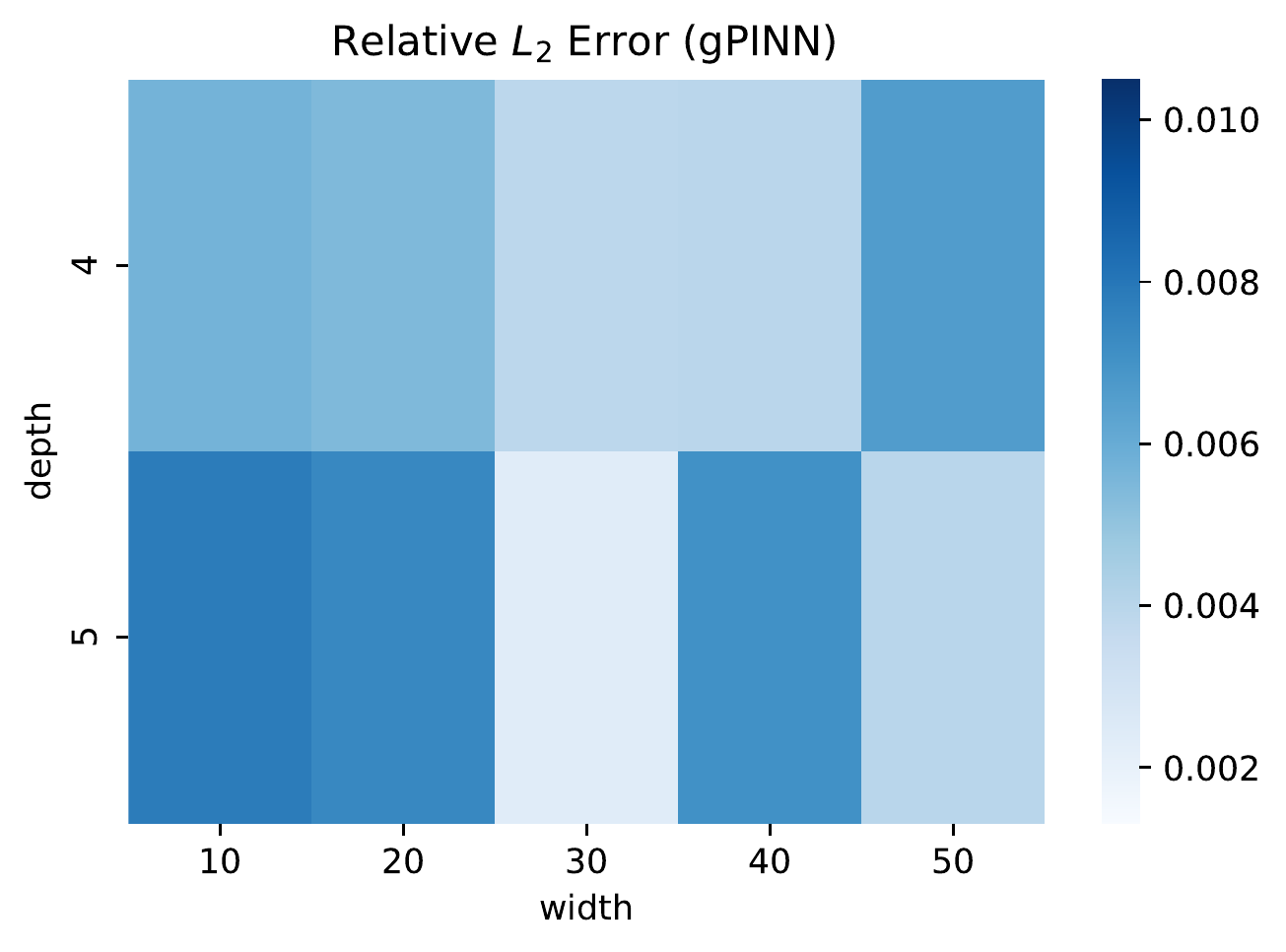}\\
$c$\\
\includegraphics[width=5.5cm,height=3.5cm]{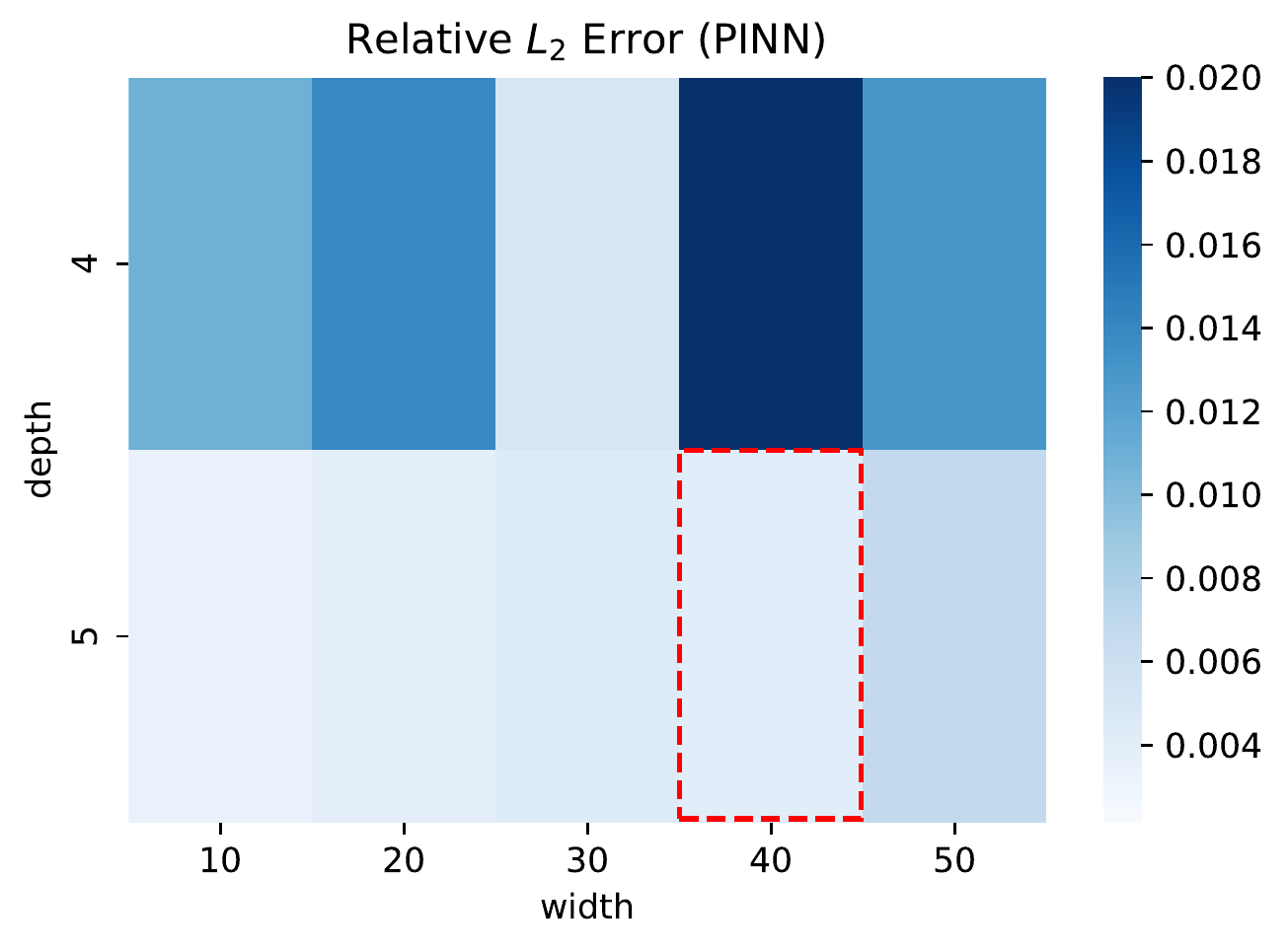}
\includegraphics[width=5.5cm,height=3.5cm]{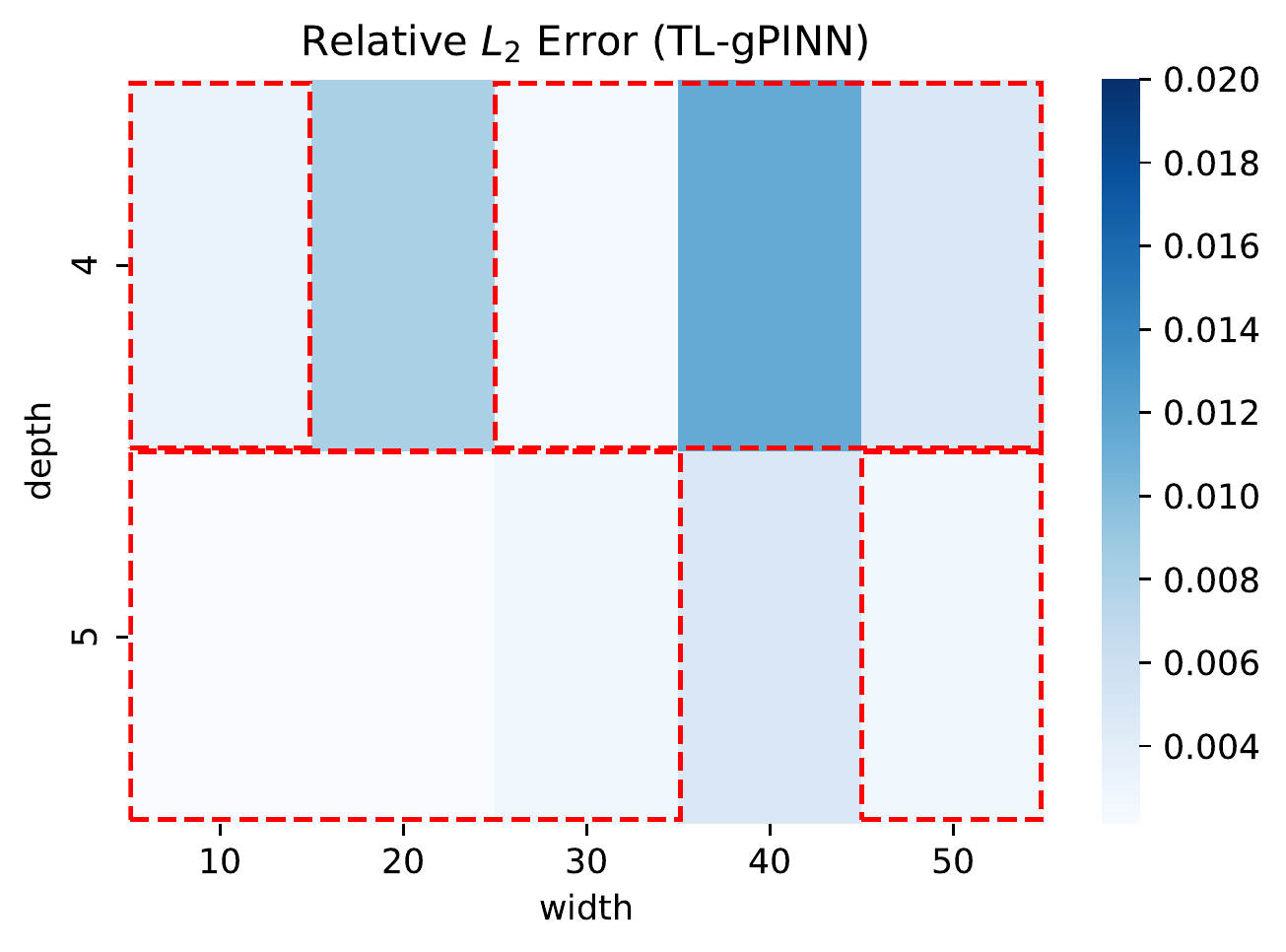}
\includegraphics[width=5.5cm,height=3.5cm]{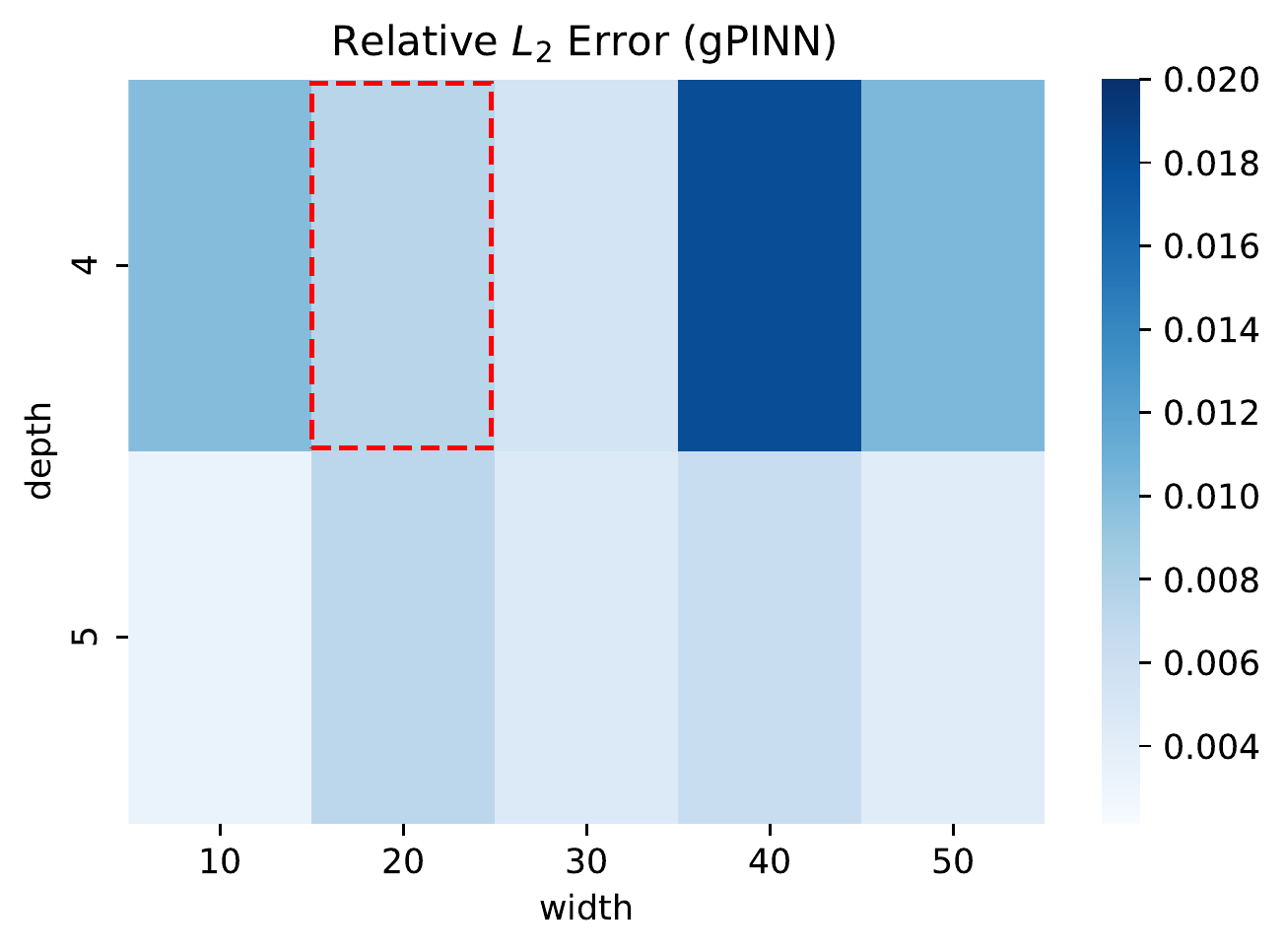}\\
$d$
\caption{(Color online) Relative $\mathbb{L}_2$ errors of nonlinear variable coefficients $\gamma(t)$ via three methods under different depth and width: (a) quadratic $\gamma(t)$; (b) sine $\gamma(t)$; (c) hyperbolic tangent $\gamma(t)$; (d) fractional $\gamma(t)$.}
\label{fig4-1}
\end{figure}

Finally, heat maps of relative $\mathbb{L}_2$ errors are shown in Fig. \ref{fig4-1} in order to display the experimental results more intuitively, and the detailed results are given in Table \ref{tableA-1} in Appendix A.

The figures in the first, second and third columns are the visualization of relative $\mathbb{L}_2$ errors given by PINNs, TL-gPINNs and gPINNs respectively. The darker the color, the greater the error. For each group of experiments, we will compare the performance of the three methods and use the red dotted line to frame the one with the smallest error in the heat maps. Evidently, the color of heat maps in the second column is the lightest on the whole. Also, the proportion of numerical experiments with the smallest error is the largest. Since the weights and biases as the initialization parameters of TL-gPINN are inherited from PINN, the color depth that reflects the value of relative $\mathbb{L}_2$ errors of PINN and TL-gPINN is highly correlated according to heat maps in Fig. \ref{fig4-1}. It may contribute to the stability of TL-gPINN in significant accuracy enhancement.

Numerically, the average (10 runs) relative $\mathbb{L}_2$ errors of nonlinear variable coefficient $\gamma(t)$ as well as the error reduction rates of TL-gPINNs and gPINNs are listed in Table \ref{table4-2}. Undoubtedly, it illustrates that our proposed method (TL-gPINN) outperforms the other two (PINN and gPINN) thoroughly.

For numerous cases above, TL-gPINN always performs well and has stable improvement of accuracy under different width and depth of the branch network for the identification of nonlinear variable coefficients.

\begin{table}[htbp]
\caption{Average performance comparison of three methods in identifying nonlinear variable coefficient $\gamma(t)$ for the vcNLS equation under different width and depth.}
\label{table4-2}  
\centering
\begin{tabular}{c|ccc}
\bottomrule
\multirow{2}{*}{\begin{tabular}[c]{@{}c@{}}Correct nonlinear $\gamma(t)$\end{tabular}} & \multicolumn{3}{c}{Relative $\mathbb{L}_2$ errors($ERR_2$)}                                         \\ \cline{2-4} 
  & \multicolumn{1}{c|}{PINNs} & \multicolumn{1}{c|}{gPINNs} & TL-gPINNs \\ \hline
 $t^2$    & 5.375348e-03     & 6.562848e-03\textbf{(-22.09\%)}           &   4.711041e-03\textbf{(12.36\%)}    \\
$\sin(t)$      & 2.603825e-03         & 2.651862e-03\textbf{(-1.84\%)}          &  1.184149e-03\textbf{(54.52\%)}      \\
$\tanh(t)$      & 7.727333e-03  &   5.430929e-03\textbf{(29.72\%)}     &    3.814903e-03\textbf{(50.63\%)}    \\
$\frac{1}{1+t^2}$    &  8.590022e-03    & 7.668054e-03\textbf{(10.73\%)}          &  4.530164e-03\textbf{(47.26\%)}     \\ \toprule
\end{tabular}
\end{table}

$\bullet$ \textbf{The size of training dataset}

The difference between the inverse problem and the forward one lies in the incorporation of some extra measurements $\{x^i_{in},t^i_{in},u^i,v^i\}^{N_{A_{in}}}_{i=1}$ of the internal region. Hence, the major consideration is the size of internal data, i.e. the value of $N_{A_{in}}$.

Considering the randomness involved in sampling and initialization, the setting of the parameter $seed$ in the codes will affect the numerical results. We perform six groups of numerical experiments for each nonlinear variable coefficient $\gamma(t)$ and the value of $N_{A_{in}}$ changes from 500 to 3000 with step size 500. Meanwhile, each group contains five experiments under the condition of different initial seeds to explore the impact of randomness on the results.

\begin{figure}[htbp]
\centering
\includegraphics[width=7cm,height=5cm]{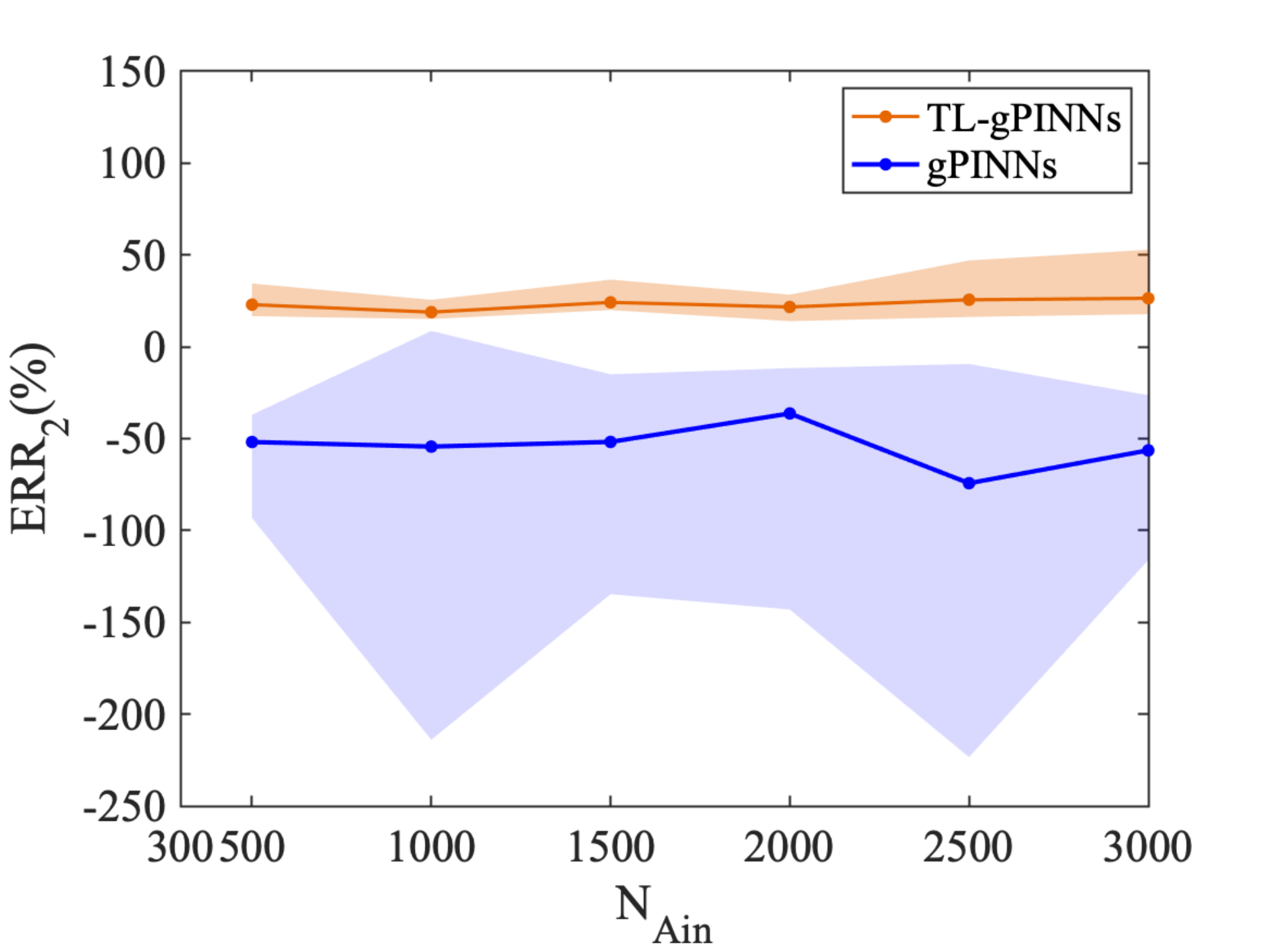}
$a$
\includegraphics[width=7cm,height=5cm]{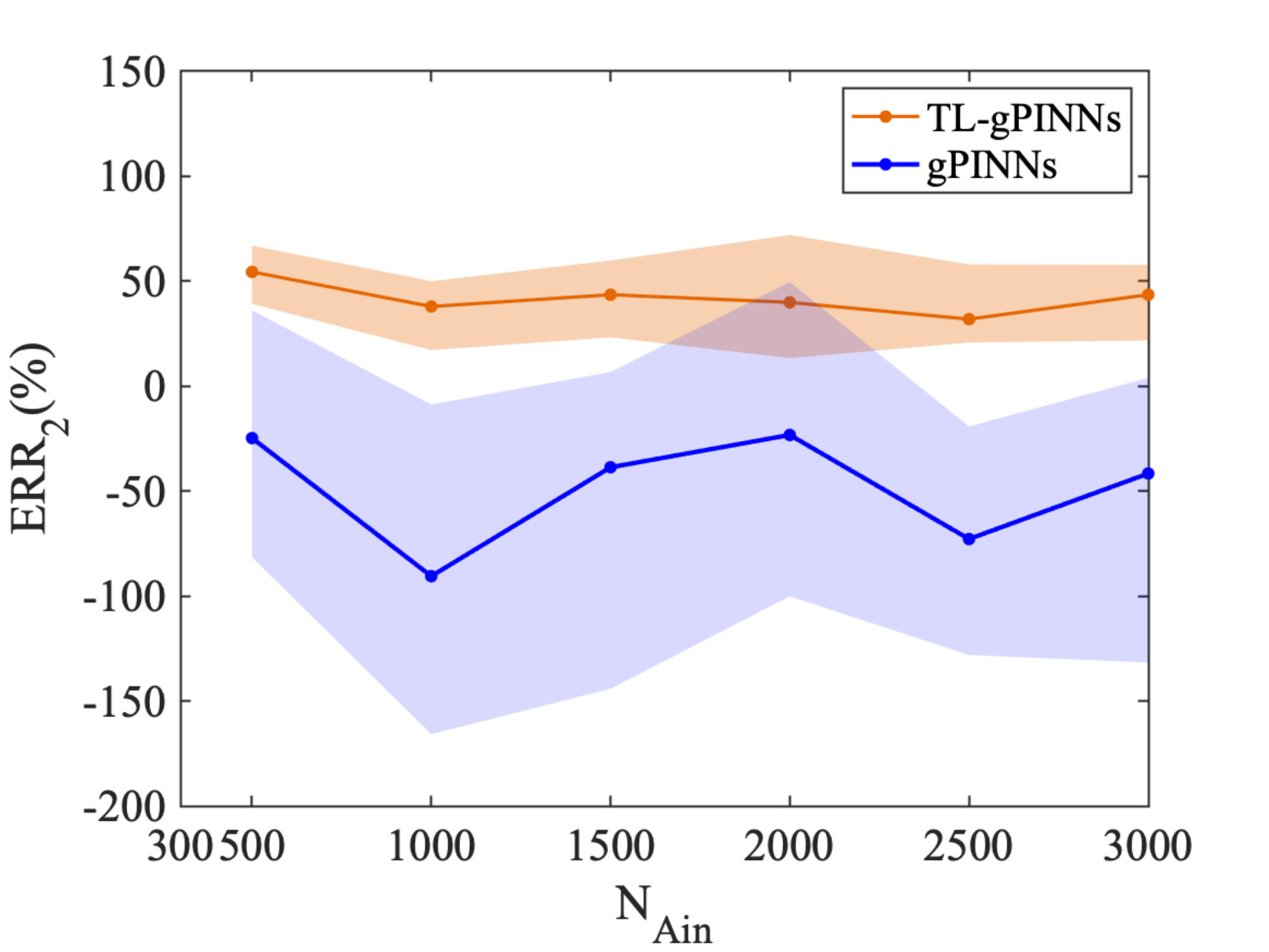}
$b$\\
\includegraphics[width=7cm,height=5cm]{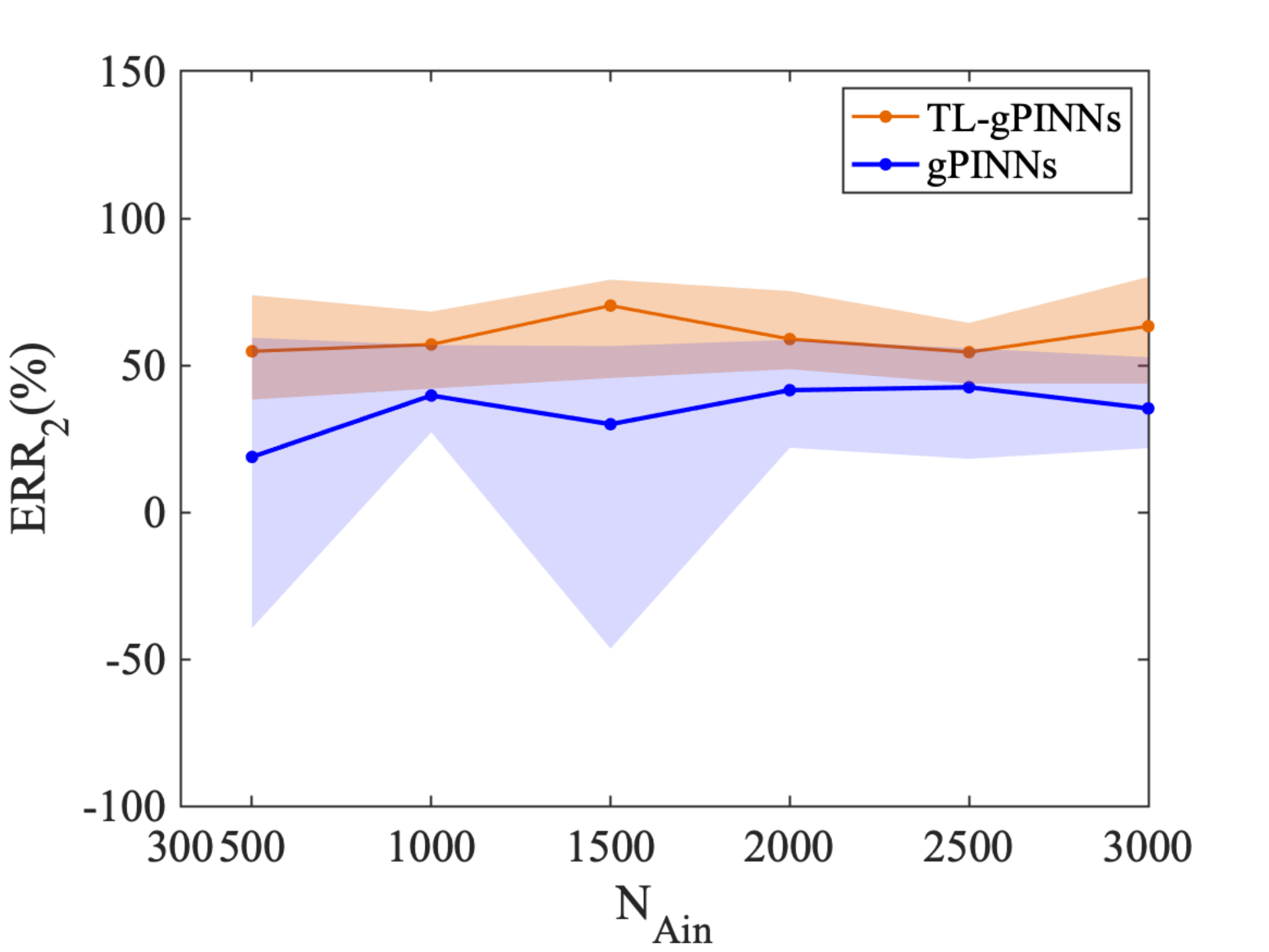}
$c$
\includegraphics[width=7cm,height=5cm]{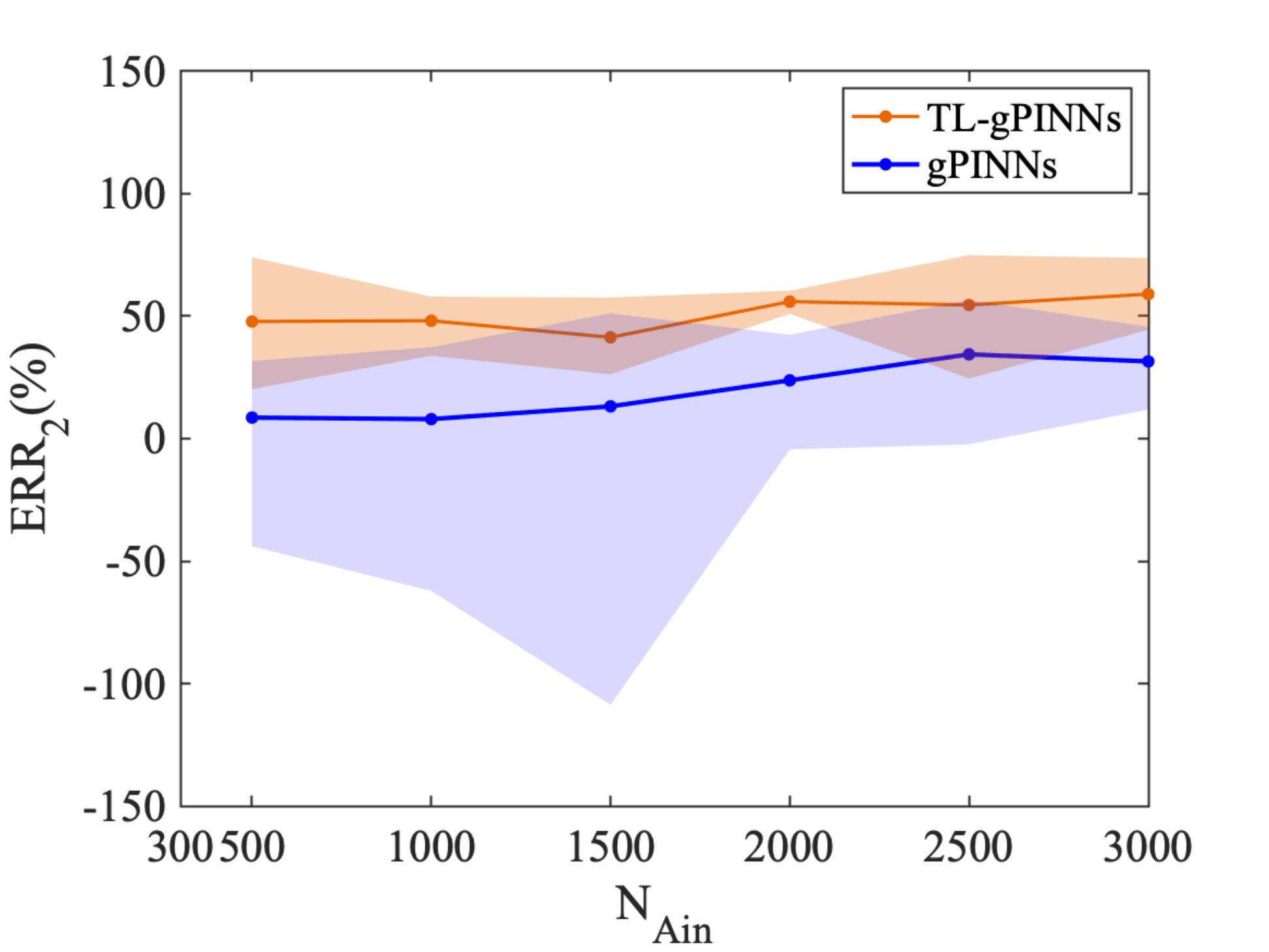}
$d$
\caption{(Color online) Error reduction rates of relative $\mathbb{L}_2$ error ($ERR_2$) in identifying nonlinear variable coefficient $\gamma(t)$ for the vcNLS equation achieved by TL-gPINNs and gPINNs compared with PINNs under different number of $N_{A_{in}}$: (a) quadratic $\gamma(t)$; (b) sine $\gamma(t)$; (c) hyperbolic tangent $\gamma(t)$; (d) fractional $\gamma(t)$.}
\label{fig4-2}
\end{figure}

Here, we are chiefly concerned with the accuracy of the nonlinear $\gamma(t)$ obtained by TL-gPINNs as well as the error reduction rates of TL-gPINNs and gPINNs compared with PINNs, which are shown in Fig. \ref{fig4-2} and Table \ref{tableA-2} in Appendix A. In Fig. \ref{fig4-2}, the orange and blue lines correspond to the mean error reduction rates ($ERR_2$) of five numerical experiments by using TL-gPINNs and gPINNs respectively, and the shade regions depict the max-min ones. It can be concluded from figures above that TL-gPINN has higher error reduction rates for each case whether in average, maximum, or minimum sense. However, $ERR_2$ of gPINNs is even less than 0\% in many examples, which means the accuracy of gPINN is reduced rather than improved compared to the traditional PINN method. Furthermore, the size of the shaded area to some extent reflects the stability of the method. Thus, TL-gPINN apparently is more stable and accurate than gPINN based on error reduction rates of relative $\mathbb{L}_2$ error ($ERR_2$) under different size of training dataset.

\section{Conclusion}

Traditional numerical methods have many limitations in solving inverse problems, especially in dealing with noisy data, complex regions, and high-dimensional problems. Moreover, the inverse problem of the function discovery is a relatively under explored field. In this paper, for the sake of overcoming deficiency of the discrete characterization of the PDE loss in neural networks and improving accuracy of function feature description, we propose gradient-enhanced PINNs based on transfer learning (TL-gPINNs) for inverse problems of inferring unknown variable coefficients and give a new viewpoint on gPINNs.

The TL-gPINN method uses a two-step optimization strategy and gradually increases the difficulty. Firstly, the original PINN is applied in the inverse problem of the variable coefficient equations. Then for further optimization, gPINN inherits the saved weight matrixes and bias vectors of PINN at the end of the iteration process as the initialization parameters and the introduction of the gradient term contributes to the accuracy enhancement of variable coefficients. Moreover, the trunk and branch networks are  established to infer the solution and variable coefficients separately in order to eliminate mutual influence.

The effectiveness of TL-gPINNs is demonstrated in identifying several types of single variable coefficients, including linear, quadratic, sine, hyperbolic tangent and fractional functions as well as multiple ones for the well-known variable coefficient nonlinear Schr\"{o}odinger (vcNLS) equation in the field of integrable systems. Meanwhile, abundant dynamic behaviors of the corresponding soliton solution can be well reproduced. Plenty of numerical experiments are carried out to compare the performance of PINNs, TL-gPINNs and gPINNs. It has been proved that gPINN learned the unknown parameters more accurately than PINN for the inverse problems in many examples, such as Poisson equation, diffusion-reaction equation, Brinkman-Forchheimer model and so on by Yu et al. However, the accuracy of gPINN is reduced rather than improved compared with the standard PINN method in inverse PDE problems of the vcNLS equation. Presumably it's because the loss function itself consists of many constraint terms even without regard to the gradient restriction and thus the result may not necessarily be better even if more constraints are imposed when solving the multi-objective optimization problems. What's worse, the computational cost of gPINN is higher than PINN unavoidably since the introduction of additional constraints on gradients gives rise to low efficiency. Consequently, one viable path towards accelerating the convergence of training could come by adopting the technique of transfer learning and thus the TL-gPINN method is put forward here. Through the comparison among the three methods, TL-gPINN has the highest prediction precision and can improve efficiency compared to gPINN. In other words, TL-gPINNs can successfully infer the unknown variable coefficients with satisfactory accuracy and outperform the PINNs in accuracy, and gPINNs in both accuracy and efficiency. Besides, we ulteriorly conduct robustness analysis and parametric sensitivity analysis. Numerical results also illustrate that the ability of TL-gPINNs to improve accuracy compared to the standard PINNs, and gPINNs also remains robust to noise and other hyper-parameters, including width and depth of the branch network and the size of training dataset.

The TL-gPINN method applied in this paper is universal and can be adapted to the inverse problems of inferring unknown high-dimensional variable coefficients. In future work, we will strive to propose more targeted improvements that enhance accuracy without sacrificing efficiency on this subject.

\section*{Acknowledgments}
The authors would like to thank Zhengwu Miao sincerely for providing with support and help. The project is supported by National Natural Science Foundation of China (No. 12175069 and No. 12235007), Science and Technology Commission of Shanghai Municipality (No. 21JC1402500 and No. 22DZ2229014) and Natural Science Foundation of Shanghai (No. 23ZR1418100).

\section*{Appendix A.} \label{Appendix}

\begin{table}[htbp]
\caption{Performance comparison of three methods in identifying variable coefficient $\gamma(t)$ for the vcNLS equation under different noise conditions.}
\label{tableA-0}  
\centering
\begin{tabular}{ccc|ccccc}
\bottomrule
\multicolumn{3}{c|}{\multirow{2}{*}{\diagbox{Results}{Solution types}}}                                                   & \multicolumn{5}{c}{Correct $\gamma(t)$} \\ \cline{4-8} 
\multicolumn{3}{c|}{}                                                                               & $t$   & $t^2$   & $\sin(t)$   & $\tanh(t)$   & $\frac{1}{1+t^2}$   \\ \hline
\multicolumn{1}{c|}{\multirow{10}{*}{0.5\% noise}} & \multicolumn{1}{c|}{\multirow{2}{*}{PINNs}}     & $MAE_{\gamma}$ & 1.322670e-04   & 9.929034e-03   &    2.189932e-03& 8.245944e-03  & 1.063432e-03 \\ \cline{3-3}
\multicolumn{1}{c|}{}                        & \multicolumn{1}{c|}{}                           & $RE_{\gamma}$  & 6.551441e-05   & 8.084222e-03   &4.279738e-03    & 1.111883e-02  & 3.326201e-03  \\ \cline{2-3}
\multicolumn{1}{c|}{}                        & \multicolumn{1}{c|}{\multirow{4}{*}{TL-gPINNs}} & $MAE_{\gamma}$ & 1.322670e-04   & 9.615011e-03   & 1.097929e-03   & 3.973503e-03  & 1.009619e-03  \\ \cline{3-3}
\multicolumn{1}{c|}{}                        & \multicolumn{1}{c|}{}                           & $RE_{\gamma}$  & 6.551441e-05   & 7.959068e-03   & 2.579889e-03   &5.486678e-03   & 2.968212e-03  \\ \cline{3-3}
\multicolumn{1}{c|}{}                        & \multicolumn{1}{c|}{}                           & $ERR_1$ & 0.00\%   &3.16\%    & 49.86\%   & 51.81\%  &5.06\%   \\ \cline{3-3}
\multicolumn{1}{c|}{}                        & \multicolumn{1}{c|}{}                           & $ERR_2$ & 0.00\%   & 1.55\%   & 39.72\%   & 50.65\%  &10.76\%   \\ \cline{2-3}
\multicolumn{1}{c|}{}                        & \multicolumn{1}{c|}{\multirow{4}{*}{gPINNs}}    & $MAE_{\gamma}$ & 2.300093e-04   & 1.886554e-02   & 2.773028e-03   & 5.021791e-03  &1.176551e-03   \\ \cline{3-3}
\multicolumn{1}{c|}{}                        & \multicolumn{1}{c|}{}                           & $RE_{\gamma}$  & 1.150104e-04   & 1.308772e-02   & 5.779078e-03   &6.675003e-03   &3.859277e-03   \\ \cline{3-3}
\multicolumn{1}{c|}{}                        & \multicolumn{1}{c|}{}                           & $ERR_1$ & -73.90\%   & -90.00\%   & -26.63\%   & 39.10\%  & -10.64\%  \\ \cline{3-3}
\multicolumn{1}{c|}{}                        & \multicolumn{1}{c|}{}                           & $ERR_2$ & -75.55\%   & -61.89\%   & -35.03\%   & 39.97\%  & -16.03\%  \\ \hline
\multicolumn{1}{c|}{\multirow{10}{*}{1\% noise}} & \multicolumn{1}{c|}{\multirow{2}{*}{PINNs}}     & $MAE_{\gamma}$ &4.643823e-04    &2.069045e-02    & 2.843303e-03   & 6.361015e-03  &2.513644e-03   \\ \cline{3-3}
\multicolumn{1}{c|}{}                        & \multicolumn{1}{c|}{}                           & $RE_{\gamma}$  &2.317933e-04    & 2.002579e-02   & 7.413523e-03   & 8.855182e-03  & 7.986546e-03  \\ \cline{2-3}
\multicolumn{1}{c|}{}                        & \multicolumn{1}{c|}{\multirow{4}{*}{TL-gPINNs}} & $MAE_{\gamma}$ &3.393921e-04    & 1.990454e-02   & 2.040646e-03   & 5.745897e-03  & 1.529080e-03  \\ \cline{3-3}
\multicolumn{1}{c|}{}                        & \multicolumn{1}{c|}{}                           & $RE_{\gamma}$  &1.693800e-04    & 2.038240e-02   & 6.616817e-03   & 7.922166e-03  & 4.553724e-03  \\ \cline{3-3}
\multicolumn{1}{c|}{}                        & \multicolumn{1}{c|}{}                           & $ERR_1$ &26.92\%    & 3.80\%   & 28.23\%   &9.67\%  & 39.17\%  \\ \cline{3-3}
\multicolumn{1}{c|}{}                        & \multicolumn{1}{c|}{}                           & $ERR_2$ & 26.93\%   & -1.78\%   & 10.75\%   &10.54\%   &42.98\%   \\ \cline{2-3}
\multicolumn{1}{c|}{}                        & \multicolumn{1}{c|}{\multirow{4}{*}{gPINNs}}    & $MAE_{\gamma}$ & 2.450279e-04   &  1.679156e-02  &4.735398e-03    & 6.096149e-03  &1.782366e-03   \\ \cline{3-3}
\multicolumn{1}{c|}{}                        & \multicolumn{1}{c|}{}                           & $RE_{\gamma}$  & 1.221250e-04   & 1.835826e-02   &1.017474e-02    & 8.625601e-03  & 5.073325e-03  \\ \cline{3-3}
\multicolumn{1}{c|}{}                        & \multicolumn{1}{c|}{}                           & $ERR_1$ & 47.24\%   & 18.84\%   & -66.55\%   &4.16\%   &29.09\%   \\ \cline{3-3}
\multicolumn{1}{c|}{}                        & \multicolumn{1}{c|}{}                           & $ERR_2$ & 47.31\%   & 8.33\%   & -37.25\%   &2.59\%   &36.48\%   \\ \hline
\multicolumn{1}{c|}{\multirow{10}{*}{3\% noise}} & \multicolumn{1}{c|}{\multirow{2}{*}{PINNs}}     & $MAE_{\gamma}$ & 7.029357e-04   & 1.459252e-02   & 3.614353e-03   & 1.109293e-02  &2.486608e-03   \\ \cline{3-3}
\multicolumn{1}{c|}{}                        & \multicolumn{1}{c|}{}                           & $RE_{\gamma}$  & 3.517501e-04   & 1.405216e-02   & 8.000196e-03   & 1.588254e-02  & 9.198019e-03  \\ \cline{2-3}
\multicolumn{1}{c|}{}                        & \multicolumn{1}{c|}{\multirow{4}{*}{TL-gPINNs}} & $MAE_{\gamma}$ & 3.642581e-04   & 1.241316e-02   & 3.276918e-03   & 1.005220e-02  & 2.160257e-03  \\ \cline{3-3}
\multicolumn{1}{c|}{}                        & \multicolumn{1}{c|}{}                           & $RE_{\gamma}$  & 1.825354e-04   & 1.403889e-02    & 8.239684e-03   & 1.449078e-02  &6.147466e-03   \\ \cline{3-3}
\multicolumn{1}{c|}{}                        & \multicolumn{1}{c|}{}                           & $ERR_1$ & 48.18\%   & 14.93\%   & 9.34\%   & 9.38\%  &13.12\%   \\ \cline{3-3}
\multicolumn{1}{c|}{}                        & \multicolumn{1}{c|}{}                           & $ERR_2$ & 48.11\%   & 0.09\%   & -2.99\%   & 8.76\%  &33.17\%   \\ \cline{2-3}
\multicolumn{1}{c|}{}                        & \multicolumn{1}{c|}{\multirow{4}{*}{gPINNs}}    & $MAE_{\gamma}$ &5.101351e-04    & 1.470036e-02   & 4.077128e-03   & 9.832798e-03  &2.296355e-03   \\ \cline{3-3}
\multicolumn{1}{c|}{}                        & \multicolumn{1}{c|}{}                           & $RE_{\gamma}$  & 2.561838e-04   & 1.436189e-02   & 8.476734e-03   & 1.415209e-02  & 9.024564e-03  \\ \cline{3-3}
\multicolumn{1}{c|}{}                        & \multicolumn{1}{c|}{}                           & $ERR_1$ & 27.43\%   & -0.74\%   & -12.80\%   & 11.36\%  &7.65\%   \\ \cline{3-3}
\multicolumn{1}{c|}{}                        & \multicolumn{1}{c|}{}                           & $ERR_2$ & 27.17\%   & -2.20\%   & -5.96\%   &10.90\%   &1.89\%   \\ \hline
\multicolumn{1}{c|}{\multirow{10}{*}{5\% noise}} & \multicolumn{1}{c|}{\multirow{2}{*}{PINNs}}     & $MAE_{\gamma}$ & 4.427668e-04   & 4.122176e-02   & 7.337520e-03   & 9.317941e-03  & 3.243209e-03  \\ \cline{3-3}
\multicolumn{1}{c|}{}                        & \multicolumn{1}{c|}{}                           & $RE_{\gamma}$  & 2.190762e-04   & 3.776385e-02   & 1.687885e-02   & 1.416438e-02  & 1.190859e-02  \\ \cline{2-3}
\multicolumn{1}{c|}{}                        & \multicolumn{1}{c|}{\multirow{4}{*}{TL-gPINNs}} & $MAE_{\gamma}$ & 3.007159e-04   & 3.890978e-02   & 5.462772e-03   & 8.254448e-03  & 2.393983e-03  \\ \cline{3-3}
\multicolumn{1}{c|}{}                        & \multicolumn{1}{c|}{}                           & $RE_{\gamma}$  & 1.475556e-04   & 3.750857e-02   & 1.526540e-02   & 1.280577e-02  & 7.430092e-03  \\ \cline{3-3}
\multicolumn{1}{c|}{}                        & \multicolumn{1}{c|}{}                           & $ERR_1$ & 32.08\%   & 5.61\%   &25.55\%    & 11.41\%  &26.18\%   \\ \cline{3-3}
\multicolumn{1}{c|}{}                        & \multicolumn{1}{c|}{}                           & $ERR_2$ & 32.65\%   & 0.68\%   & 9.56\%   & 9.59\%  & 37.61\%  \\ \cline{2-3}
\multicolumn{1}{c|}{}                        & \multicolumn{1}{c|}{\multirow{4}{*}{gPINNs}}    & $MAE_{\gamma}$ & 3.408351e-04   & 3.847756e-02   &5.860518e-03    &1.159346e-02   &4.138773e-03   \\ \cline{3-3}
\multicolumn{1}{c|}{}                        & \multicolumn{1}{c|}{}                           & $RE_{\gamma}$  & 1.679388e-04   & 3.838483e-02   & 1.568930e-02   & 1.755944e-02  &1.279140e-02   \\ \cline{3-3}
\multicolumn{1}{c|}{}                        & \multicolumn{1}{c|}{}                           & $ERR_1$ &23.02\%    & 6.66\%   & 20.13\%   &-24.42\%   &-27.61\%   \\ \cline{3-3}
\multicolumn{1}{c|}{}                        & \multicolumn{1}{c|}{}                           & $ERR_2$ & 23.34\%   & -1.64\%   &7.05\%    & -23.97\%  &-7.41\%   \\ \toprule
\end{tabular}
\end{table}

\begin{table}[htbp]
\caption{Relative $\mathbb{L}_2$ errors of three methods in identifying nonlinear variable coefficient $\gamma(t)$ for the vcNLS equation under different depth and width.}
\label{tableA-1}  
\centering
\begin{tabular}{cc|cccc}
\bottomrule
\multicolumn{2}{c|}{\multirow{2}{*}{\diagbox{Depth-Width}{Solution types}}}          & \multicolumn{4}{c}{Correct nonlinear $\gamma(t)$}                                    \\ \cline{3-6} 
\multicolumn{2}{c|}{}                                  & \multicolumn{1}{c|}{$t^2$ } & \multicolumn{1}{c|}{$\sin(t)$} & \multicolumn{1}{c|}{$\tanh(t)$} & $\frac{1}{1+t^2}$ \\ \hline
\multicolumn{1}{c|}{\multirow{3}{*}{4-10}} & PINNs     & 5.548427e-03                      & 4.313755e-03                      & 5.736964e-03                      & 1.094201e-02 \\ \cline{2-2}
\multicolumn{1}{c|}{}                      & TL-gPINNs &  5.009717e-03                     & 1.879644e-03                      & 3.354807e-03                      & 3.421307e-03 \\ \cline{2-2}
\multicolumn{1}{c|}{}                      & gPINNs    & 9.929734e-03                      & 9.520874e-04                      & 5.707187e-03                      & 9.942505e-03  \\ \cline{1-2}
\multicolumn{1}{c|}{\multirow{3}{*}{4-20}} & PINNs     & 6.289173e-03                      & 2.827434e-03                      & 7.257074e-03                      & 1.397660e-02 \\ \cline{2-2}
\multicolumn{1}{c|}{}                      & TL-gPINNs & 5.602426e-03                      &  9.669967e-04                     & 3.819827e-03                      &8.108070e-03  \\ \cline{2-2}
\multicolumn{1}{c|}{}                      & gPINNs    & 6.289173e-03                      & 2.301905e-03                      &5.461584e-03                       & 7.388081e-03 \\ \cline{1-2}
\multicolumn{1}{c|}{\multirow{3}{*}{4-30}} & PINNs     & 6.993137e-03                      &  2.703498e-03                     & 9.421331e-03                      & 5.177898e-03 \\ \cline{2-2}
\multicolumn{1}{c|}{}                      & TL-gPINNs & 5.906274e-03                      & 7.559607e-04                      & 3.178629e-03                      &2.536767e-03  \\ \cline{2-2}
\multicolumn{1}{c|}{}                      & gPINNs    & 7.971344e-03                      &  1.363048e-03                     &  3.897784e-03                     & 5.403644e-03 \\ \cline{1-2}
\multicolumn{1}{c|}{\multirow{3}{*}{4-40}} & PINNs     & 5.235672e-03                      &  3.427815e-03                     & 7.712084e-03                      & 2.004701e-02 \\ \cline{2-2}
\multicolumn{1}{c|}{}                      & TL-gPINNs & 4.238258e-03                      & 8.965765e-04                      & 3.491911e-03                      & 1.139473e-02 \\ \cline{2-2}
\multicolumn{1}{c|}{}                      & gPINNs    &5.675098e-03                       & 4.312995e-03                      & 3.935904e-03                      & 1.795422e-02 \\ \cline{1-2}
\multicolumn{1}{c|}{\multirow{3}{*}{4-50}} & PINNs     & 5.680928e-03                      & 1.049763e-03                      & 1.047242e-02                      & 1.292315e-02 \\ \cline{2-2}
\multicolumn{1}{c|}{}                      & TL-gPINNs & 4.968693e-03                      & 9.878295e-04                      &5.333929e-03                       & 4.860513e-03 \\ \cline{2-2}
\multicolumn{1}{c|}{}                      & gPINNs    & 3.899222e-03                      & 3.260390e-03                      &6.635557e-03                       & 1.025656e-02 \\ \cline{1-2}
\multicolumn{1}{c|}{\multirow{3}{*}{5-10}} & PINNs     & 2.284578e-03                      & 2.471582e-03                      & 6.157078e-03                      & 3.461629e-03 \\ \cline{2-2}
\multicolumn{1}{c|}{}                      & TL-gPINNs & 3.116250e-03                      & 1.052025e-03                      & 5.378014e-03                      & 2.140225e-03 \\ \cline{2-2}
\multicolumn{1}{c|}{}                      & gPINNs    & 3.098707e-03                      & 2.417611e-03                      & 7.805791e-03                      & 3.361290e-03 \\ \cline{1-2}
\multicolumn{1}{c|}{\multirow{3}{*}{5-20}} & PINNs     & 5.756127e-03                      & 2.225477e-03                      & 6.526232e-03                      & 4.002114e-03 \\ \cline{2-2}
\multicolumn{1}{c|}{}                      & TL-gPINNs & 4.563876e-03                      & 1.198677e-03                      & 5.002341e-03                      & 2.209469e-03 \\ \cline{2-2}
\multicolumn{1}{c|}{}                      & gPINNs    & 8.777979e-03                      & 1.554303e-03                      & 7.391697e-03                      & 7.201052e-03 \\ \cline{1-2}
\multicolumn{1}{c|}{\multirow{3}{*}{5-30}} & PINNs     & 4.810807e-03                      & 2.901531e-03                      & 7.516973e-03                      & 4.447760e-03 \\ \cline{2-2}
\multicolumn{1}{c|}{}                      & TL-gPINNs & 4.322382e-03                      & 1.523483e-03                      & 1.312732e-03                      &2.936587e-03  \\ \cline{2-2}
\multicolumn{1}{c|}{}                      & gPINNs    & 6.851686e-03                      & 1.264615e-03                      &2.413443e-03                       & 4.567902e-03 \\ \cline{1-2}
\multicolumn{1}{c|}{\multirow{3}{*}{5-40}} & PINNs     &  5.682128e-03                     & 1.865532e-03                      & 9.515989e-03                      & 4.110138e-03 \\ \cline{2-2}
\multicolumn{1}{c|}{}                      & TL-gPINNs & 4.427266e-03                      & 1.121917e-03                      & 5.257534e-03                      & 4.786998e-03 \\ \cline{2-2}
\multicolumn{1}{c|}{}                      & gPINNs    & 6.765182e-03                      & 4.029679e-03                      & 7.059347e-03                      & 6.303027e-03 \\ \cline{1-2}
\multicolumn{1}{c|}{\multirow{3}{*}{5-50}} & PINNs     & 5.472502e-03                      & 2.251865e-03                      &  6.957182e-03                     & 6.811910e-03 \\ \cline{2-2}
\multicolumn{1}{c|}{}                      & TL-gPINNs & 4.955263e-03                      & 1.458381e-03                      &  2.019302e-03                     &2.906973e-03  \\ \cline{2-2}
\multicolumn{1}{c|}{}                      & gPINNs    & 6.370352e-03                      & 5.061989e-03                      & 4.000996e-03                      & 4.302256e-03 \\\toprule
\end{tabular}
\end{table}

\begin{table}[htbp]
\caption{Error reduction rates of relative $\mathbb{L}_2$ error ($ERR_2$) in identifying nonlinear variable coefficient $\gamma(t)$ for the vcNLS equation achieved by TL-gPINNs and gPINNs compared with PINNs under different number of $N_{A_{in}}$.}
\label{tableA-2}  
\centering
\begin{tabular}{ccc|cccccc}
\bottomrule
\multicolumn{3}{c|}{\diagbox{Results}{Value of $N_{A_{in}}$}}                                                                        & 500 & 1000 & 1500 & 2000 & 2500 & 3000 \\ \hline
\multicolumn{1}{c|}{\multirow{6}{*}{$t^2$}} & \multicolumn{1}{c|}{\multirow{3}{*}{TL-gPINNs}} & max     & 34.47\%     & 25.59\%     &  36.57\%    & 28.42\%     & 47.01\%  & 52.87\%   \\ \cline{3-3}
\multicolumn{1}{c|}{}                   & \multicolumn{1}{c|}{}                           & min     & 16.71\%    & 15.16\%     & 20.03\%     & 13.99\%     & 16.29\% & 17.79\%   \\ \cline{3-3}
\multicolumn{1}{c|}{}                   & \multicolumn{1}{c|}{}                           & average &  22.97\%   &18.89\%      & 24.22\%     & 21.69\%     & 25.64\% & 26.42\%   \\ \cline{2-3}
\multicolumn{1}{c|}{}                   & \multicolumn{1}{c|}{\multirow{3}{*}{gPINNs}}    & max     & -36.98\%    & 8.64\%     & -14.95\%     & -11.62\%     & -9.35\%  & -26.32\%  \\ \cline{3-3}
\multicolumn{1}{c|}{}                   & \multicolumn{1}{c|}{}                           & min     & -93.22\%    & -213.93\%     & -134.70\%     & -143.01\%     & -223.35\% & -116.03\%   \\ \cline{3-3}
\multicolumn{1}{c|}{}                   & \multicolumn{1}{c|}{}                           & average &-51.83\%     & -54.32\%     & -51.77\%     & -36.28\%     & -74.26\% & -56.37\%   \\ \cline{1-3}
\multicolumn{1}{c|}{\multirow{6}{*}{$\sin(t)$}} & \multicolumn{1}{c|}{\multirow{3}{*}{TL-gPINNs}} & max     & 66.98\%    &  50.02\%    & 59.88\%     &  72.04\%    & 58.02\%  & 57.77\%  \\ \cline{3-3}
\multicolumn{1}{c|}{}                   & \multicolumn{1}{c|}{}                           & min     & 39.22\%    & 17.20\%     & 23.26\%     & 13.39\%     &20.81\%  & 21.79\%   \\ \cline{3-3}
\multicolumn{1}{c|}{}                   & \multicolumn{1}{c|}{}                           & average &54.45\%     & 37.96\%     &43.58\%      & 39.93\%     & 31.93\% & 43.59\%   \\ \cline{2-3}
\multicolumn{1}{c|}{}                   & \multicolumn{1}{c|}{\multirow{3}{*}{gPINNs}}    & max     & 36.15\%    & -8.73\%     &6.76\%      & 49.58\%     & -19.21\%  & 3.96\%  \\ \cline{3-3}
\multicolumn{1}{c|}{}                   & \multicolumn{1}{c|}{}                           & min     & -81.44\%    & -165.79\%     & -144.16\%     & -100.08\%     & -128.10\% & -131.65\%   \\ \cline{3-3}
\multicolumn{1}{c|}{}                   & \multicolumn{1}{c|}{}                           & average & -24.74\%    & -90.53\%     & -38.65\%     & -23.23\%     &-72.83\%  & -41.60\%   \\ \cline{1-3}
\multicolumn{1}{c|}{\multirow{6}{*}{$\tanh(t)$}} & \multicolumn{1}{c|}{\multirow{3}{*}{TL-gPINNs}} & max     & 73.85\%    & 68.23\%     & 79.10\%     & 75.23\%     & 64.39\%  &80.04\%   \\ \cline{3-3}
\multicolumn{1}{c|}{}                   & \multicolumn{1}{c|}{}                           & min     & 38.30\%    & 42.04\%     & 45.63\%     & 48.64\%     & 43.74\%  &43.71\%   \\ \cline{3-3}
\multicolumn{1}{c|}{}                   & \multicolumn{1}{c|}{}                           & average & 54.73\%    & 57.05\%     & 70.28\%     & 58.90\%     & 54.44\% &63.26\%    \\ \cline{2-3}
\multicolumn{1}{c|}{}                   & \multicolumn{1}{c|}{\multirow{3}{*}{gPINNs}}    & max     &59.33\%     & 56.81\%     & 56.50\%     &  58.63\%    & 55.75\% &52.71\%     \\ \cline{3-3}
\multicolumn{1}{c|}{}                   & \multicolumn{1}{c|}{}                           & min     &-39.43\%     & 27.21\%     & -46.44\%     & 21.92\%     & 18.13\% &21.78\%     \\ \cline{3-3}
\multicolumn{1}{c|}{}                   & \multicolumn{1}{c|}{}                           & average & 18.78\%    & 39.68\%     & 29.92\%     &41.50\%      & 42.50\% &35.28\%    \\ \cline{1-3}
\multicolumn{1}{c|}{\multirow{6}{*}{$\frac{1}{1+t^2}$}} & \multicolumn{1}{c|}{\multirow{3}{*}{TL-gPINNs}} & max     & 74.07\%    & 58.00\%     & 57.67\%     &60.42\%      & 74.94\% &73.77\%     \\ \cline{3-3}
\multicolumn{1}{c|}{}                   & \multicolumn{1}{c|}{}                           & min     & 20.26\%    & 33.87\%     & 26.32\%     & 51.01\%     & 24.59\% &44.46\%    \\ \cline{3-3}
\multicolumn{1}{c|}{}                   & \multicolumn{1}{c|}{}                           & average & 47.84\%    & 48.16\%     & 41.39\%     & 55.96\%     & 54.59\%  &59.05\%    \\ \cline{2-3}
\multicolumn{1}{c|}{}                   & \multicolumn{1}{c|}{\multirow{3}{*}{gPINNs}}    & max     & 31.69\%    & 37.35\%     & 51.14\%     &  42.45\%    & 56.68\%  &45.66\%   \\ \cline{3-3}
\multicolumn{1}{c|}{}                   & \multicolumn{1}{c|}{}                           & min     & -43.93\%    & -62.20\%     & -108.66\%     & -4.36\%     &-2.33\% &11.89\%     \\ \cline{3-3}
\multicolumn{1}{c|}{}                   & \multicolumn{1}{c|}{}                           & average & 8.60\%    &7.97\%      & 13.18\%     &  23.81\%    & 34.44\% & 31.54\%   \\ \toprule
\end{tabular}
\end{table}

\end{document}